\documentclass[11pt]{article}

\usepackage{geometry}   
\usepackage{authblk}

\usepackage[utf8]{inputenc} 
\usepackage[hidelinks]{hyperref}       
\usepackage{url}            
\usepackage{booktabs}       
\usepackage{amsfonts}       
\usepackage{nicefrac}       
\usepackage{microtype}      

\usepackage[T1]{fontenc}
\let\oldamp\&
\renewcommand{\&}{\normalfont\oldamp}

\usepackage{natbib}

\usepackage{amsmath,amsfonts,amssymb,amsthm,array}

\usepackage{algorithm}
\usepackage{caption}
\usepackage{subcaption}
\usepackage{algpseudocode}
\usepackage{mdframed} 
\usepackage{thmtools}

\usepackage{booktabs} 
\usepackage{array}
\usepackage{threeparttable}
\usepackage{makecell} 
\usepackage{booktabs,multirow,makecell}

\newcommand{\eqdef}{\; { := }\;}

\newcommand{\prox}{\operatorname{Prox}}
\newcommand{\R}{\mathbb{R}}

\newcommand{\cD}{\mathcal{D}}
\newcommand{\cX}{\mathcal{X}}
\newcommand{\cT}{\mathcal{T}}
\newcommand{\cK}{\mathcal{K}}

\newcommand{\cB}{\mathcal{B}}   

\newcommand{\Proj}{\Pi}

\newcommand{\ind}{\iota}
\newcommand{\dom}{\operatorname{dom}}  

\DeclareMathOperator*{\argmin}{argmin}

\usepackage[colorinlistoftodos,bordercolor=orange,backgroundcolor=orange!20,linecolor=orange,textsize=scriptsize]{todonotes}

\declaretheorem[style=shaded,within=section]{definition}
\declaretheorem[style=shaded,sibling=definition]{theorem}

\declaretheorem[style=shaded,sibling=definition]{assumption}
\declaretheorem[style=shaded,sibling=definition]{corollary}

\declaretheorem[style=shaded,sibling=definition]{lemma}
\declaretheorem[style=shaded,sibling=definition]{example}
\declaretheorem[style=shaded,sibling=definition]{remark}

\title{\bf Broximal Gradient Descent: A Generalization Beyond the Constraints}

\date{}

\author[1]{Xun Qian\thanks{email: \texttt{xun.qian@kaust.edu.sa}}}
\author[1]{Peter Richt\'arik\thanks{email: \texttt{peter.richtarik@kaust.edu.sa}}}

\affil[1]{King Abdullah University of Science and Technology, Thuwal, Saudi Arabia}

\begin{document}
	
	\maketitle
	
\begin{abstract}
We consider composite optimization problems of the form $\min_{x} F(x) := f(x) + g(x)$, where $f$ is differentiable and $g$ is a proper, closed, convex function that may encode a nonsmooth regularizer, a constraint, or both. We propose \emph{Broximal Gradient Descent} (BGD), which minimizes the first-order model of $f$ plus $g$ over a Euclidean ball around the current iterate, so that the radius of the ball plays the role of the stepsize. BGD extends the recently proposed BroxGD from constrained to composite problems, and contains proximal gradient descent with any stepsize as a special case, but not conversely. Our analysis rests on two master theorems: one describing a single step of the method under a condition that requires neither convexity nor smoothness of $f$, and one characterizing the radii for which this condition holds. From these we derive convergence rates: $\mathcal{O}(1/K)$ for convex $L$-smooth $f$, $\mathcal{O}(1/\sqrt{K})$ for convex $F$ with bounded subgradients, linear rates for strongly convex $f$ and under a proximal Polyak--{\L}ojasiewicz condition, and an $\mathcal{O}(1/K)$ rate for the squared gradient-mapping norm in the nonconvex smooth case, each with an idealized radius rule and an implementable one that does not require the solution. We also propose a rule that halves the radius whenever the objective has not decreased enough; it uses objective values only and attains the rates of the smooth cases without any problem constant. Experiments on LASSO, elastic net, and a box-constrained neural network show that BGD with an adaptive or a halving radius rule needs fewer matrix-vector products or network passes than proximal gradient descent with line search on most instances, at the price of a scalar search with about a dozen proximal evaluations per iteration.
\end{abstract}

\section{Introduction}

In this work we consider the composite optimization problem
\begin{equation}\label{p:p}
	\min_{x \in \mathbb{R}^d } F(x) \eqdef  f(x) + g(x), 
\end{equation}
where \( f : \mathbb{R}^d \to \mathbb{R} \) is a differentiable loss/cost/objective function, and $g : \mathbb{R}^d \to \mathbb{R} \cup \{ +\infty \}$ is proper, closed, and convex. The function $g$ may be a nonsmooth regularizer; the indicator function $\ind_{\cX}$ of a nonempty, closed, convex set $\cX$, in which case \eqref{p:p} is the constrained problem $\min_{x \in \cX} f(x)$; or a combination of the two. We denote the set of all minimizers of \eqref{p:p} by \( \mathcal{X}_\star \), and the infimum of $F$ by $F_\star$. 

Problems of the form \eqref{p:p} are ubiquitous in statistics, signal processing, and machine learning. Sparse regression models such as the LASSO \citep{tibshirani1996regression} and the elastic net \citep{zou2005regularization} take $f$ to be a least-squares loss and $g$ an $\ell_1$ or mixed $\ell_1$--$\ell_2$ regularizer; compressed sensing \citep{candes2006robust} and total-variation image denoising \citep{rudin1992nonlinear} lead to problems with the same structure, and low-rank matrix recovery replaces the $\ell_1$ norm by the nuclear norm \citep{recht2010guaranteed}. In regularized empirical risk minimization, $f$ is the average loss over the training data and $g$ a regularizer promoting sparsity or group structure.

\paragraph{Proximal gradient descent.} The workhorse first-order method for \eqref{p:p} is proximal gradient descent (PGD), also known as forward--backward splitting \citep{combettes2005signal, nesterov2013gradient, beck2009fast, parikh2014proximal}:
\begin{equation}\label{eq:pgd}
	x_{k+1} = \operatorname{Prox}_{\gamma g}\left( x_k - \gamma \nabla f(x_k) \right) = \argmin_{z \in \R^d} \left\{ \langle \nabla f(x_k), z \rangle + g(z) + \frac{1}{2\gamma}\|z - x_k\|^2 \right\},
\end{equation}
where $\operatorname{Prox}_{\gamma g}(x) \eqdef \argmin_{z} \{ g(z) + \frac{1}{2\gamma}\|z - x\|^2 \}$ is the proximal operator of $\gamma g$ \citep{moreau1965proximite}. PGD minimizes the first-order model $m_k(z) \eqdef \langle \nabla f(x_k), z\rangle + g(z)$ of $F$ at $x_k$, and keeps the minimizer close to $x_k$ by means of a quadratic penalty. The stepsize $\gamma$ has to be matched to the smoothness constant $L$ of $f$ (typically $\gamma \leq 1/L$), and the distance traveled in one step, $\|x_{k+1} - x_k\| = \gamma \|G_\gamma(x_k)\|$, where $G_\gamma$ is the gradient mapping (see Section~\ref{sec:prelim}), depends on the magnitude of the gradient in a way that is difficult to control a priori. When $L$ is unknown, or varies substantially over the domain, one resorts to a line search \citep{nesterov2013gradient, beck2009fast}, at the cost of extra function evaluations per iteration.

\paragraph{From penalties to balls.} An alternative to \emph{penalizing} the distance to the current iterate is to \emph{constrain} it. This idea underlies trust-region methods \citep{conn2000trust} and has recently resurfaced in optimizers for deep learning: the Muon optimizer \citep{jordan2024muon} and its relatives can be understood as steps of a linear minimization oracle (LMO) over a norm ball around the current iterate \citep{bernstein2024old, pethick2025training, kovalev2025understanding, riabinin2025gluon}. The nonconvex algorithm of \citet{kovalev2025understanding} is close to the method studied here, with two differences: it uses a general, typically non-Euclidean norm matched to the geometry of the parameters, whereas we use the Euclidean norm, which makes the Pythagorean distance argument at the core of our analysis available; and its radius is constant, whereas the choice of the radius sequence $t_k$ is the main object of our analysis. \citet{gruntkowska2025ball} studied the ball-proximal (``broximal'') point method, $x_{k+1} = \argmin_{z \in \cB(x_k, t_k)} f(z)$, which minimizes $f$ itself over the Euclidean ball $\cB(x_k, t_k) \eqdef \{ z \in \R^d \,|\, \|z - x_k\| \leq t_k \}$, and showed that the radius $t_k$ acts as a stepsize with remarkably clean convergence properties. Most closely related to the present work, \citet{richtarik2026local} proposed \emph{Broximal Gradient Descent} (BroxGD) for the constrained problem $\min_{x \in \cX} f(x)$, which minimizes the linear model of $f$ over $\cX \cap \cB(x_k, t_k)$, and developed a ``master theorem'' describing a single step of the method under two admissibility conditions on the radius.

\paragraph{Broximal gradient descent.} This paper extends BroxGD from constraints to general composite problems of the form \eqref{p:p}. We propose \emph{Broximal Gradient Descent} (Broximal GD, or BGD for short), which at each iteration computes
\begin{equation}\label{eq:BG}
	x_{k+1} \in \argmin_{z \in \cB(x_k, t_k)} m_k(z), \qquad m_k(z) \eqdef \langle \nabla f(x_k), z \rangle + g(z).
\end{equation}
Compared with PGD \eqref{eq:pgd}, the model $m_k$ is the same, but the quadratic penalty is replaced by a ball constraint: the amount of trust placed in the model is expressed through the radius $t_k$, which directly controls how far the iterate moves, rather than through a penalty parameter. BGD contains several methods as special cases. For $g = \ind_{\cX}$ it is BroxGD \citep{richtarik2026local}, whose name we keep for the composite method; for $g \equiv 0$ it is normalized gradient descent \citep{nesterov1984minimization, nesterov2018lectures, hazan2015beyond}, $x_{k+1} = x_k - t_k \nabla f(x_k)/\|\nabla f(x_k)\|$ (when $\nabla f(x_k) \neq 0$; if $\nabla f(x_k) = 0$, the model is constant and every point of the ball is a minimizer); and, as we show in Section~\ref{sec:BGrecover}, for every stepsize $\gamma > 0$ there is a radius for which the broximal step coincides with the PGD step \eqref{eq:pgd}, so that PGD itself is a special case of BGD. The converse is false: a broximal step need not be a PGD step for any stepsize (Remark~\ref{rem:not-prox}), and for large radii the minimizer in \eqref{eq:BG} need not even be unique (Remark~\ref{rem:strict}). BGD is therefore a strict generalization of PGD, not a reparameterization of it.

\paragraph{Related work.} The subproblem \eqref{eq:BG} is not new; what is new here is its analysis for a general convex $g$ with radii that are fixed in advance (constant, decaying, or geometric), or adjusted from the observed decrease of the objective without ever rejecting a step. For $g = \ind_\cX$, \citet{ferris1996linear} introduced the minimization of the linearized objective over $\cX$ intersected with a norm ball, showed that for the Euclidean norm a projected gradient step with an implicit stepsize solves it, selected the radius by an Armijo-type search, and proved asymptotic linear convergence under an error bound; \citet{garber2016linearly} used a local linear oracle over polytopes. Trust-region methods for nonsmooth composite problems \citep{baraldi2023proximal, aravkin2022proximal, dao2025projected} minimize a quadratic model plus $g$ over a ball, which coincides with \eqref{eq:BG} when the quadratic term vanishes and the subproblem is solved exactly; they accept steps by a ratio test, which requires a function evaluation per trial, prove $\mathcal{O}(\epsilon^{-2})$ complexity bounds for reaching stationarity, without function-value rates in the convex case, and, as we do, count proximal evaluations separately from gradient evaluations. BroxGD \citep{richtarik2026local} analyzes \eqref{eq:BG} for $g = \ind_\cX$ under either of two admissibility conditions on the radius, which are stated a priori and involve $x_\star$, with prescribed, Polyak, and backtracking radii; \citet{richtarik2026secretly} showed that, for radii not exceeding the Polyak radius, its step is the projection of $x_k$ onto $\cX$ intersected with a half-space that separates $x_k$ from the solution set. Table~\ref{tab:related} summarizes the comparison. Frank--Wolfe methods \citep{frank1956algorithm, jaggi2013revisiting, lacoste2015global} also rely on a linear minimization oracle, but over the entire feasible set, and then take a convex combination with the current iterate; BGD instead localizes the oracle to a ball, which is what allows the radius to serve as a stepsize. The idealized radius rules in our analysis, which depend on $x_\star$ or $F_\star$, play the same role as the Polyak stepsize for (sub)gradient methods \citep{polyak1969minimization, hazan2019revisiting, loizou2021stochastic}, and the proximal P{\L} condition we use in the nonconvex case is the composite analogue of the classical Polyak--{\L}ojasiewicz inequality \citep{polyak1963gradient, lojasiewicz1963topological, karimi2016linear}. Our analysis extends the master-theorem approach of \citet{richtarik2026local} from the linear model over $\cX \cap \cB(x_k,t_k)$ to the composite model $m_k$ over $\cB(x_k, t_k)$. For more related work, see \citet{richtarik2026local} and the references therein. 

\begin{table}[t!]
	\centering
	\begin{threeparttable}
	\caption{Methods based on a ball-constrained linear or quadratic model. ``Prescribed'' radii are fixed in advance; the adaptive and halving radii of this paper are adjusted from observed objective values, and no step is ever rejected.}
	\label{tab:related}
	\footnotesize
	\begin{tabular}{@{}>{\raggedright\arraybackslash}p{0.14\linewidth}>{\raggedright\arraybackslash}p{0.11\linewidth}>{\raggedright\arraybackslash}p{0.12\linewidth}>{\raggedright\arraybackslash}p{0.23\linewidth}>{\raggedright\arraybackslash}p{0.27\linewidth}@{}}
		\toprule
		Work & $g$ & Radius & One-step results & Rates \\
		\midrule
		\citet{ferris1996linear} & $\ind_\cX$ & Armijo-type search & a projected gradient step with implicit stepsize solves the subproblem & asymptotic; linear under an error bound \\
		\midrule
		\citet{aravkin2022proximal} & proper lsc, possibly nonconvex & ratio test & sufficient model decrease (Cauchy-type), checked by a ratio test & $\mathcal{O}(\epsilon^{-2})$ to stationarity \\
		\midrule
		\citet{dao2025projected} & proper closed convex & ratio test & sufficient model decrease (Cauchy-type), checked by a ratio test & $\mathcal{O}(\epsilon^{-2})$ to stationarity \\
		\midrule
		\citet{richtarik2026local} & $\ind_\cX$ & prescribed, Polyak, backtracking & full radius and $-t_k^2$ decrease under a priori conditions involving $x_\star$ & $\mathcal{O}(1/\sqrt{K})$, $\mathcal{O}(1/K)$, linear (strongly convex)\tnote{a}; nonconvex and P{\L} with projection-dependent radii\tnote{b} \\
		\midrule
		\citet{richtarik2026secretly} & $\ind_\cX$ & Polyak & step is the projection onto $\cX\cap\{\text{half-space}\}$ for radii up to the Polyak radius & universal H\"older rates (best iterate) \\
		\midrule
		This paper & proper closed convex & prescribed, adaptive, halving & full radius and $-t_k^2$ decrease under an a posteriori model-value condition\tnote{c} & $\mathcal{O}(1/K)$, $\mathcal{O}(1/\sqrt{K})$, linear (strongly convex, last iterate; proximal P{\L}), $\mathcal{O}(1/K)$ for $\|G_\gamma\|^2$, all with prescribed radii; the same rates, except $\mathcal{O}(1/\sqrt{K})$, with a halving rule that needs no problem constant \\
		\bottomrule
	\end{tabular}
	\begin{tablenotes}
		\item[a] With prescribed geometric radii, \citet[Theorem~3.3]{richtarik2026local} contract the objective gap and the squared distance by the factor $1-\frac{\mu}{4L}$ per iteration, i.e., $f(x_K)-f_\star\leq\epsilon$ after $\frac{4L}{\mu}\ln\frac{C_0}{\epsilon}$ iterations. Theorem~\ref{th:strongly-convex-prac} below contracts the squared distance by $q^2 = 1-\frac{\theta^2}{1+2\theta+3\theta^2}$ with $\theta = \frac{2\sqrt{\mu L}}{L+\mu}$, which is about $1-\frac{4\mu}{L}$ for $L\gg\mu$, i.e., $\|x_K-x_\star\|\leq\epsilon$ after about $\frac{L}{2\mu}\ln\frac{\|x_0-x_\star\|}{\epsilon}$ iterations. The idealized rule of Theorem~\ref{th:strongly-convex-theory} gives the contraction $(\frac{L-\mu}{L+\mu})^{2K}$ of \citet[Theorem~3.7]{richtarik2026local}, and their backtracking rule gives $(1+\xi^2\frac{\mu}{L})^{-K}$ with a factor $\xi\in(0,1)$ \citep[Theorem~3.6]{richtarik2026local}. With the halving rule, which requires neither $\mu$ nor $L$, Theorem~\ref{th:sc-halving} divides the objective gap by four every $\frac{4L}{\mu}+T+4$ iterations (for $c = \nicefrac{L}{2}$).
		\item[b] Under a projected-P{\L} condition, \citet[Theorem~D.16]{richtarik2026local} obtain the rate $(1-\frac{\mu}{L})^K$ with the radius $\frac{1}{L}\|G_{\nicefrac{1}{L}}(x_k)\|$, which requires a projection at every iteration. Theorem~\ref{th:non-convex-PL-theory} recovers this rate with the same radius under the weaker Assumption~\ref{as:proximalPL}, Theorem~\ref{th:non-convex-PL} obtains, with a prescribed geometric radius, a complexity within a factor of $\nicefrac{4}{3}$ of it, and Theorem~\ref{th:pl-halving} obtains a linear rate with the halving rule, without knowledge of $\mu$ or $L$.
		\item[c] The admissible radii form an interval $(0,b_k)$, and the step with radius $b_k$ is the projection of $x_k$ onto the model sublevel set (Theorem~\ref{th:Tk}); for $g = \ind_\cX$, this is the projection of \citet{richtarik2026secretly}.
	\end{tablenotes}
	\end{threeparttable}
\end{table}

\paragraph{Contributions.} Our main contributions are as follows.
\begin{itemize}
	\item \textbf{Master theorems for composite models}. Theorem~\ref{th:another-ad} describes a single broximal step relative to an arbitrary reference point under an a posteriori condition on the model value, and requires neither convexity nor smoothness of $f$; for $g = \ind_\cX$ and $y = x_\star$ its hypothesis is implied by the admissibility conditions of \citet{richtarik2026local}, whose master theorem it recovers. Theorem~\ref{th:Tk} characterizes the admissible radii as an interval $(0, b_k)$ and identifies its endpoint as the projection of $x_k$ onto the model sublevel set, which for $g=\ind_\cX$ is the projection interpretation of \citet{richtarik2026secretly}. Theorem~\ref{th:recoverPGD} shows that every PGD step is a broximal step whose radius is its displacement, the composite version of a result of \citet{ferris1996linear}, but not conversely: a broximal step need not be a PGD step for any stepsize, nor need it be unique.
	\item \textbf{Convergence rates}. From the master theorems we derive an $\mathcal{O}(1/K)$ rate for convex $L$-smooth $f$, an $\mathcal{O}(1/\sqrt{K})$ rate for convex $F$ with bounded subgradients, linear rates for strongly convex $f$ and under a proximal P{\L} condition, and an $\mathcal{O}(1/K)$ rate for the squared gradient-mapping norm for nonconvex smooth $f$. The rules in which the radius is the displacement of a PGD step (Theorems~\ref{th:convex-f} and \ref{th:non-convex-PL-theory}) recover the rates of PGD and serve as baselines; the other rules use radii that are fixed in advance (constant, decaying, or geometric) or adapted, and do not require $x_\star$. In the strongly convex case, a geometric envelope argument yields a last-iterate bound, and we show how knowledge of the condition number can be removed by restarts; under the proximal P{\L} condition, the linear rate holds with a geometric radius, whereas \citet{richtarik2026local} require a projection-dependent radius. Finally, the halving rule \eqref{eq:halving}, which uses objective values only and no problem constant, attains, for the last iterate, an $\mathcal{O}(1/K)$ rate in the convex case (Theorem~\ref{th:convex-halving}) and a linear rate with the complexity $\mathcal{O}(\frac{L}{\mu}\ln\frac{1}{\epsilon})$ in the strongly convex case (Theorem~\ref{th:sc-halving}) and under the proximal P{\L} condition (Theorem~\ref{th:pl-halving}), as well as the $\mathcal{O}(1/K)$ rate for the smallest squared gradient-mapping norm in the nonconvex case (Theorem~\ref{th:nc-halving}). 
	\item \textbf{Experiments}. On LASSO and elastic net ($20$ random instances each) and a box-constrained neural network trained on MNIST ($10$ random seeds), BGD with the adaptive rule \eqref{eq:adaptive} or the halving rule \eqref{eq:halving} is competitive with PGD with line search and clearly outperforms PGD with a constant stepsize: on the two regression problems, both rules are faster than the line search on most instances in terms of matrix-vector products, the halving rule being the fastest in median, and on the network both reach a training loss an order of magnitude smaller than the line search, with similar values for the two rules; each BGD iteration requires about a dozen proximal evaluations.
\end{itemize}

\section{Notation, Assumptions, and the Algorithm}\label{sec:prelim}

\paragraph{Notation.} Throughout, $\|\cdot\|$ is the Euclidean norm and $\langle \cdot, \cdot \rangle$ the standard inner product on $\R^d$. For $x \in \R^d$ and $t \geq 0$ we write $\cB(x, t) \eqdef \{ z \in \R^d \,|\, \|z - x\| \leq t \}$ for the closed Euclidean ball of radius $t$ centered at $x$. The effective domain of $g$ is $\dom g \eqdef \{ x \in \R^d \,|\, g(x) < +\infty \}$, and $\partial g(x)$ is the subdifferential of $g$ at $x$. With a slight abuse of notation that is standard in this context, we also write $\partial g(x)$ for a particular element of the subdifferential whenever the choice is clear from the context (for instance, the subgradient appearing in the optimality conditions of a proximal subproblem). The proximal operator of $\gamma g$, $\gamma > 0$, is
$$
\operatorname{Prox}_{\gamma g}(x) \eqdef \argmin_{z \in \R^d} \left\{ g(z) + \frac{1}{2\gamma}\|z - x\|^2 \right\},
$$
and the associated gradient mapping is
$
G_\gamma(x) \eqdef \nicefrac{1}{\gamma}\left( x - \operatorname{Prox}_{\gamma g}\left( x - \gamma \nabla f(x) \right) \right).
$
We say that $f$ is $L$-smooth on a set $S \subseteq \R^d$ if $\|\nabla f(x) - \nabla f(y)\| \leq L \|x - y\|$ for all $x, y \in S$, and $\mu$-strongly convex on $S$ if $f(y) \geq f(x) + \langle \nabla f(x), y - x \rangle + \frac{\mu}{2}\|y - x\|^2$ for all $x, y \in S$. We write $F_\star \eqdef \inf_x F(x)$ and $x_\star$ for an arbitrary element of $\cX_\star$.

\paragraph{Standing assumption.} The following assumption is in force throughout the paper. 

\begin{assumption}\label{as:basic}
(i) $f : \R^d \to \R$ is differentiable and not necessarily convex;
(ii) $g : \R^d \to \R \cup \{+\infty\}$ is proper, closed, and convex;
(iii) whenever $f$ is assumed to be convex, the set of minimizers $\cX_\star$ of problem \eqref{p:p} is nonempty; 
(iv) the infimum $F_\star \eqdef \inf_x F(x)$ is finite. 
\end{assumption}

\paragraph{The algorithm.} BGD is summarized in Algorithm~\ref{alg:BG}. Given $x_k$ and a radius $t_k \geq 0$, the next iterate minimizes the first-order model $m_k(z) = \langle \nabla f(x_k), z \rangle + g(z)$ over the ball $\cB(x_k, t_k)$.

\begin{algorithm}[t]
	\caption{Broximal Gradient Descent (BGD)}\label{alg:BG}
	\begin{algorithmic}[1]
		\Require initial point $x_0 \in \dom g$; radii $t_k \geq 0$, $k = 0, 1, 2, \dots$
		\For{$k = 0, 1, 2, \dots$}
		\State Compute $\nabla f(x_k)$
		\State $x_{k+1} \in \argmin_{z \in \cB(x_k, t_k)} \left\{ \langle \nabla f(x_k), z \rangle + g(z) \right\}$
		\EndFor
	\end{algorithmic}
\end{algorithm}

The subproblem in Algorithm~\ref{alg:BG} is always well posed. Indeed, since $g$ is proper, closed, and convex, so is $m_k$, and the ball $\cB(x_k, t_k)$ is compact; as we assume $x_0 \in \dom  g$, the feasible set $\cB(x_0, t_0) \cap \dom g$ is nonempty and the subproblem has a solution $x_1 \in \dom g$. By induction, the same holds at every iteration, and hence for any $k\geq 0$, we have $x_k \in \dom  g$.

\section{Master theorems}

The following theorem describes a single step of BGD relative to an arbitrary reference point $y \in \dom g$. It is the basic tool of the paper: in the convex case we will apply it with $y = x_\star$, and in the nonconvex case with $y$ equal to a proximal gradient point. 

\begin{theorem}\label{th:another-ad}
	Let Assumption \ref{as:basic} hold. Fix \( x_k \in \dom  g \), \( y \in  \dom  g \), and let $t_k>0$. Let \( x_{k+1} \) be obtained from \( x_k \) after one step of Broximal GD with radius \( t_k \), i.e.,
	\[
	x_{k+1} \in \argmin_{z \in  \mathcal{B}(x_k, t_k)} m_k(z) \eqdef \langle \nabla f(x_k), z \rangle + g(z).
	\]
	Assume $m_k(x_{k+1}) > m_k(y)$.
	Then
	\begin{itemize}
		\item[(i)] \( y \) does not belong to the ball \( \mathcal{B}(x_k, t_k) \), i.e.,
		$
		t_k < \|x_k - y \|;
		$
		\item[(ii)] \( x_{k+1} \) lies on the boundary of the ball \( \mathcal{B}(x_k, t_k) \), i.e., $\|x_{k+1} - x_k\| = t_k$; 
		
		\item[(iii)] the squared distance to \( y \) decreases by at least the square of the radius, i.e.,
		\begin{equation}\label{eq:distance-tk2}
		\|x_{k+1} - y\|^2 < \|x_k - y\|^2 - t_k^2; 
		\end{equation}
		
		\item[(iv)] $x_{k+1}$ is the unique solution of $\min_{z \in \mathcal{B}(x_k, t_k)} m_k(z)$. 
		
		\item[(v)] $x_{k+1} = \operatorname{Prox}_{\gamma_k g} \left(  x_k - \gamma_k \nabla f(x_k)  \right)$ for $\gamma_k = \frac{t_k}{\|\nabla f(x_k) + \partial g(x_{k+1})\|}>0$, where $\partial g(x_{k+1})$ denotes a suitable subgradient of $g$ at $x_{k+1}$. 
	\end{itemize}
\end{theorem}

The hypothesis $m_k(x_{k+1}) > m_k(y)$ means that the ball is too small to reach the model value of $y$; the step then exhausts the full radius and, by (ii) and (iii), lands strictly inside the ball with diameter $[x_k, y]$, with no smoothness, convexity, or Lipschitz property of $f$ involved. Under the hypothesis, the broximal step is a proximal gradient step whose stepsize is determined a posteriori by the radius (v); once $m_k(x_{k+1}) \leq m_k(y)$, it may fail to be unique (Remark~\ref{rem:strict}) and may fail to be a proximal gradient step for any stepsize (Remark~\ref{rem:not-prox}), so BGD is not a special case of PGD, while the reverse inclusion holds (Section~\ref{sec:BGrecover}).

\begin{example}[Comparison with the admissibility conditions of \citet{richtarik2026local}]\label{ex:admissibility}
	Consider the constrained case $g = \ind_{\cX}$, in which BGD reduces to BroxGD \citep{richtarik2026local}. The following example shows that the upper bounds on $t_k$ imposed by the Type-I and Type-II admissibility conditions of \citet{richtarik2026local} can be much smaller than the radii for which the hypothesis of Theorem~\ref{th:another-ad} holds, i.e., for which $\langle \nabla f(x_k), x_{k+1} - x_\star \rangle > 0$. Let $\mathcal{X} \eqdef \{  (u, v) | u=0, \ 0\leq v\leq 1  \}$ and $f(u, v) = uv + (v-1)^4$. The optimal solution on $\mathcal{X}$ is $x_\star = (0, 1)$. For $x_k = (0, 1-\beta)$ with $0<\beta <1$, $\nabla f(x_k) = (1-\beta,  - 4\beta^3)$. The upper bounds on $t_k$ in the Type-I and Type-II admissibility conditions are $t_{k, 1} = \frac{4\beta^4}{\sqrt{(1-\beta)^2 + 16\beta^6}}$ and $t_{k, 2} = \frac{4\beta^3}{\sqrt{1+16\beta^4}}$, respectively. In contrast, for every $t_k \in (0, \beta)$ one has $x_{k+1} = (0, 1-\beta +t_k)$ and $\langle \nabla f(x_k), x_{k+1} - x_\star \rangle = 4\beta^3(\beta  - t_k) > 0$, so that the hypothesis of Theorem~\ref{th:another-ad} holds. Choosing $t_k = \frac{\beta}{2}$, we obtain $\lim\limits_{\beta \to 0^+} \frac{t_{k, 1}}{t_k} = 0$ and $\lim\limits_{\beta \to 0^+} \frac{t_{k, 2}}{t_k} = 0$.
\end{example}

Theorem~\ref{th:another-ad} applies for a given radius $t_k$ as long as $m_k(x_{k+1}) > m_k(y)$. We now characterize the set of radii with this property. Define $x_{k+1}(t) \in \argmin_{z \in  \mathcal{B}(x_k, t)} m_k(z)$ and $\cT_k(y) \eqdef \{  t>0 | m_k(x_{k+1}(t)) - m_k(y) > 0  \}$ for any $y\in\dom g$. For simplicity, we denote $\cT_k \eqdef \cT_k(x_\star)$, where $x_\star \in {\cX}_\star$. The properties of $\cT_k(y)$ are collected in the following theorem.

\begin{theorem}\label{th:Tk}
Let Assumption \ref{as:basic} hold. Fix \( x_k \in \dom  g \), \( y\in \dom  g \). Then the following results hold. 
	\begin{itemize}
	\item[(i)] $\cT_k(y) \neq \emptyset$ if and only if $m_k(x_k) - m_k(y) > 0$. If ${\bar t} \in \cT_k(y)$, then $t \in \cT_k(y)$ for any $0<t<{\bar t}$. 
	
	\item[(ii)] Assume $\cT_k(y) \neq \emptyset$. Then $\cT_k(y)$ is an open interval $(0, b)$ for some constant $0< b \leq \|x_k - y\|$, and $m_k(x_{k+1}(b)) = m_k(y)$.  Moreover, $\|x_{k+1}(b) - x_k\| = b$, $x_{k+1}(b)$ is the unique solution of $\min_{z \in \mathcal{B}(x_k, b)} m_k(z)$, and 
	$$
		\|x_{k+1}(b) - y\|^2 \leq \|x_k - y\|^2 - b^2. 
	$$
	Furthermore, we have $\cT_k(x_{k+1}(b)) = \cT_k(y)=(0, b)$, and for any $0<t<b$,
	$$
	\|x_{k+1}(t) - x_{k+1}(b)\|^2 < b^2 - t^2. 
	$$
	
	\item[(iii)] Denote $\cT_k^\star \eqdef \{  t>0 |  \langle \nabla f(x_k) - \nabla f(x_\star), x_{k+1}(t) - x_\star \rangle >0  \}$, where $x_\star \in \cX_\star$. Then $\cT_k^\star \subseteq \cT_k$. 
\end{itemize}

\end{theorem}

Geometrically, $b$ is the distance from $x_k$ to the sublevel set $\{z \,|\, m_k(z) \leq m_k(y)\}$ and $x_{k+1}(b)$ is the projection of $x_k$ onto it, which is the basis of the proof: the admissible radii are those for which the ball misses the sublevel set, and at $t=b$ the distance decrease follows from the projection inequality. For $g = \ind_\cX$, the sublevel set is $\cX$ intersected with a half-space with normal $\nabla f(x_k)$, so the endpoint step is a projection onto a ``localized feasible set'' in the sense of \citet{richtarik2026secretly}, who showed that every step with a radius not exceeding the Polyak radius is such a projection. The last statement in (ii) says that the steps with smaller radii are nested toward $x_{k+1}(b)$, which can therefore replace $y$ as the reference point without changing the admissible radii.

\begin{remark}
Theorem \ref{th:another-ad} and claims (i) and (ii) of Theorem \ref{th:Tk} remain valid if $m_k(\cdot)$ is replaced by any proper, closed, and convex function $m$, provided that $x_k \in \dom m$ and $y \in \dom m$, so that $m(x_k)$ and $m(y)$ are finite; in that case, $\dom g$ in the statements is replaced by $\dom m$. The proofs are identical, since they only use differences of finite model values. When $g$ is the indicator function of a nonempty, closed, and convex set, Theorems \ref{th:another-ad} and \ref{th:Tk} with $y=x_\star$ recover the master theorem of BroxGD \citep[Theorem~3.1]{richtarik2026local}.
\end{remark}

\begin{remark}\label{rem:strict}
	The strict inequality in the definition of $\cT_k(y)$ cannot be relaxed. The following example shows that for $t \in \{  t>0 |  m_k(x_{k+1}(t)) \geq m_k(x_\star)  \}$, $x_{k+1}(t)$ may fail to satisfy the inequality $\|x_{k+1}(t) - x_\star\|^2 \leq \|x_k - x_\star\|^2 - t^2$. Let $g$ be the indicator function of $\cX \eqdef \{  (u, v) | u\geq -1  \}$ and $f(u, v) = \frac{1}{2}(u+v)^2 + u$. The optimal solution on $\cX$ is $x_\star = (-1, 1)$. For $x_k = (0, 0)$, $\nabla f(x_k) = (1, 0)$. Choosing $t=2$ and the minimizer $x_{k+1}(2) = (-1, -1)$, we have $ m_k(x_{k+1}(t)) - m_k(x_\star) = \langle \nabla f(x_k), x_{k+1}(t) - x_\star \rangle =0$. However, $\|x_{k+1}(2) - x_\star\|^2 = 4 > 2 = \|x_k - x_\star\|^2$. Note that here $\cT_k = (0, 1)$, whereas $\{ t > 0 \,|\, m_k(x_{k+1}(t)) \geq m_k(x_\star) \} = (0, +\infty)$, and that the minimizer $x_{k+1}(2)$ is not unique.
\end{remark}

\begin{remark}\label{rem:not-prox}
	It is natural to ask whether $x_{k+1}(b)$ equals $\operatorname{Prox}_{\gamma g} \left(  x_k - \gamma \nabla f(x_k)  \right)$ for some $\gamma>0$, as is the case for every $x_{k+1}(t)$ with $t < b$ by Theorem~\ref{th:another-ad}(v). We show that this does not hold in general by providing a counterexample. Let $g$ be the indicator function of $\cX \eqdef \{  (u,v) | u\geq v^2-2v  \}$ and $f(u,v) = (u+2)^2$. In this case, the operator $\operatorname{Prox}_{\gamma g} (\cdot)$ is the projection onto $\cX$. The optimal point on $\cX$ is $x_\star=(-1, 1)$. For $y=x_\star$ and $x_k=(0, 0)$, we have $\nabla f(x_k) = (4, 0)$, $\cT_k = (0, \sqrt{2})$, and $x_{k+1}(b) = x_\star$. Then $x_k - \gamma \nabla f(x_k) = (-4\gamma, 0)$. The set $\cX$ is convex, and since the constraint $v^2-2v-u\leq0$ is active at $x_\star$ with gradient $(-1,0)$, its normal cone is $\mathcal{N}_\cX(x_\star) = \{(-a,0) \,|\, a\geq0\}$. As $(-4\gamma,0) - x_\star = (1-4\gamma,-1)\notin\mathcal{N}_\cX(x_\star)$, the projection of $(-4\gamma,0)$ onto $\cX$ differs from $x_\star = x_{k+1}(b)$ for every $\gamma>0$. Hence the boundary radius $b$ is genuinely different from the radii $t < b$: here $x_{k+1}(b)$ is the exact solution $x_\star$, which no proximal gradient step from $x_k$ can reach. In particular, BGD cannot be obtained from PGD by a special choice of stepsize.
\end{remark}

\section{BGD recovers PGD with specific radii}\label{sec:BGrecover}

Let $y_k \eqdef \operatorname{Prox}_{\gamma g} \left( x_k - \gamma \nabla f(x_k) \right)$ for some $\gamma>0$. If $y_k = x_k$, then $x_k$ is a stationary point of problem \eqref{p:p}, since $y_k = x_k$ if and only if $0 \in \nabla f(x_k) + \partial g(x_k)$, and Broximal GD with $t_k = 0$ returns $x_{k+1} = y_k = x_k$. The following theorem treats the case $y_k \neq x_k$. 

\begin{theorem}\label{th:recoverPGD}
		Let Assumption \ref{as:basic} hold. Fix some $\gamma >0$ and let \( x_{k+1} \) be obtained from \( x_k \) after one step of Broximal GD with radius \( t_k>0 \), i.e.,
	\[
	x_{k+1} \in \argmin_{z \in \mathcal{B}(x_k, t_k)} m_k(z) = g(z) + \langle \nabla f(x_k), z \rangle.
	\]
	Denote $y_k \eqdef \operatorname{Prox}_{\gamma g} \left( x_k - \gamma \nabla f(x_k) \right)$. Assume $y_k \neq x_k$. Then the following results hold. 
	
	\begin{itemize}
		
		\item[(i)] $\cT_k(y_k) = (0, b_k)$ with $b_k = \gamma \|G_\gamma(x_k)\| = \|x_k-y_k\|>0$, and by choosing $t_k=b_k$, we have $x_{k+1} = y_k$. 
		
		\item[(ii)] For $y_k^{(\alpha)} = \operatorname{Prox}_{\alpha \gamma g} \left(  x_k - \alpha \gamma \nabla f(x_k)  \right)$ with any $0<\alpha<1$, we have $x_{k+1}(t^{(\alpha)}) = y_k^{(\alpha)}$, where $t^{(\alpha)} = \|y_k^{(\alpha)} - x_k\| \geq \alpha \|y_k-x_k\|$. Furthermore, if $t^{(\alpha)} = \alpha \|y_k-x_k\|$, then $y_k^{(\alpha)} = (1-\alpha) x_k + \alpha y_k$. 
		
		\item[(iii)] For $t = \alpha\|y_k-x_k\|$ with any $0<\alpha<1$, we have $x_{k+1}(t) = \operatorname{Prox}_{\gamma_t g}\left(  x_k - \gamma_t \nabla f(x_k) \right)$, where $0<\gamma_t \leq \alpha \gamma$. 
		
	\end{itemize}
	
\end{theorem}

Thus PGD with stepsize $\gamma$ is BGD with the radius $t_k = \gamma\|G_\gamma(x_k)\|$, the displacement of the proximal gradient step (for $g = \ind_\cX$, this is \citet[Theorem~2.10]{ferris1996linear}), and by (ii) and (iii) decreasing the radius below $b_k$ traces the proximal gradient path back toward $x_k$ with implicit stepsizes $\gamma_t \leq \alpha\gamma$. The inequality $t^{(\alpha)} \geq \alpha \|y_k - x_k\|$ says that the displacement of a proximal gradient step, which is nondecreasing in the stepsize (Lemma~\ref{lm:prox-mono}), grows at most linearly with it, with equality if and only if $y_k^{(\alpha)} = (1-\alpha)x_k + \alpha y_k$, as in the unconstrained case $g \equiv 0$. The relationship is one-directional: every proximal gradient step is a broximal step but not conversely, and parameterizing the step by the radius, which lives in the same units as $\|x_k - x_\star\|$, allows the analysis below to be carried out entirely in terms of distances.

\section{Convex case}

\subsection{$f$ is $L$-smooth}\label{sec:convex-smooth}

We consider the case where $f$ is $L$-smooth and convex, and begin with the idealized radius $t_k = b_k$, the supremum of $\cT_k = \cT_k(x_\star)$.

\begin{theorem}\label{th:convex-f}
	Let Assumption \ref{as:basic} hold, and assume that $f$ is $L$-smooth and convex on $\dom  g$. Let \(\{x_k\}_{k \geq 0}\) be the iterates generated by Broximal GD. Further, let $x_\star \in \cX_\star$, and assume that whenever $F(x_k) > F(x_\star)$, we choose the radius $t_k = b_k$, where $b_k$ is the supremum of the set $\cT_k \eqdef \{  t>0 |  m_k(x_{k+1}(t)) > m_k(x_\star)  \}$, and $t_k=0$ when $F(x_k) = F(x_\star)$. Then $\frac{\sqrt{2(F(x_{k+1}) - F(x_\star))}}{\sqrt{L}} \leq t_k = b_k \leq \|x_k - x_\star\|$ for all $k$ such that $F(x_k) > F(x_\star)$, and for every $K\geq 1$, 
	$$
	\frac{1}{K} \sum_{k=1}^K (F(x_k) - F(x_\star)) \leq \frac{L\|x_0 - x_\star\|^2}{2K}. 
	$$
	Moreover, the averaged iterate ${\hat x}_K \eqdef \frac{1}{K} \sum_{k=1}^K x_k$ satisfies 
	$$
	F({\hat x}_K) - F(x_\star)  \leq \frac{L\|x_0 - x_\star\|^2}{2K}. 
	$$
\end{theorem}

The rate $\mathcal{O}(L\|x_0 - x_\star\|^2/K)$ coincides with that of PGD with stepsize $1/L$; it follows from the lower bound $b_k \geq \sqrt{2(F(x_{k+1}) - F_\star)/L}$ combined with the decrease of $\|x_k - x_\star\|^2$ by $b_k^2$ given by Theorem~\ref{th:Tk}(ii). The rule is not implementable, since $b_k$ depends on $x_\star$; the next result replaces it by an implementable schedule.

\begin{remark}
The radius rule in Theorem \ref{th:convex-f} is equivalent to the following one: $t_k = b_k$ if $\cT_k \neq \emptyset$, and $t_k=0$ otherwise. In fact, when $f$ is convex, $\cT_k = \emptyset \iff m_k(x_k) = m_k(x_\star) \iff F(x_k) = F(x_\star)$. 
\end{remark}

\begin{theorem}\label{th:convex-f-prac}
Let Assumption \ref{as:basic} hold, and assume that $f$ is $L$-smooth and convex on $\dom  g$. Let \(\{x_k\}_{k \geq 0}\) be the iterates generated by Broximal GD. Further, let $x_\star \in \cX_\star$, and let the radii be $t_k = \frac{\sqrt{2}}{\sqrt{L}} a_{k+1}$, where $\{  a_k  \}_{k=0}^{+\infty}$ is any positive sequence such that $\sum_{k=0}^{+\infty} a_k^2 = +\infty$. Then 
\begin{itemize}
	\item[(i)] for any $k\geq 1$ with $k \notin \cK$, where $\cK \eqdef \{  k |  F(x_k) - F(x_\star) \leq a_k^2  \}$, we have $t_{k-1} \in \cT_{k-1}$ and 
	$$
	\|x_k - x_\star\|^2 < \|x_{k-1} - x_\star\|^2 - \frac{2}{L} a_k^2; 
	$$
	
	\item[(ii)] there exists a subsequence $\{  x_{j_k}  \}_{k=0}^{+\infty}$ with $0\leq j_0<j_1<\dots$, such that $F(x_{j_k}) - F(x_\star) \leq a_{j_k}^2$ for all $k\geq 0$. 
\end{itemize}

\end{theorem}

The proof rests on a dichotomy that recurs throughout the paper: at each iteration, either $F(x_k) - F_\star \leq a_k^2$, or the radius $t_{k-1}$ is admissible and the squared distance to $x_\star$ decreases by $2a_k^2/L$; a smaller $a_k$ thus buys a more accurate subsequence at the price of a slower decrease of the distance.

For instance, we can choose the sequence $\{  a_k  \}_{k=0}^{+\infty}$ in Theorem \ref{th:convex-f-prac} as $a_k = \sqrt{\frac{L}{(k+2) \ln (k+2)}}$, for which $a_k \to 0$ while $\sum_k a_k^2 = +\infty$. By setting $a_k$ to a constant, we obtain the following corollary. 

\begin{corollary}\label{co:convex-f-prac}
Under the premise of Theorem \ref{th:convex-f-prac}, let $a_k \equiv a > 0$ for all $k\geq 0$. Then $\min_{0 \leq k \leq K} F(x_k) - F(x_\star) \leq a^2$ for any $K \geq \frac{L\|x_0-x_\star\|^2}{2a^2}$. In particular, for a given $\epsilon>0$, by choosing $t_k \equiv \frac{\sqrt{2\epsilon}}{\sqrt{L}}$ for all $k\geq 0$, we have $\min_{0 \leq k \leq K} F(x_k) - F(x_\star) \leq \epsilon$ for any $K \geq \frac{L\|x_0-x_\star\|^2}{2\epsilon}$. 
\end{corollary}

\paragraph{Adaptive and halving radii.} Corollary~\ref{co:convex-f-prac} exhibits a trade-off: a higher precision requires a smaller radius, while too small a radius slows down the convergence. In our experiments (Section~\ref{sec:exp}), the following adaptive rule, which follows a decaying schedule and halves its constant whenever no progress has been made, performs well:
\begin{equation}\label{eq:adaptive}
	t_k = \max\{ t_0 2^{-h_k},\, l_t \}\, s_k, \qquad s_k \eqdef \sqrt{\tfrac{2\ln 2}{(k+2)\ln(k+2)}},
\end{equation}
where $t_0>0$, $l_t>0$ is a lower bound on the constant, and $h_k$ counts the checkpoints $j\leq k$, taken every $T$ iterations, at which $F(x_j) \geq F(x_{j-T})$. Since $l_t s_k \leq t_k \leq t_0 s_k$, we have $\sum_k t_k^2 = +\infty$ and $t_k \to 0$, so the rule is covered by Theorem~\ref{th:convex-f-prac} with $a_{k+1} = \sqrt{\nicefrac{L}{2}}\, t_k$, whose claim (ii) gives $F(x_k) - F_\star \leq \frac{L}{2}t_{k-1}^2$ for infinitely many $k$. This guarantee holds along a subsequence only, and the prescribed decay $s_k$ limits the attainable speed. To obtain an $\mathcal{O}(1/K)$ rate for the last iterate, we drop the prescribed decay and halve the radius whenever the decrease of the objective over the last $T$ iterations is insufficient: given $t_0>0$, an integer $T\geq1$ and $c>0$, we set
\begin{equation}\label{eq:halving}
	t_k = \tau_{\lfloor k/T \rfloor}, \qquad \tau_0 = t_0, \qquad
	\tau_{i+1} = \begin{cases} \nicefrac{\tau_i}{2}, & \text{if } F(x_{iT}) - F(x_{(i+1)T}) < c\,T\,\tau_i^2, \\ \tau_i, & \text{otherwise}. \end{cases}
\end{equation}
The radius is thus constant between two consecutive checkpoints. The rule uses objective values only: it requires neither $L$, nor any information on $x_\star$, nor a target accuracy.

\begin{theorem}\label{th:convex-halving}
	Let Assumption \ref{as:basic} hold, and assume that $f$ is $L$-smooth and convex on an open convex set containing $\dom g$. Let \(\{x_k\}_{k \geq 0}\) be the iterates generated by Broximal GD with the radii \eqref{eq:halving}, and let $L_c \eqdef \max\{L, 2c\}$ and
	$$
	\bar R \eqdef \operatorname{dist}(x_0, \cX_\star) + T\left( 4 + \frac{L_c}{c} \right) t_0,
	$$
	where $\operatorname{dist}(x, \cX_\star) \eqdef \min_{z \in \cX_\star}\|x - z\|$. Then
	\begin{itemize}
		\item[(i)] for every $J \geq 0$, the radius $t_0 2^{-J}$ is used in at most $N_J$ iterations, where $N_0 \eqdef T + \frac{F(x_0) - F_\star}{c\, t_0^2}$ and $N_J \eqdef T + \frac{2L_c}{c} + \frac{2^{J+1} L_c \bar R}{c\, t_0}$ for $J \geq 1$;
		\item[(ii)] for every $J \geq 0$ and every $K \geq \sum_{j=0}^{J} N_j$, we have
		$$
		F(x_K) - F_\star \leq L_c\, t_0 2^{-J} \left( \bar R + \tfrac{1}{2} t_0 2^{-J} \right);
		$$
		moreover, $\|x_k - x_\star\| \leq \bar R$ for all $k \geq 0$, where $x_\star$ is the projection of $x_0$ onto $\cX_\star$;
		\item[(iii)] for every $0 < \epsilon \leq L_c\bar R t_0$, we have $F(x_K) - F_\star \leq \epsilon$ for all
		$$
		K \geq T + \frac{F(x_0) - F_\star}{c\, t_0^2} + J_\epsilon \left(T + \frac{2L_c}{c}\right) + \frac{9 L_c^2 \bar R^2}{c\, \epsilon}, \qquad J_\epsilon \eqdef \left\lceil \log_2 \frac{9L_c \bar R t_0}{8\epsilon} \right\rceil.
		$$
	\end{itemize}
\end{theorem}

The proof distinguishes two kinds of steps. With $y_k \eqdef \operatorname{Prox}_{g/L_c}(x_k - \nabla f(x_k)/L_c)$, a step with $t_k \leq \|y_k - x_k\|$ is a proximal gradient step with a stepsize of at most $\nicefrac{1}{L_c}$ (Theorem~\ref{th:recoverPGD}) and decreases $F$ by at least $\frac{L_c}{2}t_k^2$, whereas otherwise $F(x_{k+1}) - F_\star \leq L_c t_k (\operatorname{dist}(x_k, \cX_\star) + \nicefrac{t_k}{2})$; hence the radius is halved only when the objective gap is already of the order of the radius. No boundedness of $\cX_\star$ or of the sublevel sets of $F$ is required, since the distance to a solution can only increase in a step after which the gap is at most $\frac{L}{2}t_k^2$. For $c = \nicefrac{L}{2}$, which minimizes $\nicefrac{L_c^2}{c}$, we have $\bar R = \operatorname{dist}(x_0, \cX_\star) + 6Tt_0$ and the leading term in (iii) is $\nicefrac{18L\bar R^2}{\epsilon}$, of the same order as the bound $\nicefrac{L\|x_0-x_\star\|^2}{2\epsilon}$ of PGD and of Corollary~\ref{co:convex-f-prac} when $t_0 = \mathcal{O}(\operatorname{dist}(x_0, \cX_\star)/T)$, but for the last iterate and without knowledge of $L$ or of the target accuracy; any $c>0$ is allowed, and a misspecified $c$ costs a factor of order $\max\{\nicefrac{L}{c}, \nicefrac{c}{L}\}$.

\subsection{$F$ with $G$-bounded subgradients} 

We next drop the smoothness of $f$ and assume instead that $F$ has bounded subgradients on $\dom g$. This is the composite counterpart of the classical setting of the (projected) subgradient method.

\begin{theorem}\label{th:convex-G-prac}
	Let Assumption \ref{as:basic} hold, let $f$ be convex on $\dom g$, and assume that for all $x \in \dom  g$, there exists a subgradient of $g$ at $x$, and there exists $G>0$ such that $\|\nabla f(x) + \partial g(x)\| \leq G$ for all $x\in \dom  g$, where $\partial g(x)$ denotes a fixed choice of subgradient of $g$ at $x$. Let \(\{x_k\}_{k \geq 0}\) be the iterates generated by Broximal GD. Further, let $x_\star \in \cX_\star$, and assume that whenever $\|\nabla f(x_k) + \partial g(x_k)\| >0$, we choose $t_k = \frac{a_k}{\|\nabla f(x_k) + \partial g(x_k)\|}$, where $\{  a_k  \}_{k=0}^{+\infty}$ is any positive sequence such that $\sum_{k=0}^{+\infty} a_k^2 = +\infty$, and $t_k = 0$ when $\|\nabla f(x_k) + \partial g(x_k)\| =0$. Then 
	\begin{itemize}
		\item[(i)] for any $k \notin \cK$, where $\cK \eqdef \{  k |  F(x_k) - F(x_\star) \leq a_k  \}$, we have $t_k \in \cT_k$ and 
		$$
		\|x_{k+1} - x_\star\|^2 < \|x_k - x_\star\|^2 - \frac{a_k^2}{G^2}; 
		$$
		
		\item[(ii)] there exists a subsequence $\{  x_{j_k}  \}_{k=0}^{+\infty}$ with $0\leq j_0<j_1<\dots$, such that $F(x_{j_k}) - F(x_\star) \leq a_{j_k}$ for all $k\geq 0$; 
		
		\item[(iii)] if instead $a_k = F(x_k) - F(x_\star)$ for all $k \geq 0$ (a choice that need not be positive or satisfy $\sum_{k=0}^{+\infty} a_k^2 = +\infty$), we have $\frac{1}{K} \sum_{k=0}^{K-1} (F(x_k) - F(x_\star))^2 \leq \frac{G^2\|x_0 - x_\star\|^2}{K}$ for any $K\geq 1$. 
	\end{itemize}

\end{theorem}

Theorem~\ref{th:convex-G-prac} is the composite analogue of the normalized subgradient method: the radius needs neither a smoothness constant nor $G$, which enters only the analysis, and the choice $a_k = F(x_k) - F_\star$ in (iii) is the analogue of the Polyak stepsize \citep{polyak1969minimization}. For instance, we can choose the sequence $\{  a_k  \}_{k=0}^{+\infty}$ in Theorem \ref{th:convex-G-prac} as $a_k = \sqrt{\frac{1}{(k+2) \ln (k+2)}}$. The following corollary shows that a constant $a_k$ recovers the optimal $\mathcal{O}(G\|x_0 - x_\star\|/\sqrt{K})$ rate of the subgradient method \citep{nesterov2018lectures}.

\begin{corollary}\label{co:convex-G-prac}
	Under the premise of Theorem \ref{th:convex-G-prac}, let $a_k \equiv a > 0$ for all $k\geq 0$. Then $\min_{0 \leq k \leq K-1} F(x_k) - F(x_\star) \leq a$ for any $K \geq \frac{G^2\|x_0-x_\star\|^2}{a^2}$. In particular, for a given $\epsilon>0$, by choosing $t_k \equiv \frac{\epsilon}{\|\nabla f(x_k) + \partial g(x_k)\|}$ for all $k\geq 0$, we have $\min_{0 \leq k \leq K-1} F(x_k) - F(x_\star) \leq \epsilon$ for any $K \geq \frac{G^2\|x_0-x_\star\|^2}{\epsilon^2}$. 
\end{corollary}

\section{Strongly convex case}

We now assume that $f$ is $\mu$-strongly convex and $L$-smooth. As before, we first analyze an idealized radius rule that depends on $\|x_k - x_\star\|$, then an implementable geometric schedule, and finally the halving rule \eqref{eq:halving}.

\begin{theorem}\label{th:strongly-convex-theory}
Let Assumption \ref{as:basic} hold, and assume that $f$ is $L$-smooth and $\mu$-strongly convex on an open convex set containing $\dom  g$. Let $\{x_k\}_{k \ge 0}$ be the iterates generated by Broximal GD. Further, let $x_\star \in \mathcal{X}_\star$ and assume that whenever $\nabla f(x_k) \neq \nabla f(x_\star)$, we choose the radius 
$$
t_k = \theta \|x_k - x_\star\|, \quad \theta = \frac{2\sqrt{\mu L}}{L+\mu}, 
$$
and $t_k=0$ when $\nabla f(x_k) = \nabla f(x_\star)$. Then $\|x_{k+1} - x_\star\|^2 \leq (1 - \theta^2) \|x_k - x_\star\|^2 = \left( \frac{L-\mu}{L+\mu} \right)^2 \|x_k - x_\star\|^2$ for any $k\geq 0$. Hence, for any $K \geq 0$, 
$$
\|x_K - x_\star\|^2 \leq \left( \frac{L-\mu}{L+\mu} \right)^{2K} \|x_0 - x_\star\|^2. 
$$
\end{theorem}

Note that $\nabla f(x_k) = \nabla f(x_\star)$ if and only if $x_k = x_\star$ under strong convexity, so the case $t_k = 0$ is only reached at the solution.

\begin{theorem}\label{th:strongly-convex-prac}
Let Assumption \ref{as:basic} hold, and assume that $f$ is $L$-smooth and $\mu$-strongly convex on an open convex set containing $\dom  g$. Let $\{x_k\}_{k \ge 0}$ be the iterates generated by Broximal GD. Further, let $x_\star \in \mathcal{X}_\star$, and let the radii be $t_k = \theta C q^k$, where $\theta = \frac{2\sqrt{\mu L}}{L + \mu}$, $C>0$, $q^2 = 1-(1-\alpha) \theta^2$ with $0<\alpha<1$. Assume $1-\alpha < \frac{1}{1+2\theta+2\theta^2}$. Then the following results hold. 
\begin{itemize}
	
	\item[(i)] For any $k \notin \cK$, where $\cK \eqdef \{  k |  \|x_k - x_\star\| \leq Cq^k  \}$, we have $t_k \in \cT_k$ and 
	$$
	\|x_{k+1} - x_\star\|^2 < \|x_k - x_\star\|^2 - \theta^2 C^2 q^{2k}. 
	$$
	
	\item[(ii)] For $k\geq 1$ such that $k-1 \in \cK$, we have 
	$$
	\|x_k - x_\star\|^2 \leq (1+2\theta+\theta^2) C^2 q^{2k-2}.
	$$
	
	\item[(iii)] If $C > \sqrt{1-\alpha} \|x_0-x_\star\|$, then there exists a subsequence $\{  x_{j_k}  \}_{k=0}^{+\infty}$ with $0\leq j_0<j_1<\dots$, such that $\|x_{j_k} - x_\star\| \leq Cq^{j_k}$ for all $k\geq 0$. 
	
	\item[(iv)] Let $\epsilon>0$ and $C \geq \sqrt{1-\alpha} \|x_0-x_\star\|$. Then we have $\|x_K - x_\star\| \leq \epsilon$ for any $K \geq \frac{1}{1-q^2} \ln \frac{\theta^2 C^2}{\epsilon^2(1-q^2)}$. In particular, by choosing $\alpha = \frac{2\theta + 3\theta^2}{1+2\theta+3\theta^2}$ and $C = \sqrt{1-\alpha} \|x_0-x_\star\|$, we have $\|x_K - x_\star\| \leq \epsilon$ for any $K \geq \left(  \frac{1}{\theta^2} + \frac{2}{\theta} + 3  \right) \ln \frac{\|x_0-x_\star\|^2}{\epsilon^2}$. 
	
\end{itemize}

\end{theorem}

Theorem~\ref{th:strongly-convex-prac} replaces the unknown distance $\|x_k - x_\star\|$ by the geometric envelope $Cq^k$, which requires only $\theta$ and an upper bound $C$ on $\sqrt{1-\alpha}\,\|x_0 - x_\star\|$: at each iteration, either the iterate is inside the envelope, or the radius $\theta Cq^k$ is admissible and the squared distance decreases by $\theta^2C^2q^{2k}$, while (ii) bounds how far an iterate can leave the envelope, which is what makes (iv) a last-iterate statement. The price is small: the bound in (iv) has the same leading term $\frac{1}{\theta^2}$ as the idealized rule of Theorem~\ref{th:strongly-convex-theory}.

\paragraph{What if the parameter $\theta$ is unknown?} Suppose that $L$ is known but $\mu$, and hence $\theta$, is not, and that an upper bound ${\tilde C} \geq \|x_0-x_\star\|$ is available. The gradient mapping provides a computable certificate. Let $x \in \dom g$, $\gamma = \nicefrac{1}{L}$ and $y = \operatorname{Prox}_{\gamma g}(x - \gamma\nabla f(x))$. Since $s \eqdef G_\gamma(x) - \nabla f(x) + \nabla f(y) \in \partial F(y)$ with $\|s\| \leq 2\|G_\gamma(x)\|$, the $\mu$-strong convexity of $F$ gives
\begin{equation}\label{eq:certificate}
\|G_{\nicefrac{1}{L}}(x)\| \leq \epsilon_G \quad\Longrightarrow\quad \|x - x_\star\| \leq \frac{3\epsilon_G}{\mu} \quad\text{and}\quad F(y) - F_\star \leq \frac{2\epsilon_G^2}{\mu},
\end{equation}
and, conversely, $\|G_\gamma(x)\| \leq 3L\|x-x_\star\|$, so the certificate holds at every $x$ with $\|x - x_\star\| \leq \epsilon \eqdef \nicefrac{\epsilon_G}{3L}$. One can therefore run BGD repeatedly from $x_0$, the $i$-th run using $t_k = {\tilde C}q_i^k$ with $1-q_i^2 = 2^{-i}$ for $K_i = \lceil 2^i \ln \frac{2^i {\tilde C}^2}{\epsilon^2} \rceil$ iterations, and stop as soon as $\|G_{\nicefrac{1}{L}}(x_{K_i})\| \leq \epsilon_G$. The $i$-th run corresponds to $C = \nicefrac{{\tilde C}}{\theta}$ and $1-\alpha = \nicefrac{2^{-i}}{\theta^2}$ in Theorem~\ref{th:strongly-convex-prac}, whose claim (iv) guarantees a passed test whenever $2^{-i} < \frac{\theta^2}{1+2\theta+2\theta^2}$. The procedure therefore stops after at most $n \eqdef \lfloor \log_2 \frac{1+2\theta+2\theta^2}{\theta^2} \rfloor + 1$ runs and
$$
\sum_{i=1}^{n} K_i \leq n + 4\left(\frac{1}{\theta^2} + \frac{2}{\theta} + 3\right) \ln \frac{2^n {\tilde C}^2}{\epsilon^2}
$$
iterations (for ${\tilde C} > \epsilon$), i.e., about four times the bound of Theorem~\ref{th:strongly-convex-prac}~{\em (iv)} for known $\theta$. It still requires $L$, ${\tilde C}$ and a target $\epsilon_G$, while the resulting accuracy \eqref{eq:certificate} depends on the unknown $\mu$.

\paragraph{Halving radii.} The halving rule \eqref{eq:halving} requires neither $\mu$, nor $L$, nor a bound on $\|x_0 - x_\star\|$, and it adapts to strong convexity.

\begin{theorem}\label{th:sc-halving}
	Let Assumption \ref{as:basic} hold, and assume that $f$ is $L$-smooth and $\mu$-strongly convex on an open convex set containing $\dom g$. Let \(\{x_k\}_{k \geq 0}\) be the iterates generated by Broximal GD with the radii \eqref{eq:halving}, and let $L_c \eqdef \max\{L, 2c\}$,
	$$
	N_0 \eqdef T + \frac{F(x_0) - F_\star}{c\, t_0^2}, \qquad N \eqdef T + \frac{2L_c}{c}\left( 1 + \frac{L_c}{\mu} \right).
	$$
	Then the radius $t_0$ is used in at most $N_0$ iterations and every radius $t_0 2^{-J}$ with $J \geq 1$ in at most $N$ iterations, and for every $J \geq 0$ and every $K \geq N_0 + JN$,
	$$
	\frac{\mu}{2}\|x_K - x_\star\|^2 \leq F(x_K) - F_\star \leq \frac{L_c t_0^2}{2} \left( 1 + \frac{L_c}{\mu} \right) 4^{-J}.
	$$
\end{theorem}

Thus the objective gap is divided by four every $N$ iterations. For $c = \nicefrac{L}{2}$, which minimizes $\nicefrac{L_c^2}{c}$, we have $N = \nicefrac{4L}{\mu} + T + 4$, and, for every $0 < \epsilon \leq \frac{1}{2}L t_0^2 (1 + \nicefrac{L}{\mu})$, $F(x_K) - F_\star \leq \epsilon$ for all
$$
K \geq N_0 + \left( \frac{4L}{\mu} + T + 4 \right) \left\lceil \log_4 \frac{L t_0^2 (1 + L/\mu)}{2\epsilon} \right\rceil,
$$
whose leading term is $\frac{4}{\ln 4}\frac{L}{\mu}\ln\frac{1}{\epsilon} \approx 2.9\,\frac{L}{\mu}\ln\frac{1}{\epsilon}$, against $\frac{L}{\mu}\ln\frac{1}{\epsilon}$ for PGD with stepsize $\nicefrac{1}{L}$. The rule itself involves no problem constant, and, as in the convex case, a misspecified $c$ costs a factor of order $\max\{\nicefrac{L}{c}, \nicefrac{c}{L}\}$. Compared with the restart procedure above, no upper bound on $\|x_0 - x_\star\|$ and no stopping certificate are needed.

\section{Non-convex case}

We adopt the proximal P{\L} condition of \citet{karimi2016linear}. For $x \in \dom g$ and $\gamma > 0$, define
\begin{equation}\label{eq:Dg}
\cD_g(x, \gamma) \eqdef -\frac{2}{\gamma} \min_{y \in \R^d} \left\{ \langle \nabla f(x), y - x \rangle + \frac{1}{2\gamma}\|y - x\|^2 + g(y) - g(x) \right\}.
\end{equation}
The objective in \eqref{eq:Dg} is strongly convex, and its unique minimizer is $y = \operatorname{Prox}_{\gamma g}(x - \gamma \nabla f(x))$; choosing $y = x$ shows that $\cD_g(x, \gamma) \geq 0$. 

\begin{assumption}[Proximal P\L]\label{as:proximalPL}
Let Assumption~\ref{as:basic} hold. Furthermore, assume that there exist constants $\mu > 0$ and $\gamma > 0$ such that
$$
\frac{1}{2}\cD_g(x, \gamma) \geq \mu \left(F(x) - F_\star \right), \qquad \forall x \in \dom g. 
$$
\end{assumption}

If $g$ is constant, then $\cD_g(x, \gamma) = \|\nabla f(x)\|^2$, and Assumption~\ref{as:proximalPL} reduces to the classical P{\L} inequality. In general, $\cD_g(x, \gamma) \geq \|G_\gamma(x)\|^2$ for all $x \in \dom g$ (see Lemma~\ref{lm:non-convex-step}(i) in Appendix~\ref{sec:proof-nonconvex}), where $G_\gamma$ is the gradient mapping introduced in Section~\ref{sec:prelim}, whose norm $\|G_\gamma(x)\|$ is the standard measure of stationarity for composite problems \citep{nesterov2013gradient}; thus Assumption~\ref{as:proximalPL} is weaker than requiring $\frac{1}{2}\|G_\gamma(x)\|^2 \geq \mu (F(x) - F_\star)$. For instance, for $d = 1$, $f(x) = x$ and $g = \ind_{[0,1]}$, we have $\cD_g(x, 1) = 2x - x^2 \geq F(x) - F_\star$ and $\|G_1(x)\|^2 = x^2$ for $x \in [0, 1]$, so Assumption~\ref{as:proximalPL} holds with $\gamma = 1$ and $\mu = \nicefrac{1}{2}$, while $\frac{1}{2}\|G_1(x)\|^2 \geq \mu (F(x) - F_\star)$ fails near $x = 0$ for every $\mu > 0$.

\subsection{Proximal P\L}

\begin{theorem}\label{th:non-convex-PL-theory}
Let Assumption \ref{as:proximalPL} hold for some $\gamma \in (0, \nicefrac{1}{L}]$, and let $f$ be $L$-smooth on an open set containing $\dom  g$. Let $\{x_k\}_{k \geq 0}$ be the iterates generated by Broximal GD with radii $t_k = \gamma \|G_\gamma (x_k)\|$. Then for any $k\geq 0$, 
$$
F(x_k) - F_\star \leq (1 - \gamma \mu)^k (F(x_0) - F_\star). 
$$
\end{theorem}

Since the radius $t_k = \gamma\|G_\gamma(x_k)\|$ makes the broximal step coincide with the proximal gradient step (Theorem~\ref{th:recoverPGD}(i)), this is the linear rate of PGD under the proximal P{\L} condition \citep{karimi2016linear}, restated in terms of radii; it serves as a baseline for the following result, in which the gradient-dependent radius is replaced by a prescribed geometric schedule.

\begin{theorem}\label{th:non-convex-PL}
Let Assumption \ref{as:proximalPL} hold for some $\gamma \in (0, \nicefrac{1}{L}]$, and let $f$ be $L$-smooth on an open set containing $\dom  g$. Let $\{x_k\}_{k \geq 0}$ be the iterates generated by Broximal GD with radii $t_k = \gamma \sqrt{\nicefrac{\mu Cq^k}{2}}$, where $C>0$, $q = 1-\frac{3}{4}(1-\alpha) \gamma \mu$ with $0<\alpha<1$. Assume $1-\alpha < \frac{1}{1+\gamma \mu}$. Then the following results hold. 
\begin{itemize}
	
	\item[(i)] For any $k \notin \cK$, where $\cK \eqdef \{  k |  F(x_k) - F_\star \leq Cq^k  \}$, we have 
	$$
	F(x_{k+1}) - F_\star < F(x_k) - F_\star - \frac{3}{4}\gamma \mu Cq^k. 
	$$
	
	\item[(ii)] For $k\geq 1$ such that $k-1 \in \cK$, we have 
	$$
	F(x_k) - F_\star \leq \left(1+ \frac{\gamma \mu}{4} \right) Cq^{k-1}. 
	$$
	
	\item[(iii)] If $C > (1-\alpha) (F(x_0) - F_\star)$, then there exists a subsequence $\{  x_{j_k}  \}_{k=0}^{+\infty}$ with $0\leq j_0<j_1<\dots$, such that $F(x_{j_k}) - F_\star \leq Cq^{j_k}$ for all $k\geq 0$. 
	
	\item[(iv)] Let $\epsilon>0$ and $C \geq  (1-\alpha) (F(x_0) - F_\star)$. Then we have $F(x_K) - F_\star \leq \epsilon$ for any $K \geq \frac{1}{1-q} \ln \frac{3\gamma \mu C}{4\epsilon(1-q)}$. In particular, by choosing $\alpha = \frac{3\gamma \mu}{2+3\gamma \mu}$ and $C = (1-\alpha) (F(x_0) - F_\star)$, we have $F(x_K) - F_\star \leq \epsilon$ for any $K \geq \left(  \frac{4}{3\gamma \mu}+ 2  \right) \ln \frac{F(x_0) - F_\star}{\epsilon}$. 
	
\end{itemize}

\end{theorem}

Here progress is measured in function value, and a prescribed radius may exceed the proximal gradient displacement $\gamma\|G_\gamma(x_k)\|$, in which case $F$ may increase, by at most $\frac{t_k^2}{2\gamma}$. Lemma~\ref{lm:non-convex-step}(ii) in Appendix~\ref{sec:proof-nonconvex} quantifies the trade-off, $F(x_{k+1}) \leq F(x_k) - \min\{ t_k \sqrt{\cD_g(x_k,\gamma)}, \frac{\gamma}{2}\cD_g(x_k,\gamma) \} + \frac{t_k^2}{2\gamma}$, and the radius $t_k = \gamma\sqrt{\nicefrac{\mu C q^k}{2}}$ balances the two terms against the penalty, giving the per-step decrease $\frac{3}{4}\gamma\mu Cq^k$. The resulting complexity is within a factor of $\nicefrac{4}{3}$ of that of Theorem~\ref{th:non-convex-PL-theory}, up to a lower-order term, with radii fixed in advance and without ever evaluating $G_\gamma(x_k)$ or $\cD_g(x_k, \gamma)$.

The halving rule \eqref{eq:halving} also attains a linear rate under the proximal P{\L} condition, without knowledge of $\mu$, $L$ or $\gamma$.

\begin{theorem}\label{th:pl-halving}
	Let Assumption \ref{as:proximalPL} hold for some $\gamma \in (0, \nicefrac{1}{L}]$, and let $f$ be $L$-smooth on an open set containing $\dom g$. Let $\{x_k\}_{k \geq 0}$ be the iterates generated by Broximal GD with the radii \eqref{eq:halving}, and let $s \eqdef c + \frac{1}{2\gamma}$, $\kappa \eqdef \max\{ s^2, \nicefrac{2s}{\gamma} \}$,
	$$
	N_0 \eqdef T + \frac{F(x_0) - F_\star}{c\, t_0^2}, \qquad N \eqdef T + \frac{2}{c}\left( \frac{\kappa}{\mu} + \frac{1}{\gamma} \right).
	$$
	Then the radius $t_0$ is used in at most $N_0$ iterations and every radius $t_0 2^{-J}$ with $J \geq 1$ in at most $N$ iterations, and for every $J \geq 0$ and every $K \geq N_0 + JN$,
	$$
	F(x_K) - F_\star \leq \left( \frac{\kappa}{2\mu} + \frac{1}{2\gamma} \right) t_0^2\, 4^{-J}.
	$$
\end{theorem}

For $\gamma = \nicefrac{1}{L}$ and $c \leq \nicefrac{3L}{2}$, we have $\kappa = 2L(c + \nicefrac{L}{2})$ and $N = T + \frac{4L}{\mu} + \frac{2L^2}{c\mu} + \frac{2L}{c}$, which decreases in $c$, whereas $N$ increases in $c$ for $c \geq \nicefrac{3L}{2}$; the choice $c = \nicefrac{3L}{2}$ gives $N = \frac{16L}{3\mu} + T + \frac{4}{3}$, i.e., a leading term $\frac{16}{3\ln 4}\frac{L}{\mu}\ln\frac{1}{\epsilon} \approx 3.8\,\frac{L}{\mu}\ln\frac{1}{\epsilon}$.

\subsection{Without the proximal P{\L} condition}

Finally, we drop the proximal P{\L} condition and only assume that $f$ is $L$-smooth. In this setting one can only hope to find an approximately stationary point, i.e., a point with a small gradient mapping.

\begin{theorem}\label{th:non-convex-withoutPL}
	Let Assumption \ref{as:basic} hold, and assume that $f$ is $L$-smooth on an open set containing $\dom  g$. Let \(\{x_k\}_{k \geq 0}\) be the iterates generated by Broximal GD. Further, let $\gamma \in (0, \nicefrac{1}{L}]$, and let the radii be $t_k =\gamma a_k$, where $\{  a_k  \}_{k=0}^{+\infty}$ is any positive sequence such that $\sum_{k=0}^{+\infty} a_k^2 = +\infty$. Denote $y_k \eqdef \operatorname{Prox}_{\gamma g} \left( x_k - \gamma \nabla f(x_k) \right)$. Then 
	\begin{itemize}
		\item[(i)] for any $k \notin \cK$, where $\cK \eqdef \{  k |  \|G_{\gamma} (x_k)\| \leq a_k  \}$, we have $t_k<\gamma \|G_{\gamma} (x_k)\| $, i.e., $t_k \in \cT_k(y_k)$, and 
		$$
		F(x_{k+1}) - F_\star < F(x_k) - F_\star - \frac{\gamma}{2} a_k^2; 
		$$
		
		\item[(ii)] there exists a subsequence $\{  x_{j_k}  \}_{k=0}^{+\infty}$ with $0\leq j_0<j_1<\dots$, such that $\|G_{\gamma} (x_{j_k})\|  \leq a_{j_k}$ for all $k\geq 0$; 
		
		\item[(iii)] if instead $a_k =  \|G_{\gamma} (x_k)\|$ for all $k \geq 0$ (a choice that need not be positive or satisfy $\sum_{k=0}^{+\infty} a_k^2 = +\infty$), we have $\frac{1}{K} \sum_{k=0}^{K-1}  \|G_{\gamma} (x_k)\|^2 \leq \frac{2(F(x_0) - F_\star)}{\gamma K}$ for any $K\geq 1$. 
	\end{itemize}
	
\end{theorem}

For instance, one can take $a_k = \sqrt{\frac{1}{(k+2) \ln (k+2)}}$. The following corollary shows that a constant $a_k$ attains the standard $\mathcal{O}(L(F(x_0) - F_\star)/\epsilon^2)$ complexity of first-order methods for reaching an $\epsilon$-stationary point of smooth nonconvex composite problems \citep{nesterov2013gradient}.

\begin{corollary}\label{co:non-convex-withoutPL}
	Under the premise of Theorem \ref{th:non-convex-withoutPL}, let $a_k \equiv a > 0$ for all $k\geq 0$. Then $\min_{0 \leq k \leq K-1} \|G_{\gamma} (x_k)\| \leq a$ for any $K \geq  \frac{2(F(x_0) - F_\star)}{\gamma a^2}$. In particular, for a given $\epsilon>0$, by choosing $t_k \equiv \frac{\epsilon}{L}$ for all $k\geq 0$, we have $\min_{0 \leq k \leq K-1}  \|G_{1/L} (x_k)\|  \leq \epsilon$ for any $K \geq \frac{2L(F(x_0) - F_\star)}{\epsilon^2}$. 
\end{corollary}

The halving rule \eqref{eq:halving} attains the same complexity without knowledge of $L$ or of the target accuracy.

\begin{theorem}\label{th:nc-halving}
	Let Assumption \ref{as:basic} hold, and assume that $f$ is $L$-smooth on an open convex set containing $\dom g$ and that $F_\star \eqdef \inf F > -\infty$. Let \(\{x_k\}_{k \geq 0}\) be the iterates generated by Broximal GD with the radii \eqref{eq:halving}, and let $L_c \eqdef \max\{L, 2c\}$ and $\bar\Phi \eqdef F(x_0) - F_\star + \frac{2}{3} T L_c t_0^2$. Then, for every $J \geq 0$,
	$$
	\min_{0 \leq k < K} \|G_{\nicefrac{1}{L}}(x_k)\| < L_c t_0 2^{-J} \qquad \text{for all } K \geq (J+1)T + \frac{4^{J}\,\bar\Phi}{c\, t_0^2}.
	$$
	In particular, for every $0 < \epsilon \leq L_c t_0$, we have $\min_{0\leq k<K}\|G_{\nicefrac{1}{L}}(x_k)\| < \epsilon$ for all $K \geq (J_\epsilon + 1)T + \frac{4L_c^2\bar\Phi}{c\,\epsilon^2}$, where $J_\epsilon \eqdef \lceil \log_2 \frac{L_c t_0}{\epsilon} \rceil$; for $c = \nicefrac{L}{2}$, the leading term is $\nicefrac{8L\bar\Phi}{\epsilon^2}$.
\end{theorem}

\section{Experiments}\label{sec:exp}

We compare BGD with PGD with a constant stepsize and with the line search of \citet{nesterov2013gradient} (with increase and decrease factors $\gamma_u=\gamma_d=2$) on a convex problem (LASSO), a strongly convex problem (elastic net), and a nonconvex problem (a neural network under a box constraint); since BGD, like PGD, is not accelerated, we do not compare it with accelerated methods. For BGD, ``stepsize'' refers to the radius $t_k$, and the subproblem is solved by the procedure of Appendix~\ref{sec:subproblem}, which guarantees that each computed step is an exact BGD step with a radius in $[0.99\,t_k,t_k]$.\footnote{BGD is thus run with radii $\tilde t_k\in[0.99\,t_k,t_k]$, and the convergence results apply to $\{\tilde t_k\}$ directly. Since the radii for which the master theorem applies form an interval $(0,b)$ (Theorem~\ref{th:Tk}), a radius slightly smaller than an admissible one is admissible; the schedules of Theorems~\ref{th:convex-f-prac} and \ref{th:non-convex-withoutPL} keep their form with $\tilde a_k\in[0.99\,a_k,a_k]$, and the guaranteed decrease per step, which is proportional to the squared radius, is reduced by a factor of at most $0.99^2\approx0.98$ (for Theorem~\ref{th:strongly-convex-prac}, provided the condition on $\alpha$ holds with this margin). For the halving rule \eqref{eq:halving}, whose test is applied with the nominal radius $\tau_i$, the proofs of Theorems~\ref{th:convex-halving}, \ref{th:sc-halving} and \ref{th:nc-halving} go through with $L_c$ replaced by $\max\{L, \nicefrac{2c}{0.99^2}\}$, and that of Theorem~\ref{th:pl-halving} with $s$ replaced by $\nicefrac{s}{0.99}$, since a step with radius $\tilde t_k \geq 0.99\,\tau_i$ still makes $0.99^2$, respectively $0.99$, times the progress used there, while all upper bounds are nondecreasing in the radius. The effect on the convergence speed is therefore negligible.} Except for Figure~\ref{fig:convex-nofax}(a), results are reported over $20$ random instances (LASSO and elastic net, started from $x_0=0$) or $10$ random seeds, each determining the data subsets and the initialization (network); the parameters of all methods are selected by grid search on a separate instance (seed $0$) and then frozen, and all methods use the same instances and initial points, so that the comparisons are paired. The figures show the median and the $25\%$--$75\%$ range over instances,\footnote{At each value of the horizontal axis, the median and the quartiles are taken over the instances of the value recorded at the last iteration up to that point; a run that has stopped, after reaching the target accuracy or exhausting its budget, keeps its last value, and runs that diverge keep their last finite value. For logarithmic axes, the quantiles are computed on the logarithm. The median curve is thus not the trajectory of a particular instance; the thin lines show the individual instances.} and Appendix~\ref{sec:addexp} reports the paired comparisons. In all experiments, BGD uses about $12$ to $13$ proximal evaluations per iteration in median, against about $2$ for PGD with line search.

\subsection{LASSO}

For the convex case, we test on randomly generated LASSO instances \citep{tibshirani1996regression, nesterov2013gradient}:
\begin{equation}\label{p:lasso}
	\min_{x \in \mathbb{R}^n } F(x) \eqdef  \frac{1}{2}\|Ax-b\|^2 + \|x\|_1, 
\end{equation}
where $A = (a_1, ..., a_n) \in \R^{m\times n}$ with $m<n$. The instances are generated as in \citet{nesterov2013gradient} (cf.\ Appendix~\ref{sec:elanet}), where $m^*\leq m$ is the number of nonzero components of the optimal solution $x^*$\footnote{In the description of the test problems, we write $x^*$ and $y^*$ as in \citet{nesterov2013gradient}; elsewhere, $x_\star$ denotes a minimizer.}, and $\rho>0$ controls the size of $x^*$. We first compare PGD with constant stepsize $\gamma$ and with line search (initial estimate $L_0$) with BGD with a constant radius $t$ and with the decaying radii $\nicefrac{a}{\sqrt{k+2}}$ and $\nicefrac{a}{\sqrt{(k+2)\ln(k+2)}}$, with parameters selected by grid search over $\gamma \in \{  3^{-3}, \dots, 3^{-12}  \}$, $L_0 \in \max_i\|a_i\|^2 \times \{ 1, 3^{-1}, \dots, 3^{-9}  \}$, $t \in 10^{-2} \times \{ 1, 3^{-1}, \dots, 3^{-9} \}$, and $a \in \{ 1, 3^{-1}, \dots, 3^{-9} \}$.

\begin{figure}[t!]
	\centering
	
	\includegraphics[width=0.24\linewidth]{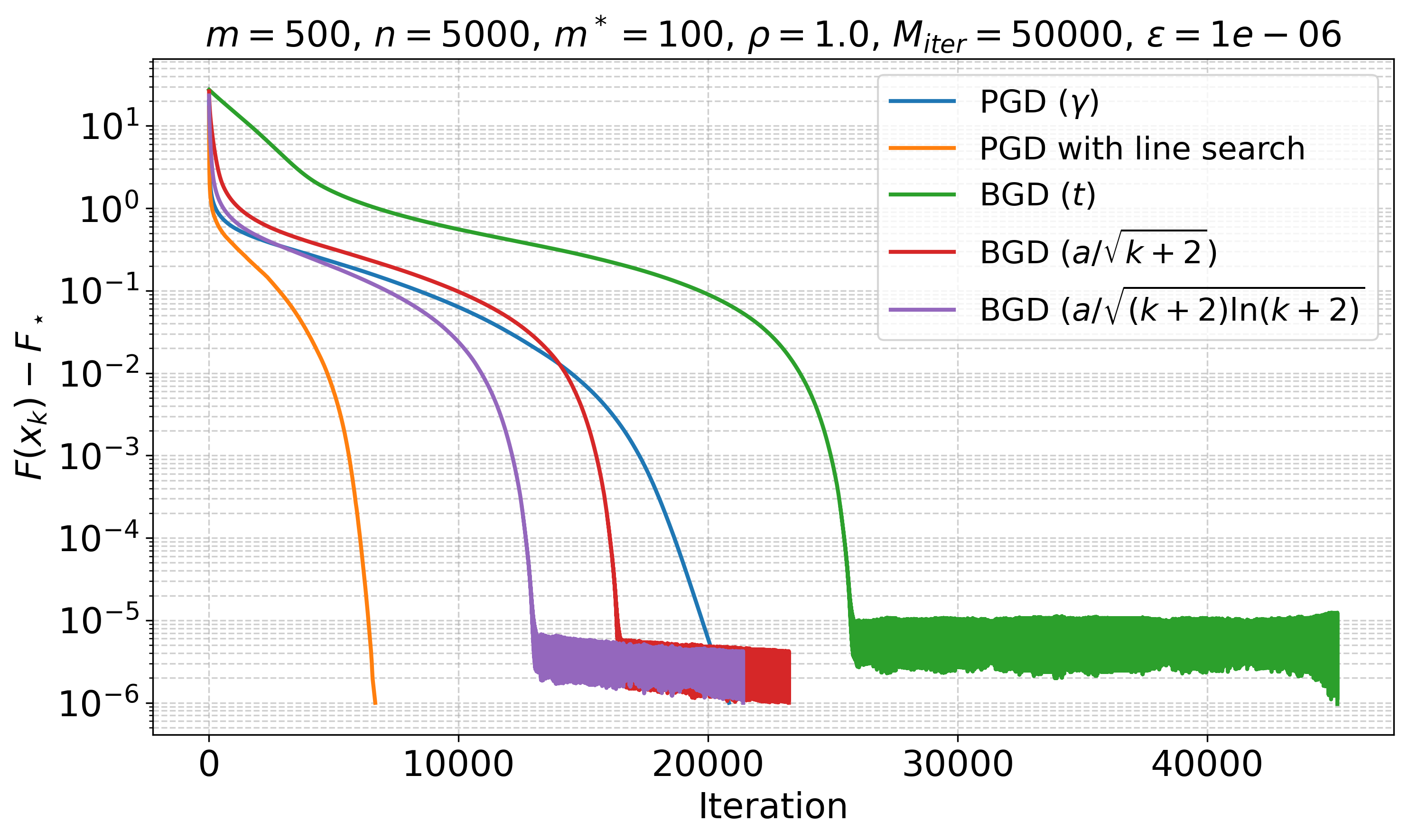}
	\hfill
	\includegraphics[width=0.24\linewidth]{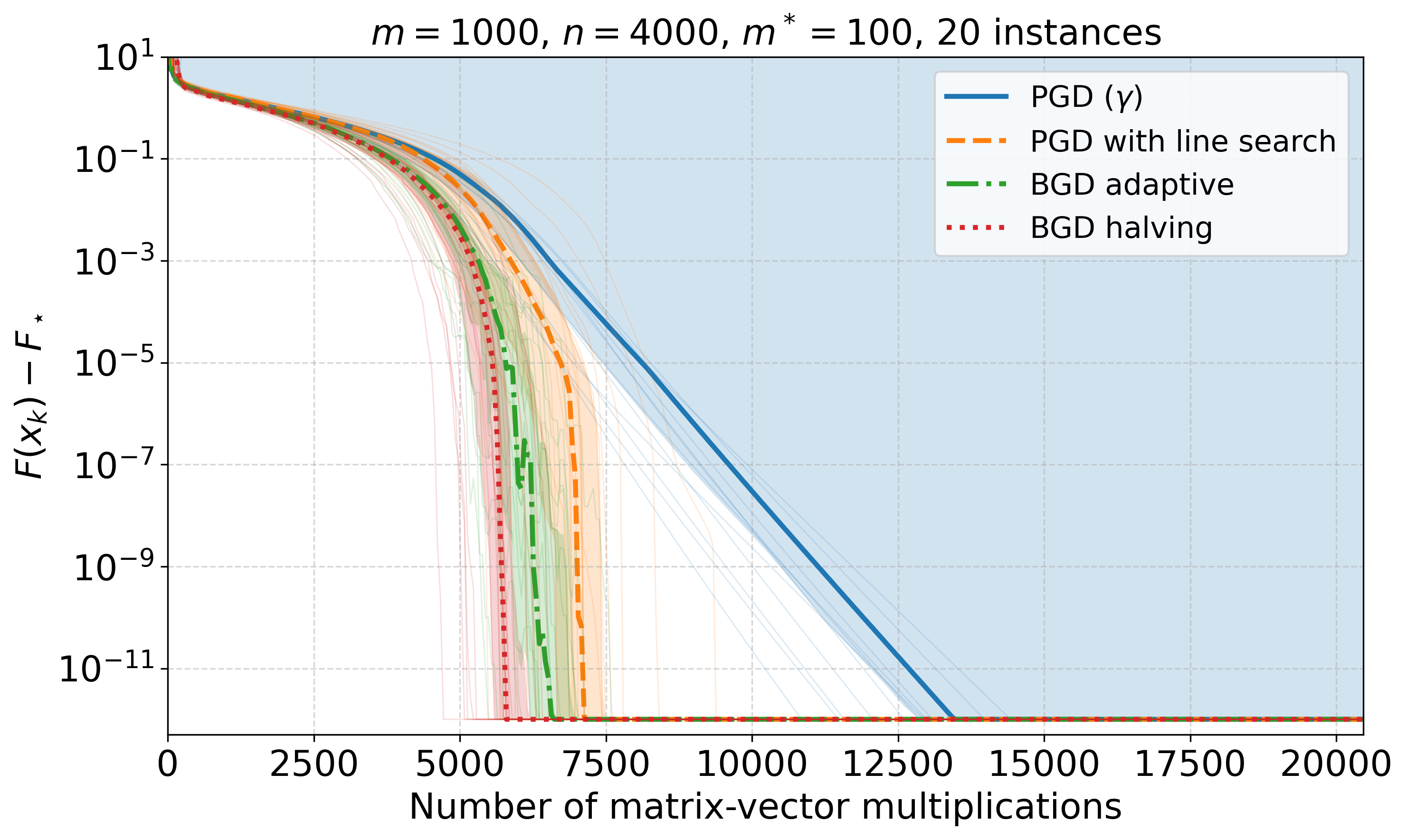}
	\hfill 
	\includegraphics[width=0.24\linewidth]{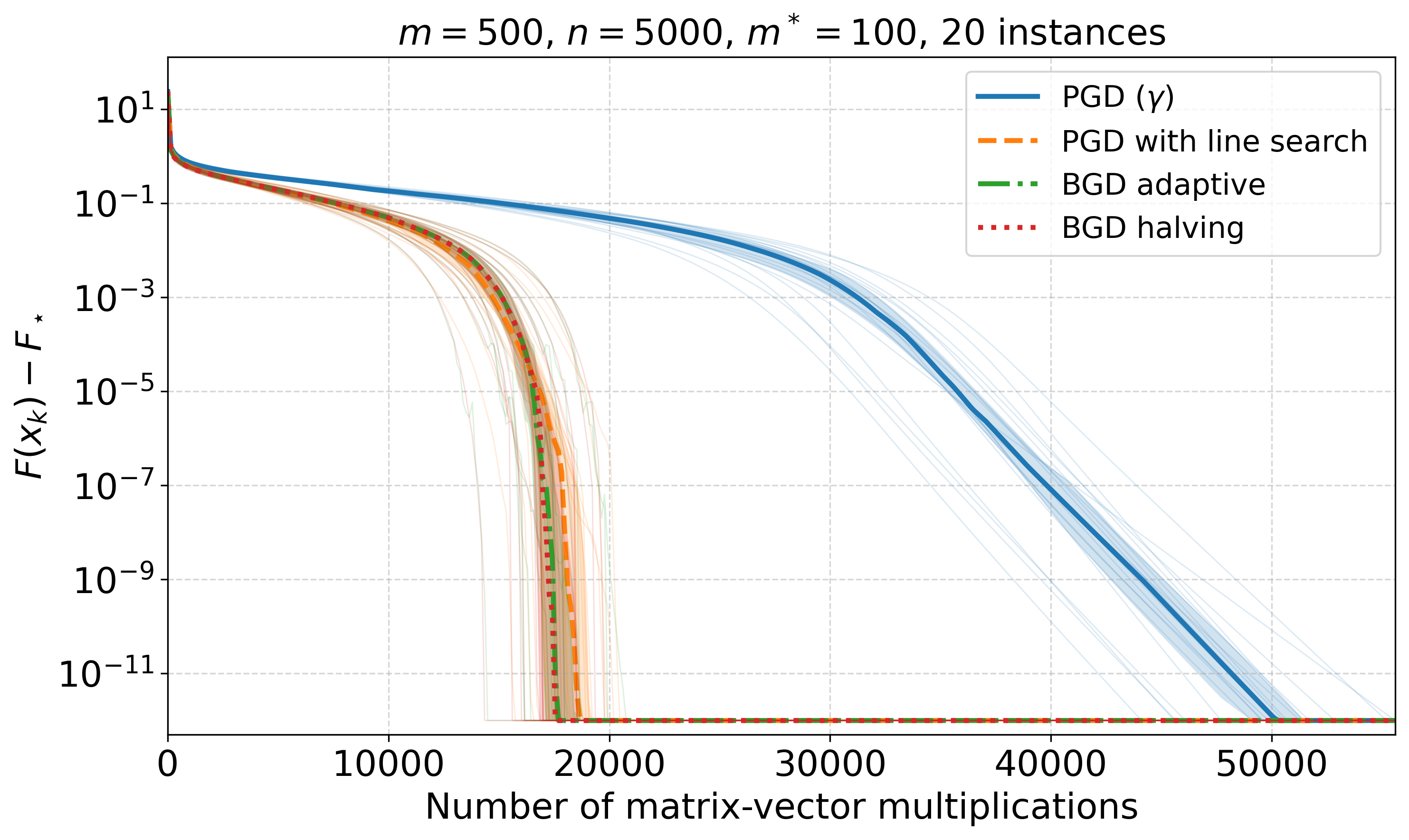}
	\hfill
	\includegraphics[width=0.24\linewidth]{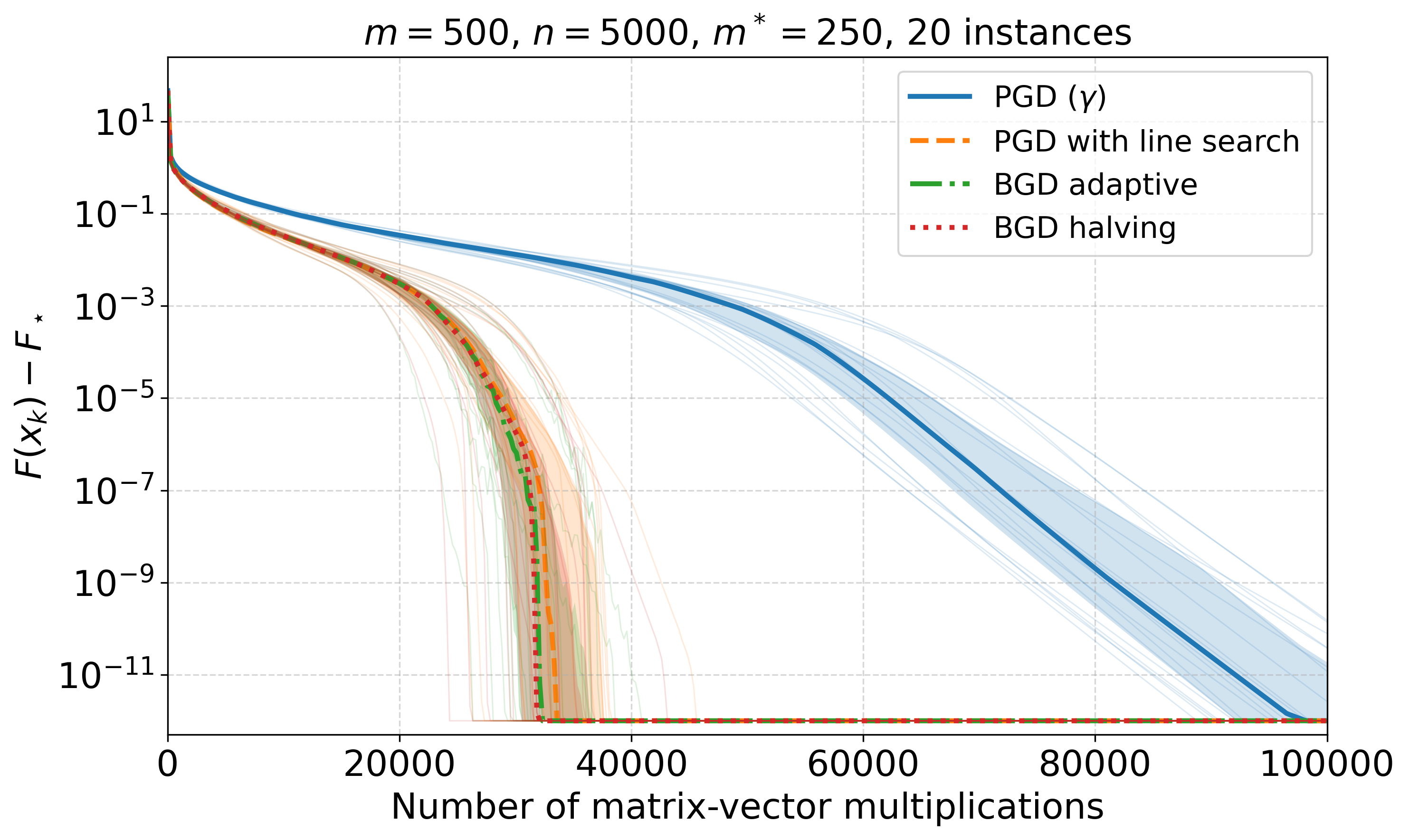}
	\caption{LASSO, float64. (a) $F(x_k) - F_\star$ vs.\ iteration $k$ on a single instance, for PGD with constant stepsize $\gamma$ and with line search, and BGD with constant radius $t$ and with decaying radii $\nicefrac{a}{\sqrt{k+2}}$ and $\nicefrac{a}{\sqrt{(k+2)\ln(k+2)}}$. (b)--(d) $F(x_k) - F_\star$ vs.\ number of matrix-vector multiplications over $20$ instances (median and $25\%$--$75\%$ range), for PGD with constant stepsize and with line search, and BGD with the adaptive radii \eqref{eq:adaptive} and with the halving radii \eqref{eq:halving}. PGD with constant stepsize diverges on $7$ of the $20$ instances in (b), which is why the upper quartile of its band lies above the plotted range.}
	\label{fig:convex-nofax}
\end{figure}

Figure~\ref{fig:convex-nofax}(a) shows a single instance, with parameters selected on it and target $\epsilon=10^{-6}$. PGD with line search is the fastest; BGD with a decaying radius needs about as many iterations as PGD with constant stepsize, but its iterates oscillate near the end, and with a constant radius they stall at $F(x_k)-F_\star\approx10^{-5}$. This reflects the trade-off of Corollary~\ref{co:convex-f-prac}: a higher precision requires a smaller radius, while too small a radius slows convergence. It motivates the adaptive rule \eqref{eq:adaptive} of Section~\ref{sec:convex-smooth}, for which we take $T = 10$ and $l_t=10^{-6}$, and select $t_0 \in 0.1 \times  \{ 1, 3^{-1}, \dots, 3^{-5} \}$ by grid search, as well as the halving rule \eqref{eq:halving}, for which we take $T = 10$ and select $t_0 \in \{1, 0.1\}$ and $c \in \frac{1}{2}\max_i\|a_i\|^2 \times \{3^{-1}, \dots, 3^{-5}\}$ by grid search.

Since the line search requires additional function evaluations, Figure~\ref{fig:convex-nofax}(b)--(d) measures the cost in matrix-vector multiplications, over $20$ instances per problem size, with $\epsilon=10^{-12}$ and at most $50000$ iterations. BGD with the adaptive stepsize requires fewer matrix-vector multiplications than PGD with line search in median for all three problem sizes, but its advantage is small: the median of the paired ratios is $1.09$, $1.03$ and $1.03$, and BGD is faster on $19$, $17$ and $18$ of the $20$ instances, respectively. Both are much faster than PGD with constant stepsize (paired ratio $2.1$ to $3.3$), which, with the stepsize selected on the tuning instance, diverges on $7$ instances in the first setting and does not reach $\epsilon$ within the budget on $8$ instances in the third. The selected $t_0$ is $0.1/3$ for all three sizes, and the halving safeguard is active, with $9$ to $10$ halvings per run in median and at most $13$. BGD with the halving radii \eqref{eq:halving} requires the fewest matrix-vector multiplications in median for all three sizes: the median of the paired ratios to BGD with the adaptive stepsize is $0.89$, $0.99$ and $0.99$, and it is faster than PGD with line search on $20$, $19$ and $20$ of the $20$ instances (see Appendix~\ref{sec:addexp} for the selected parameters).

\subsection{Elastic net problem}

In the strongly convex case, we test on the elastic net problem \citep{zou2005regularization}:
\begin{equation}\label{p:elanet-1}
	\min_{x \in \mathbb{R}^n } F(x) \eqdef  \frac{1}{2}\|Ax-b\|^2 + \frac{\lambda_2}{2}\|x\|^2 + \|x\|_1, 
\end{equation}
with $\lambda_2=0.1$, where $A$, $b$ and the solution $x^*$ are generated similarly to \citet{nesterov2013gradient} (Appendix~\ref{sec:elanet}). For BGD, we use the geometric radius $t_k = q^k$ suggested by Theorem~\ref{th:strongly-convex-prac}, the adaptive rule \eqref{eq:adaptive}, and the halving rule \eqref{eq:halving}, with grids $\gamma \in \{ 2^{-10}, \dots, 2^{-19} \}$, $L_0 \in \max_i\|a_i\|^2 \times \{ 1, 2^{-1}, \dots, 2^{-9} \}$, $q \in \{  \sqrt{1-2^{-5}}, \dots, \sqrt{1-2^{-14}}  \}$, $t_0$ and $T$ as for the LASSO, and $c \in c_0 \times \{3^{-1}, \dots, 3^{-5}\}$ with $c_0 \eqdef \frac{1}{2}(\max_i\|a_i\|^2 + \lambda_2)$; the target is $\epsilon=10^{-12}$, with at most $20000$ iterations for the first two problem sizes and $50000$ for the third.

\begin{figure}[t!]
	\centering
	
	\includegraphics[width=0.3\linewidth]{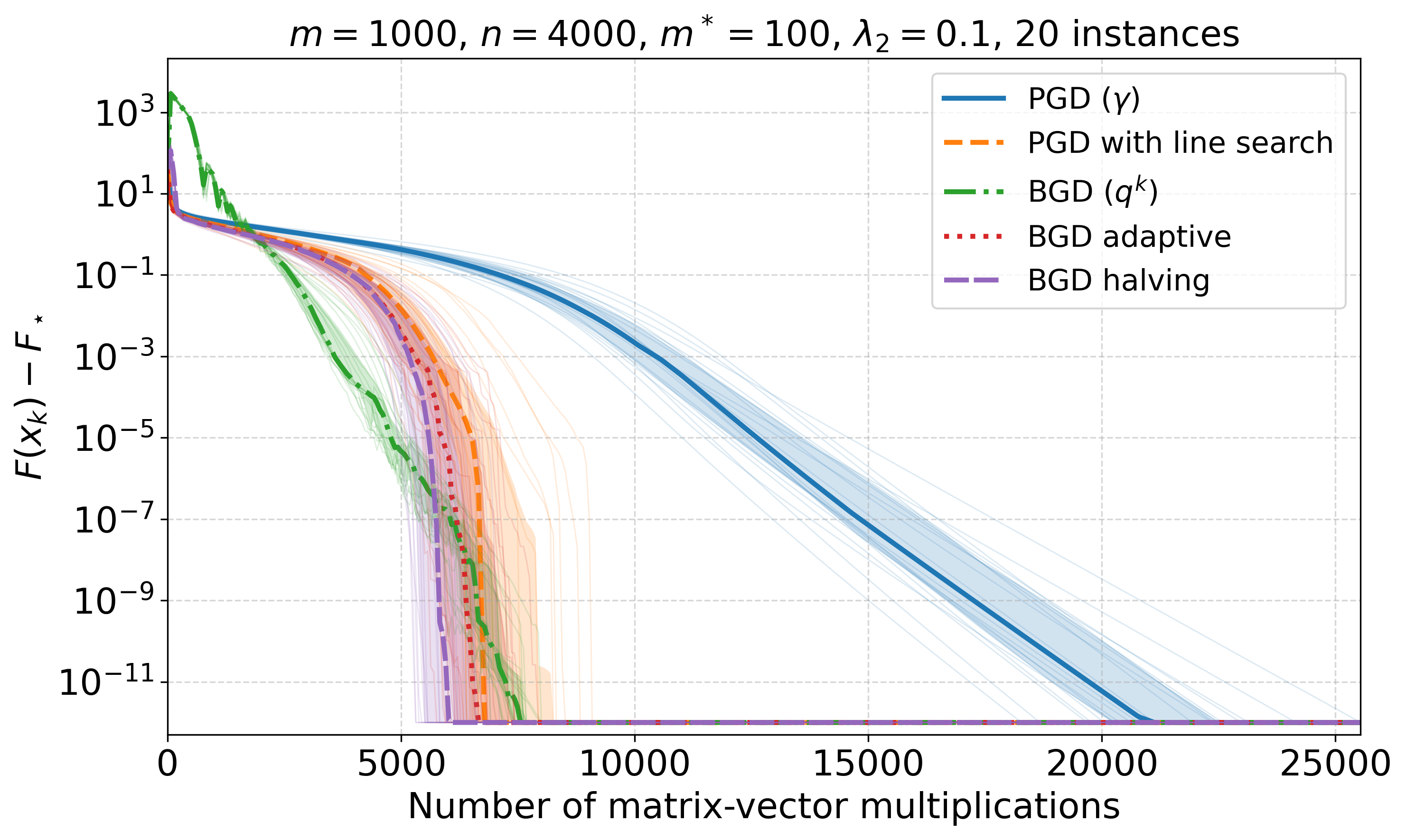}
	\hfill 
	\includegraphics[width=0.3\linewidth]{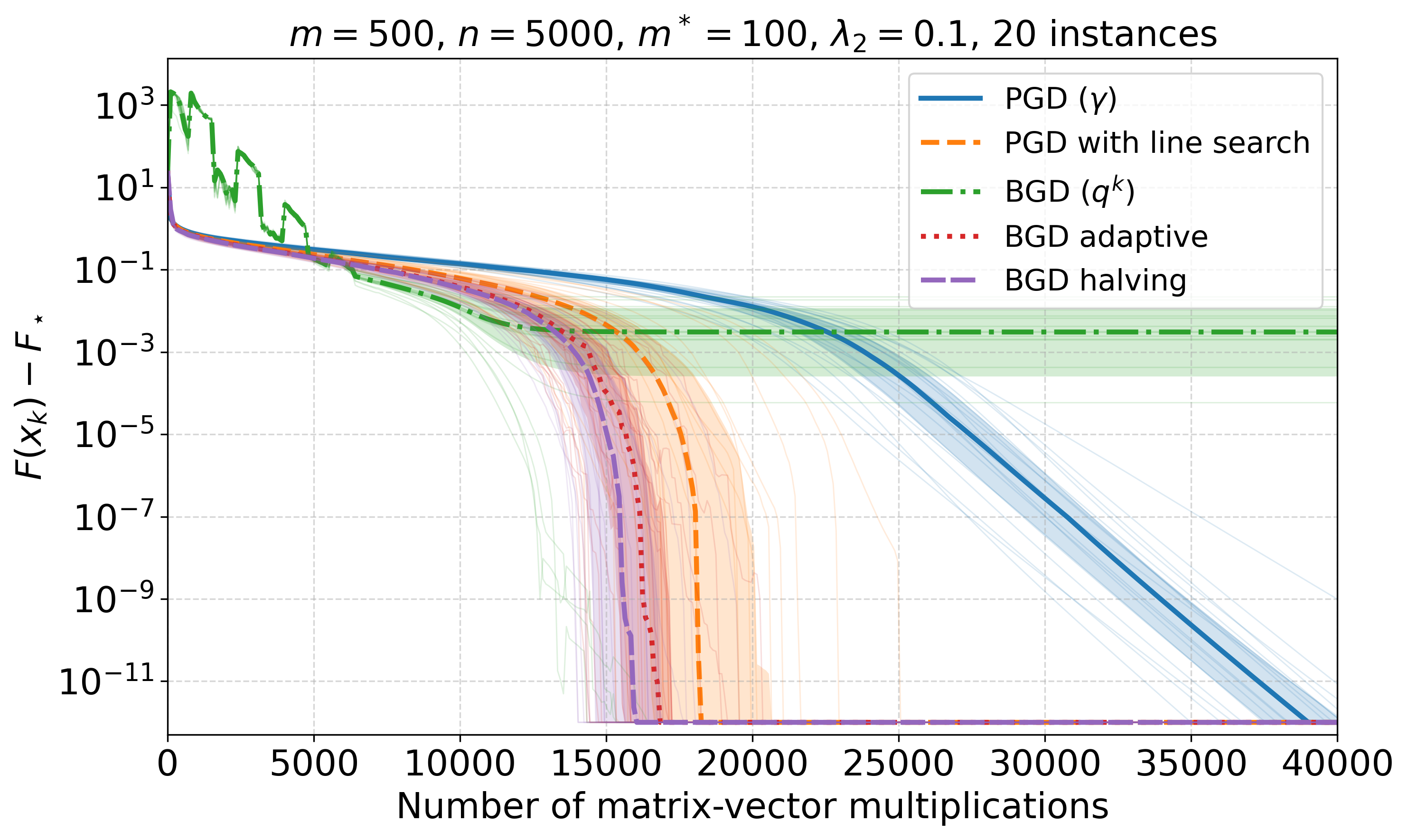}
	\hfill
	\includegraphics[width=0.3\linewidth]{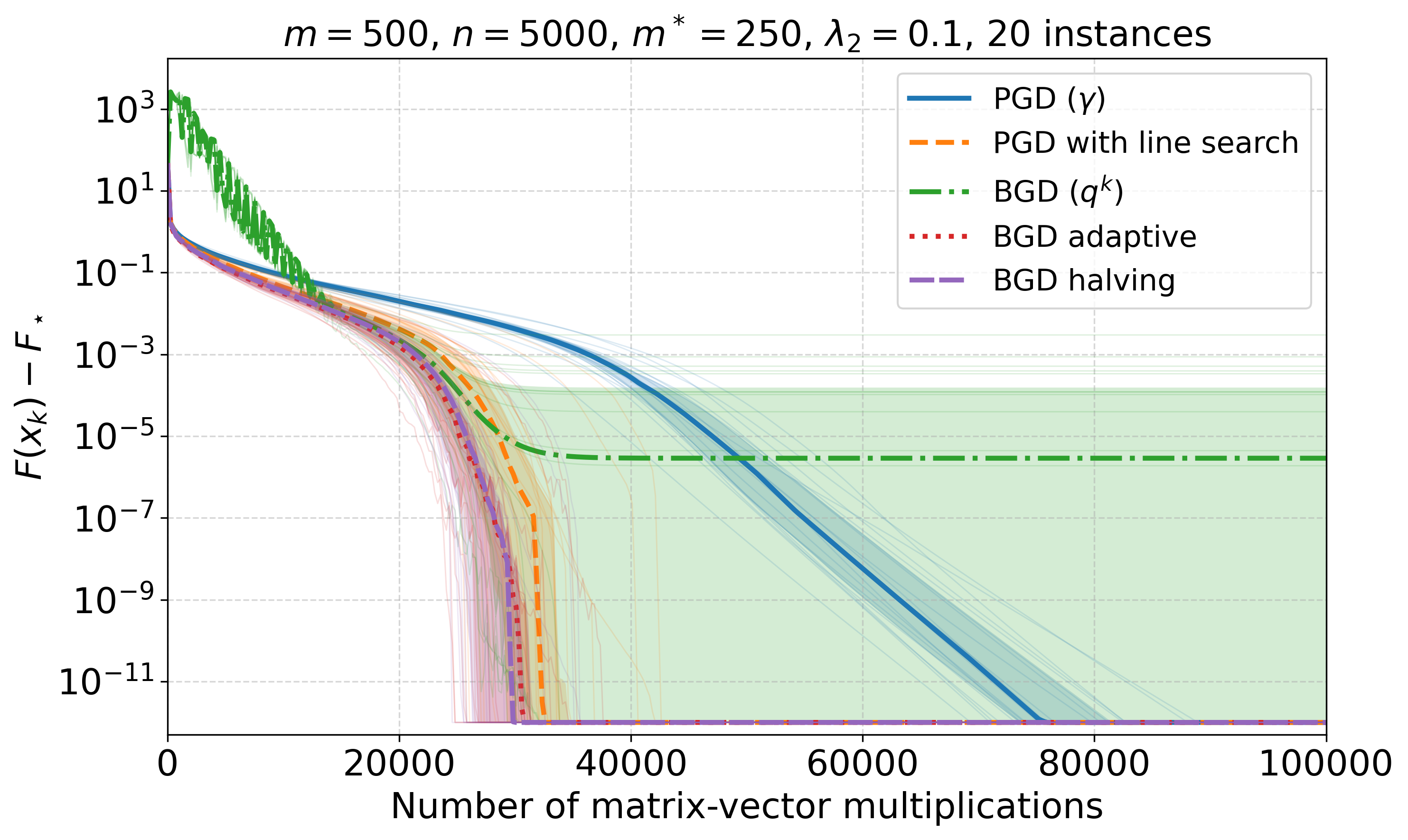}
	\caption{$F(x_k) - F_\star$ vs. number of matrix-vector multiplications for elastic net with $\lambda_2 = 0.1$, over $20$ instances (median and $25\%$--$75\%$ range). Algorithms: PGD with constant stepsize $\gamma$ and with line search, BGD with stepsize $t_k=q^k$, with the adaptive radii \eqref{eq:adaptive}, and with the halving radii \eqref{eq:halving}. Float64 is used.}
	\label{fig:sc-nofax-01}
\end{figure}

Figure~\ref{fig:sc-nofax-01} shows the results. BGD with the adaptive rule requires fewer matrix-vector multiplications than PGD in median for all three problem sizes and reaches $\epsilon$ on every instance; compared with PGD with line search, the median of the paired ratios is $1.07$, $1.13$ and $1.09$, and BGD is faster on $17$, $19$ and $19$ of the $20$ instances. PGD with constant stepsize is slower by a factor of $2.5$ to $3.2$ on every instance, and reaches $\epsilon$ within the budget on only $13$ instances in the second setting. The geometric schedule behaves very differently: on the first problem size it is competitive (paired ratio $1.12$, all instances solved), but with the factor $q$ selected on the tuning instance it stalls, at $F(x_k)-F_\star$ between $2\cdot10^{-6}$ and $2\cdot10^{-2}$, on $16$ and $11$ of the $20$ instances in the other two settings. Once $q^k$ becomes small compared with the distance to the solution, the remaining radii can no longer bring the iterates to $x_\star$; this is consistent with Theorem~\ref{th:strongly-convex-prac}, which requires $q$ to be close enough to $1$ in terms of the condition number, and the latter varies across instances. The adaptive rule, which does not prescribe the decay rate in advance, does not suffer from this, and neither does the halving rule, in line with Theorem~\ref{th:sc-halving}. BGD with the halving rule reaches $\epsilon$ on every instance and requires the fewest matrix-vector multiplications in median for all three sizes: the median of the paired ratios to BGD with the adaptive rule is $0.90$, $0.96$ and $0.98$, and it is faster than PGD with line search on $20$, $19$ and $20$ of the $20$ instances (Appendix~\ref{sec:addexp}).

\subsection{MLP on MNIST}

For the nonconvex case, we train a one-hidden-layer network on MNIST \citep{lecun1998gradient} under a box constraint, which confines the iterates to a compact set on which the loss, not globally $L$-smooth in general, is smooth:
\begin{equation}
	\min_{w\in\R^{d}}\; f(w)+g(w),
	\qquad
	f(w)=\frac1N\sum_{i=1}^{N}\ell\bigl(\phi(w;x_i),y_i\bigr),
	\qquad
	g=\ind_{[-r,r]^{d}},
	\label{p:mnist-problem}
\end{equation}
where $\phi(w;x)=W_2\tanh(W_1x+b_1)+b_2$ maps a flattened $784$-dimensional image to $10$ logits, $\ell$ is the softmax cross-entropy loss, $w=(W_1,b_1,W_2,b_2)$ collects all parameters ($d=101\,770$ for $128$ hidden units), and $\ind_C$ is the indicator function of $C$. From the official training split we draw two disjoint class-balanced subsets, $N=10000$ images for training and $10000$ for validation, with pixel values scaled to $[0,1]$. We take $r=0.5$, start all algorithms from the PyTorch default initialization \citep{paszke2019pytorch}, and use a budget of $B=5000$ forward/backward passes, which have almost the same cost for this architecture (Appendix~\ref{app:mlp-cost}). We test BGD with the adaptive rule \eqref{eq:adaptive} ($T=10$, $l_t=10^{-6}$), PGD with constant stepsize $\gamma$, and PGD with line search, with grids $\gamma \in 10\times \{ 3^{-i}\}_{i=0}^9$, $L_0 \in 10^{-3} \times \{ 3^i \}_{i=0}^9$ and $t_0 \in 100\times \{ 3^{-i} \}_{i=0}^9$, selecting the smallest training loss along the trajectory, which gives $\gamma=1.111$, $L_0=0.081$ and $t_0=11.11$. Selecting by the smallest validation loss or validation error yields the same $\gamma$ and $t_0$; for the line search, the validation loss selects a different $L_0$, which changes the median training loss below by about $4\%$ and the median test error by $0.04$ percentage points. We also test BGD with the halving rule \eqref{eq:halving}, with $T=10$, $t_0 \in \{10, 1\}$ and $c \in \{3^{-1}, \dots, 3^{-5}\}$, selected in the same way (Appendix~\ref{sec:appmlp}).

\begin{figure}[t!]
	\centering
	\includegraphics[width=0.32\linewidth]{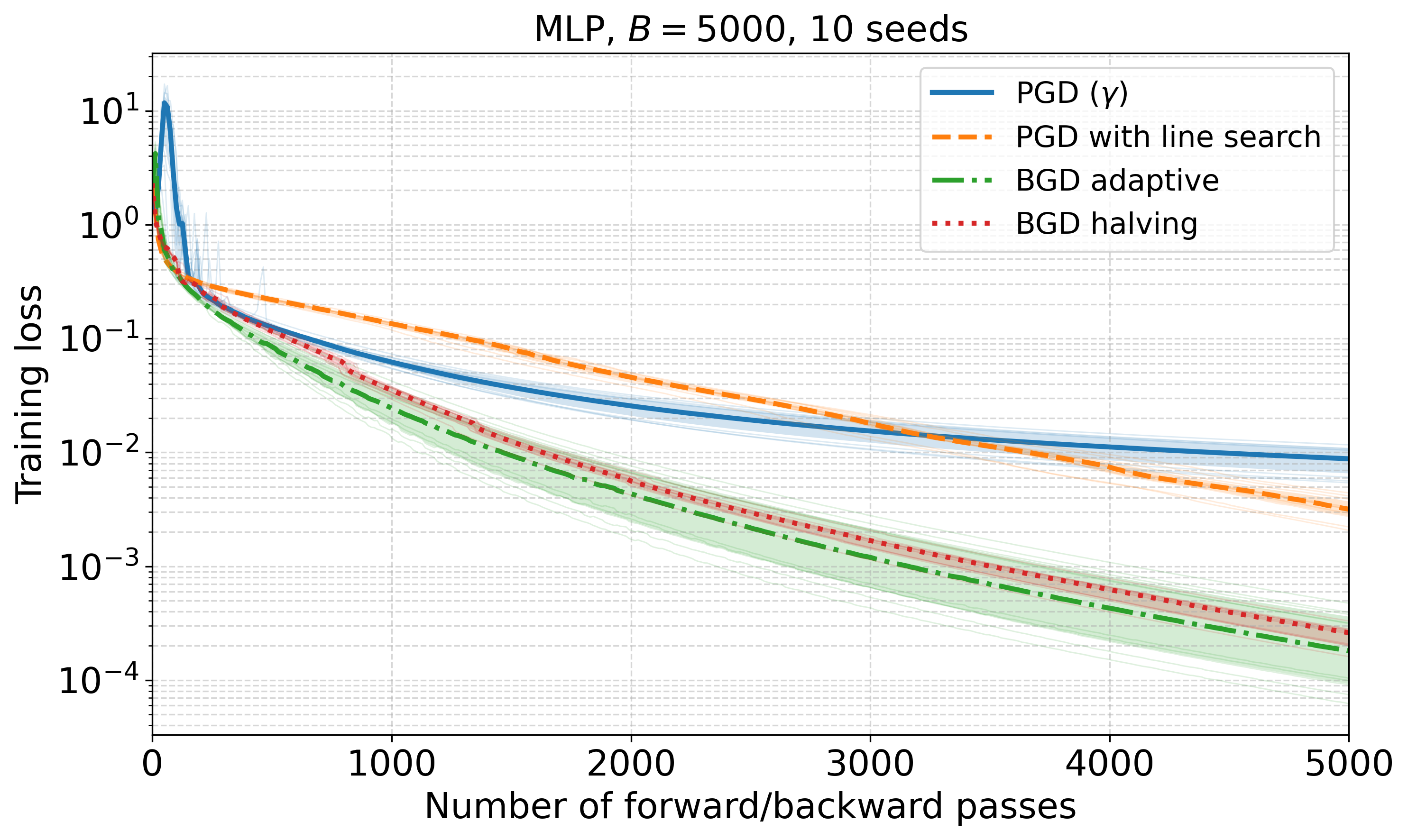}
	\hfill
	\includegraphics[width=0.32\linewidth]{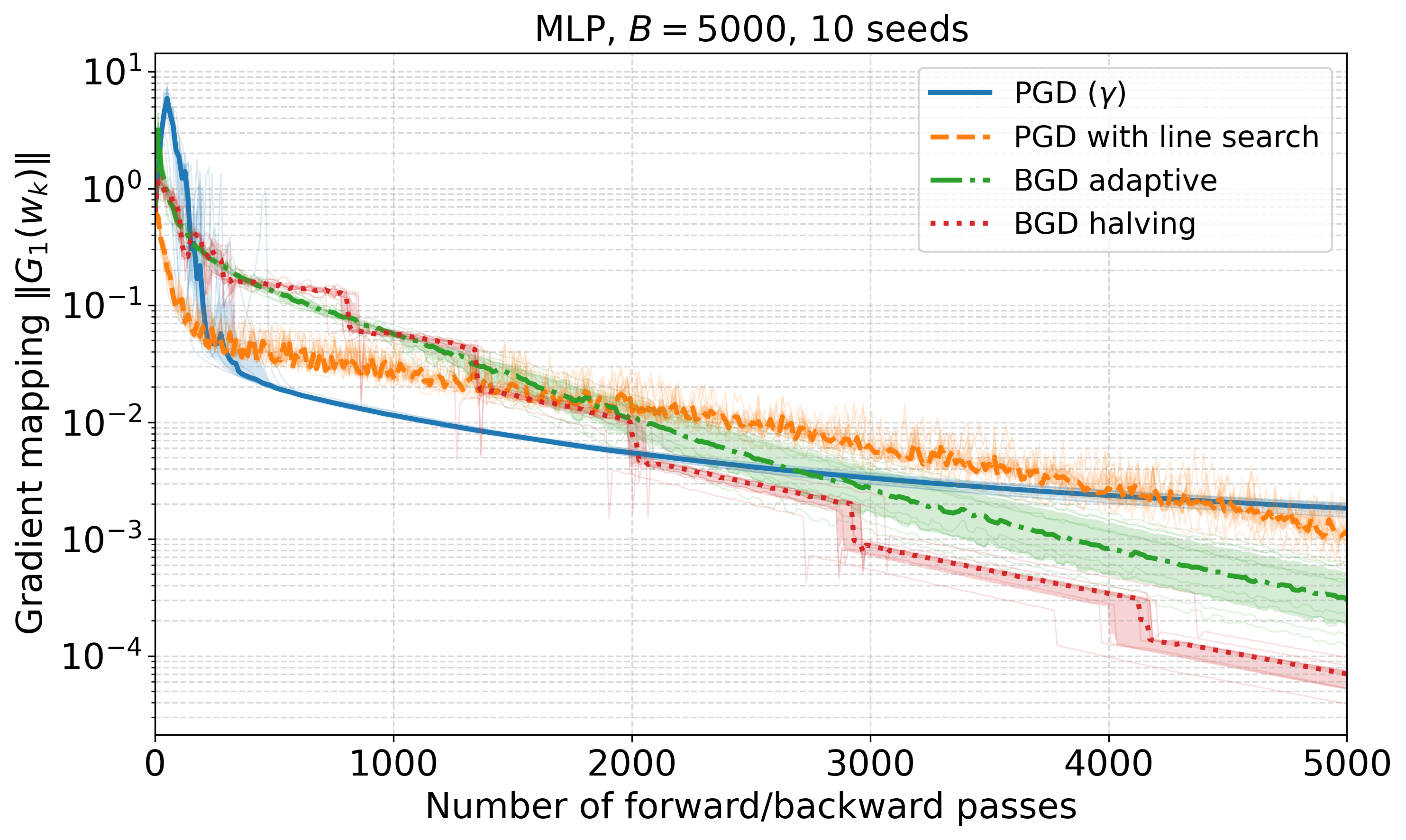}
	\hfill
	\includegraphics[width=0.32\linewidth]{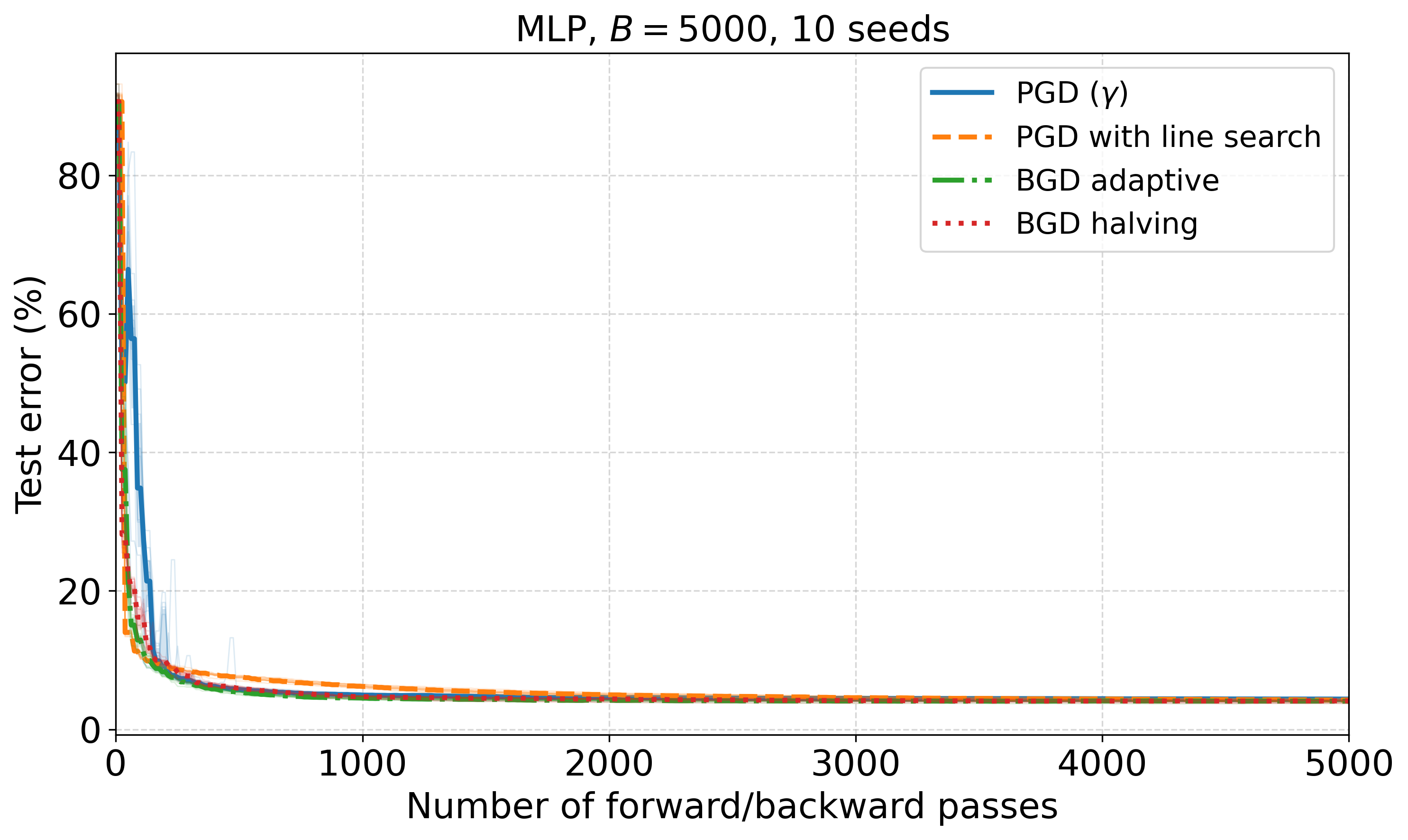}
	\caption{Training loss (left), gradient mapping $\|G_1(w_k)\|$ (middle), and test error (right) vs. number of forward/backward passes, over $10$ seeds (median and $25\%$--$75\%$ range). Algorithms: PGD with constant stepsize $\gamma$ and with line search, BGD with the adaptive radii \eqref{eq:adaptive} and with the halving radii \eqref{eq:halving}. Float32 is used.}
	\label{fig:mnist}
\end{figure}

Figure~\ref{fig:mnist} shows the results. BGD with the adaptive radii attains a smaller training loss than both PGD variants on all $10$ seeds, with a median of $2.1\cdot10^{-4}$ at the end of the budget, against $3.2\cdot10^{-3}$ for PGD with line search and $8.8\cdot10^{-3}$ for PGD with constant stepsize. The gradient mapping $\|G_1(w_k)\|$, the stationarity measure of Theorem~\ref{th:non-convex-withoutPL}, decreases more slowly for this method in the first $2000$ passes, during which it makes fast progress in $f$, but is smaller than for both PGD variants from about $3000$ passes on; at the end of the budget it is smaller on all $10$ seeds, with a median of $3.2\cdot10^{-4}$ against $1.2\cdot10^{-3}$ and $1.8\cdot10^{-3}$. The faster decrease of the training loss is not obtained at the expense of generalization: the method has a lower test error than each of the two PGD variants on $9$ of the $10$ seeds, with a median of $4.04\%$ against $4.24\%$ and $4.37\%$, although these differences are small in absolute terms. Since the training loss is essentially monotone, its minimum along the trajectory, which is used for the selection, coincides with its value at the end of the budget, at which all reported quantities are measured; we make no claim about early-stopped generalization. BGD with the halving radii behaves similarly: its median training loss at the end of the budget, $2.6\cdot10^{-4}$, is of the same order as with the adaptive radii, each of the two being smaller on $5$ of the $10$ seeds; its gradient mapping is the smallest of all methods on every seed, with a median of $7.0\cdot10^{-5}$, in line with Theorem~\ref{th:nc-halving}; and its median test error, $4.11\%$, lies between that of BGD with the adaptive radii and those of the two PGD variants.

\section{Conclusion}\label{sec:conclusion}

We introduced Broximal Gradient Descent, which minimizes the first-order model of a composite objective over a Euclidean ball around the current iterate, so that the radius of the ball, rather than a stepsize, controls how far the method moves. Two master theorems describe a single step and the radii for which it is admissible; from them we obtained convergence rates in the convex, strongly convex, and nonconvex settings that match those of proximal gradient descent, together with implementable radius schedules and a halving rule that requires no problem constant. BGD reduces to BroxGD \citep{richtarik2026local} when $g$ is an indicator function; PGD is a special case of BGD, but not conversely. In our experiments, at matched numbers of matrix-vector products or network passes, the adaptive and the halving radius rules made BGD competitive with PGD with line search: it was faster on most instances of the two regression problems, and reached a training loss an order of magnitude smaller on the network; each BGD step additionally requires a scalar search with about a dozen proximal evaluations. Natural directions for future work include accelerated and stochastic variants, and non-Euclidean balls.


	\bibliographystyle{plainnat} 
	\bibliography{ref.bib}

	\appendix

	\section{Solving the BGD subproblem}
\label{sec:subproblem}

Algorithm~\ref{alg:BG} requires an exact minimizer of the linearized
model of \eqref{p:p} over a ball,
\begin{equation}
	x_{k+1}\in\argmin_{z\in\cB(x_k,t_k)}m_k(z)
	\eqdef\langle\nabla f(x_k),z\rangle+g(z).
	\label{eq:sub-general}
\end{equation}
In all our experiments $g$ is separable and \eqref{eq:sub-general} is solved by a one-dimensional search whose relative accuracy on the radius is guaranteed. We first describe the procedure for
a general proper, closed, and convex $g$, and then specialize it to $g=\|\cdot\|_1$, used
for LASSO~\eqref{p:lasso} and for the elastic net~\eqref{p:elanet-1}, and to
$g=\ind_{[-r,r]^{d}}$, used for the box-constrained network
\eqref{p:mnist-problem}.

\subsection{The general procedure}
\label{sec:sub-general}

Assume $t_k>0$. As shown in the proof of Theorem~\ref{th:another-ad}, a minimizer $x_{k+1}$ of \eqref{eq:sub-general}
exists, and it is characterized by the KKT conditions \eqref{eq:kkt}: there is a
multiplier $\lambda_k\ge0$ with
\begin{equation*}
	0\in\nabla f(x_k)+\partial g(x_{k+1})+\lambda_k\,(x_{k+1}-x_k),
	\quad
	\|x_{k+1}-x_k\|\le t_k,
	\quad
	\lambda_k\bigl(\|x_{k+1}-x_k\|-t_k\bigr)=0 .
\end{equation*}
The solution procedure distinguishes two cases according to the value of the
term $\lambda_k(x_{k+1}-x_k)$.

\paragraph{Case 1: $\lambda_k(x_{k+1}-x_k)=0$.}
Then \eqref{eq:kkt} reduces to $0\in\nabla f(x_k)+\partial g(x_{k+1})$,
i.e.\ $x_{k+1}$ minimizes $m_k$ over the whole of $\R^{d}$, the ball playing no
role. This case therefore occurs exactly when $m_k$ attains its infimum at a
point of $\cB(x_k,t_k)$. Accordingly, we first solve
$\min_{z\in\R^{d}}m_k(z)$. If its infimum is $-\infty$ we move to Case 2.
Otherwise, denote by $Z^{0}$ the set of minimizers, which is closed and
convex. If $Z^{0}=\emptyset$, i.e., the infimum is finite but not attained, then
no point of the ball minimizes $m_k$ over $\R^{d}$, and we move to Case 2. If
$Z^{0}\neq\emptyset$, we compute its element nearest to $x_k$,
\begin{equation}
	z^{0}\eqdef\argmin_{z\in Z^{0}}\|z-x_k\| ,
	\label{eq:sub-nearest}
\end{equation}
and accept $x_{k+1}=z^{0}$ whenever $\|z^{0}-x_k\|\le t_k$; if
$\|z^{0}-x_k\|>t_k$, no minimizer of $m_k$ lies in the ball and we move to
Case~2. Deciding whether the infimum is finite and attained, and computing
$z^{0}$, in general requires information on $g$ beyond its proximal operator;
for the two choices of $g$ considered below, both are available in closed form.

\paragraph{Case 2: $\lambda_k(x_{k+1}-x_k)\neq0$.}
Then $\lambda_k>0$, so complementary slackness forces $\|x_{k+1}-x_k\|=t_k$:
the solution lies on the boundary of the ball. Writing
$\gamma\eqdef1/\lambda_k>0$, the stationarity condition in \eqref{eq:kkt}
becomes
\begin{equation*}
	0\in\partial g(x_{k+1})
	+\tfrac{1}{\gamma}\Bigl(x_{k+1}-\bigl(x_k-\gamma\nabla f(x_k)\bigr)\Bigr),
\end{equation*}
which is precisely the optimality condition of the proximal subproblem. Hence
\begin{equation}
	x_{k+1}=z(\gamma)\eqdef\prox_{\gamma g}\bigl(x_k-\gamma\nabla f(x_k)\bigr),
	\label{eq:sub-zgamma}
\end{equation}
that is, \emph{a proximal gradient step from $x_k$ with stepsize $\gamma$}, and
$\gamma$ is determined by the boundary condition
\begin{equation}
	\varphi(\gamma)\eqdef\bigl\|\prox_{\gamma g}
	\bigl(x_k-\gamma\nabla f(x_k)\bigr)-x_k\bigr\|=t_k .
	\label{eq:sub-phi}
\end{equation}
In this case a BGD step and a proximal gradient step produce the same point;
what distinguishes them is that here the stepsize is not prescribed but
determined implicitly by the radius $t_k$, as the reciprocal of the multiplier
$\lambda_k$.

\paragraph{Solving \eqref{eq:sub-phi}.}
The search for $\gamma$ relies on two monotonicity properties of $\varphi$.
First, $\varphi$ is continuous and nondecreasing \citep[see, e.g.,][]{beck2017first}.
Second, by Theorem~\ref{th:recoverPGD}(ii), $\varphi(\alpha\gamma)\ge\alpha\varphi(\gamma)$
for all $\gamma>0$ and $\alpha\in(0,1)$, i.e., $\gamma\mapsto\varphi(\gamma)/\gamma$ is
nonincreasing; equivalently, $\ln\varphi$ is a nondecreasing and $1$-Lipschitz
function of $\ln\gamma$.

The search starts from a bracket $\gamma_{\mathrm{lo}}<\gamma_{\mathrm{hi}}$ with
$\varphi(\gamma_{\mathrm{lo}})\le t_k<\varphi(\gamma_{\mathrm{hi}})$. Since
$\varphi(1/\lambda_k)=t_k$, such a bracket exists: an upper endpoint is found by
doubling $\gamma$ until $\varphi(\gamma)>t_k$, and a lower endpoint is available in
closed form, since nonexpansiveness of the proximal operator gives
\begin{equation}
	\varphi(\gamma)\le\gamma\|\nabla f(x_k)\|
	+\bigl\|\prox_{\gamma g}(x_k)-x_k\bigr\| ,
	\label{eq:sub-bracket}
\end{equation}
whose right-hand side is $O(\gamma)$ for both of our choices of $g$; requiring
it to be at most $t_k$ yields an explicit $\gamma$ with $\varphi(\gamma)\le t_k$.
Each iteration then evaluates $\varphi$ at the geometric mean
$\gamma=\sqrt{\gamma_{\mathrm{lo}}\gamma_{\mathrm{hi}}}$, stops if
$(1-\varepsilon_{\mathrm{rel}})t_k\le\varphi(\gamma)\le t_k$, and otherwise
replaces $\gamma_{\mathrm{lo}}$ by $\gamma$ if $\varphi(\gamma)<(1-\varepsilon_{\mathrm{rel}})t_k$
and $\gamma_{\mathrm{hi}}$ by $\gamma$ if $\varphi(\gamma)>t_k$. The search returns
$z(\gamma)$ for the evaluated $\gamma$ with $\varphi(\gamma)\le t_k$ and the largest
$\varphi(\gamma)$, so the point it produces is feasible for \eqref{eq:sub-general}.

The accuracy $\varepsilon_{\mathrm{rel}}$ is guaranteed, not only requested.
Indeed, \eqref{eq:sub-phi} has a solution $\gamma^\star\in[\gamma_{\mathrm{lo}},\gamma_{\mathrm{hi}}]$
throughout, and the Lipschitz property of $\ln\varphi$ gives
$\ln t_k-\ln\varphi(\gamma_{\mathrm{lo}})\le\ln(\gamma_{\mathrm{hi}}/\gamma_{\mathrm{lo}})$,
so the lower endpoint satisfies the accuracy requirement as soon as
$\ln(\gamma_{\mathrm{hi}}/\gamma_{\mathrm{lo}})\le\ln\frac{1}{1-\varepsilon_{\mathrm{rel}}}$.
Since each iteration halves $\ln(\gamma_{\mathrm{hi}}/\gamma_{\mathrm{lo}})$, the search
terminates with $(1-\varepsilon_{\mathrm{rel}})t_k\le\varphi(\gamma)\le t_k$ after at most
\begin{equation}
	N_\varphi\eqdef \left\lceil\log_2\frac{\ln(\gamma_{\mathrm{hi}}^{0}/\gamma_{\mathrm{lo}}^{0})}{\ln\frac{1}{1-\varepsilon_{\mathrm{rel}}}}\right\rceil
	\label{eq:sub-N}
\end{equation}
evaluations of $\varphi$, where $\gamma_{\mathrm{lo}}^{0}$ and $\gamma_{\mathrm{hi}}^{0}$
denote the initial bracket. For $\varepsilon_{\mathrm{rel}}=10^{-2}$, this gives
$N_\varphi\le13$ whenever $\gamma_{\mathrm{hi}}^{0}/\gamma_{\mathrm{lo}}^{0}\le e^{82}$, and
$N_\varphi\le18$ for every bracket of positive double-precision numbers, whose ratio is
below $e^{1455}$. With the cap $T_m=20$ used below, the cap is therefore never
reached, and $(1-\varepsilon_{\mathrm{rel}})t_k\le\|x_{k+1}-x_k\|\le t_k$ at every
iteration. By Theorem~\ref{th:recoverPGD}(i), applied with stepsize $\gamma$, the
returned point is the exact solution of \eqref{eq:sub-general} with radius
$\|x_{k+1}-x_k\|$, so that each step of the implementation is an exact BGD step with a
radius in $[(1-\varepsilon_{\mathrm{rel}})t_k,t_k]$. The only exception occurs when
$t_k$ falls below the floating-point resolution, $t_k\le10^{-14}\max\{1,\|x_k\|\}$,
which happens for geometric schedules after many iterations; we then set $x_{k+1}=x_k$
and report such steps separately. In the implementation the search is carried out
on $\lambda=1/\gamma$, an equivalent, order-reversing reparametrization.

Each evaluation of $\varphi$ costs one proximal operator and one norm, both
$O(d)$ for the two separable regularizers considered here; no matrix--vector product
with $A$ and no pass through the network is involved, against $O(md)$ for a
matrix--vector product with $A$ and far more for a pass through the network. We
therefore measure the budget in matrix--vector products or network passes, and
report the number of proximal evaluations per iteration separately.

\subsection{Implementation for the two regularizers}
\label{sec:sub-l1}\label{sec:sub-box}

\paragraph{$g=\|\cdot\|_1$ (LASSO and elastic net).} Here $\prox_{\gamma g}=\mathcal{S}_\gamma$ is the componentwise soft-thresholding, $\mathcal{S}_\gamma(u)_i=\mathrm{sign}(u_i)\max(|u_i|-\gamma,0)$, and $z(\gamma)=\mathcal{S}_\gamma(x_k-\gamma\nabla f(x_k))$. For Case 1, $m_k(z)=\sum_{i}(v_iz_i+|z_i|)$ with $v\eqdef\nabla f(x_k)$ is bounded below if and only if $\|v\|_\infty\le1$; if not, the ball is active and Case 2 is used directly. Otherwise the nearest minimizer \eqref{eq:sub-nearest} is assembled coordinatewise: $z^0_i=0$ if $|v_i|<1$, $z^0_i=\min(x_{k,i},0)$ if $v_i=1$, and $z^0_i=\max(x_{k,i},0)$ if $v_i=-1$, the projections of $x_{k,i}$ onto the half-lines on which $m_k$ is flat (the equalities $|v_i|=1$ are tested up to $10^{-12}$). For Case 2, soft-thresholding moves each coordinate by at most $\gamma$, so \eqref{eq:sub-bracket} gives $\varphi(\gamma)\le\gamma(\|\nabla f(x_k)\|+\sqrt{d})$, and the bracket starts at $\gamma_{\mathrm{lo}}=t_k/(\|\nabla f(x_k)\|+\sqrt{d})$.

\paragraph{$g=\ind_{[-r,r]^{d}}$ (box-constrained network).} Here $\prox_{\gamma g}=\Proj$ is the componentwise clipping to $[-r,r]$, and $z(\gamma)=\Proj(x_k-\gamma\nabla f(x_k))$ is a projected gradient step. For Case 1, the box being bounded, $m_k$ always attains its infimum, and the nearest minimizer is the corner $z^0_i=-r\,\mathrm{sign}(\nabla_if(x_k))$ if $\nabla_if(x_k)\ne0$ and $z^0_i=x_{k,i}$ otherwise (tested up to $10^{-12}$). For Case 2, the iterates are feasible, so \eqref{eq:sub-bracket} gives $\varphi(\gamma)\le\gamma\|\nabla f(x_k)\|$, and the bracket starts at $\gamma_{\mathrm{lo}}=\min\{1,t_k/(1.001\|\nabla f(x_k)\|)\}$.

In both cases $\gamma_{\mathrm{lo}}$ is halved while $\varphi(\gamma_{\mathrm{lo}})>t_k$, as a safeguard against rounding, $\gamma_{\mathrm{hi}}$ starts at $10^{12}$ and is doubled while $\varphi(\gamma_{\mathrm{hi}})\le t_k$, and we use $T_m=20$ and $\varepsilon_{\mathrm{rel}}=10^{-2}$.

\section{Paired comparisons over instances}\label{sec:addexp}

For each problem size and each method, the tables below report the number of matrix--vector products needed to reach $F(x_k)-F_\star\le\epsilon$ (median and quartiles over the instances), the median over instances of its ratio to that of BGD with the adaptive radii \eqref{eq:adaptive} on the same instance, the number of instances on which the latter is faster, the number of instances on which the accuracy is not reached within the budget, and the median number of proximal evaluations per iteration. For the network, the corresponding quantities are the training loss and the test error at the budget. For BGD with the halving radii \eqref{eq:halving}, with $c$ expressed in units of $c_0 \eqdef \frac{1}{2}\max_i\|a_i\|^2$ for the LASSO and $c_0 \eqdef \frac{1}{2}(\max_i\|a_i\|^2 + \lambda_2)$ for the elastic net, the selected $(t_0, \nicefrac{c}{c_0})$ are $(1, 0.012)$, $(0.1, 0.037)$ and $(0.1, 0.037)$ for the three LASSO sizes, and $(1, 0.012)$, $(0.1, 0.012)$ and $(0.1, 0.004)$ for the elastic net, the last value being the smallest of the grid; the radius is halved $22$, $18$ and $23$ times per run in median for the LASSO, and $22$, $18$ and $20$ times for the elastic net.

In all LASSO, elastic net and network runs, including the whole tuning grids, the search of Appendix~\ref{sec:subproblem} never reached the cap $T_m=20$: it used about $10$ to $11$ evaluations of $\varphi$ per iteration on average and at most $13$ (LASSO and elastic net) and $12$ (network), in each case equal to the bound \eqref{eq:sub-N} for the brackets encountered, and every computed radius was within $1\%$ of the prescribed one. The only exceptions are the steps of the geometric schedule whose radius has fallen below the floating-point resolution (Appendix~\ref{sec:subproblem}), at which the iterate is kept unchanged; they occur for the elastic net, after the iterates have stalled, for the factors $q$ that are too small for the instance at hand.

\begin{table}[h!]
	\centering
	\caption{MNIST network, $B=5000$, parameters selected by the training loss: values at the end of the budget over $10$ paired seeds (median and quartiles), and the number of seeds on which BGD with the adaptive radii attains a smaller value.}
	\label{tab:mnist-paired}
	\resizebox{\textwidth}{!}{%
	\begin{tabular}{lcccccc}
		\toprule
		Method & Training loss & \makecell{BGD adaptive\\lower} & Test error (\%) & \makecell{BGD adaptive\\lower} & $\|G_1(w_K)\|$ & Prox per iter. \\
		\midrule
		PGD ($\gamma$) & $8.8\cdot10^{-3}$ [$6.6\cdot10^{-3}$, $1.1\cdot10^{-2}$] & 10/10 & $4.37$ [$4.10$, $4.51$] & 9/10 & $1.8\cdot10^{-3}$ & $1.0$ \\
		PGD with line search & $3.2\cdot10^{-3}$ [$2.8\cdot10^{-3}$, $3.7\cdot10^{-3}$] & 10/10 & $4.24$ [$4.17$, $4.34$] & 9/10 & $1.2\cdot10^{-3}$ & $2.0$ \\
		BGD adaptive & $2.1\cdot10^{-4}$ [$9.3\cdot10^{-5}$, $3.6\cdot10^{-4}$] & -- & $4.04$ [$3.94$, $4.06$] & -- & $3.2\cdot10^{-4}$ & $12.4$ \\
		BGD halving & $2.6\cdot10^{-4}$ [$2.0\cdot10^{-4}$, $2.8\cdot10^{-4}$] & 5/10 & $4.11$ [$3.95$, $4.21$] & 7/10 & $7.0\cdot10^{-5}$ & $12.4$ \\
		\bottomrule
	\end{tabular}}
\end{table}

\begin{table}[h!]
	\centering
	\caption{LASSO (Figure~\ref{fig:convex-nofax}) and elastic net (Figure~\ref{fig:sc-nofax-01}): number of matrix--vector products to reach $F(x_k)-F_\star\le10^{-12}$ over $20$ paired instances (median and quartiles; ``$>$budget'' if not reached within the iteration budget, which is $50000$ except for the first two elastic net sizes, where it is $20000$), median over instances of its ratio to that of BGD with the adaptive stepsize on the same instance, number of instances on which the latter is faster, number of instances on which $\epsilon$ is not reached, and median number of proximal evaluations per iteration.}
	\label{tab:convex-paired}
	\resizebox{\textwidth}{!}{%
	\begin{tabular}{lccccc}
		\toprule
		Method & \makecell{Matrix--vector products\\to reach $\epsilon$} & \makecell{Ratio to\\BGD adaptive} & \makecell{BGD adaptive\\faster} & \makecell{Not\\reached} & \makecell{Prox\\per iter.} \\
		\midrule
		\multicolumn{6}{l}{LASSO, $m=1000$, $n=4000$, $m^*=100$} \\
		\quad PGD ($\gamma$) & 13462 [12802, $>$budget] & 2.14 & 20/20 & 7 & 1.0 \\
		\quad PGD with line search & 7062 [6578, 7424] & 1.09 & 19/20 & 0 & 2.0 \\
		\quad BGD adaptive & 6462 [6134, 6904] & 1.00 & -- & 0 & 12.7 \\
		\quad BGD halving & 5782 [5497, 6022] & 0.89 & 0/20 & 0 & 12.6 \\
		\midrule
		\multicolumn{6}{l}{LASSO, $m=500$, $n=5000$, $m^*=100$} \\
		\quad PGD ($\gamma$) & 50220 [48755, 51452] & 2.81 & 20/20 & 0 & 1.0 \\
		\quad PGD with line search & 18380 [17136, 19086] & 1.03 & 17/20 & 0 & 2.0 \\
		\quad BGD adaptive & 17552 [17067, 18368] & 1.00 & -- & 0 & 12.8 \\
		\quad BGD halving & 17365 [16752, 18500] & 0.99 & 4/20 & 0 & 12.9 \\
		\midrule
		\multicolumn{6}{l}{LASSO, $m=500$, $n=5000$, $m^*=250$} \\
		\quad PGD ($\gamma$) & 97370 [92188, $>$budget] & 3.34 & 20/20 & 8 & 1.0 \\
		\quad PGD with line search & 33542 [32147, 37560] & 1.03 & 18/20 & 0 & 2.0 \\
		\quad BGD adaptive & 32192 [30287, 36256] & 1.00 & -- & 0 & 12.6 \\
		\quad BGD halving & 31849 [30266, 35344] & 0.99 & 6/20 & 0 & 13.0 \\
		\midrule
		\multicolumn{6}{l}{Elastic net, $\lambda_2=0.1$, $m=1000$, $n=4000$, $m^*=100$} \\
		\quad PGD ($\gamma$) & 20952 [20205, 22521] & 3.25 & 20/20 & 0 & 1.0 \\
		\quad PGD with line search & 6723 [6439, 7988] & 1.07 & 17/20 & 0 & 2.0 \\
		\quad BGD ($q^k$) & 7547 [7353, 7672] & 1.12 & 16/20 & 0 & 11.6 \\
		\quad BGD adaptive & 6633 [6127, 7287] & 1.00 & -- & 0 & 12.7 \\
		\quad BGD halving & 5912 [5492, 6602] & 0.90 & 0/20 & 0 & 12.5 \\
		\midrule
		\multicolumn{6}{l}{Elastic net, $\lambda_2=0.1$, $m=500$, $n=5000$, $m^*=100$} \\
		\quad PGD ($\gamma$) & 39000 [37524, $>$budget] & 2.50 & 20/20 & 7 & 1.0 \\
		\quad PGD with line search & 18164 [16494, 20220] & 1.13 & 19/20 & 0 & 2.0 \\
		\quad BGD ($q^k$) & $>$budget [$>$budget, $>$budget] & $>$budget & 19/20 & 16 & 10.8 \\
		\quad BGD adaptive & 16530 [15474, 17212] & 1.00 & -- & 0 & 12.7 \\
		\quad BGD halving & 15752 [14527, 16812] & 0.96 & 0/20 & 0 & 12.8 \\
		\midrule
		\multicolumn{6}{l}{Elastic net, $\lambda_2=0.1$, $m=500$, $n=5000$, $m^*=250$} \\
		\quad PGD ($\gamma$) & 75509 [73521, 81545] & 2.50 & 20/20 & 0 & 1.0 \\
		\quad PGD with line search & 32178 [28937, 33909] & 1.09 & 19/20 & 0 & 2.0 \\
		\quad BGD ($q^k$) & $>$budget [31100, $>$budget] & $>$budget & 20/20 & 11 & 8.6 \\
		\quad BGD adaptive & 30522 [26677, 31252] & 1.00 & -- & 0 & 12.6 \\
		\quad BGD halving & 29644 [26082, 30748] & 0.98 & 3/20 & 0 & 12.7 \\
		\bottomrule
	\end{tabular}}
\end{table}

\section{The elastic net problem}\label{sec:elanet}
In the strongly convex case, we test the algorithms on the following elastic net problem:
\begin{equation}\label{p:elanet}
	\min_{x \in \mathbb{R}^n } F(x) \eqdef  \frac{1}{2}\|Ax-b\|^2 + \frac{\lambda_2}{2}\|x\|^2 + \|x\|_1, 
\end{equation}
where $A \in \R^{m\times n}$, $b\in \R^m$, and $\lambda_2 \geq 0$. The data are generated as in the LASSO generator of \citet{nesterov2013gradient}, with one modification that accounts for the quadratic term. Given $m^* \leq m$, the number of nonzero components of the solution $x^*$, and $\rho > 0$, which controls its size, we draw a matrix $B \in \R^{m \times n}$ with entries uniformly distributed in $[-1, 1]$ and a vector $v^* \in \R^m$ with entries uniformly distributed in $[0, 1]$, set $y^* = v^*/\|v^*\|$, and order the columns $b_i$ of $B$ by decreasing $|\langle b_i, y^* \rangle|$. We then set $[x^*]^{(i)} = \xi_i \operatorname{sign}(\langle b_i, y^* \rangle)$ with $\xi_i$ uniformly distributed in $[0, \nicefrac{\rho}{\sqrt{m^*}}]$ for $i \leq m^*$ and $[x^*]^{(i)} = 0$ otherwise, define the columns of $A$ as $a_i = \alpha_i b_i$ with
\[
\alpha_i = 
\begin{cases} 
	\frac{\operatorname{sign}(\langle b_i, y^* \rangle) + \lambda_2 [x^*]^{(i)} }{\langle b_i, y^*\rangle}, & \text{for } i \leq m^*, \\[6pt]
	1, & \text{if } |\langle b_i, y^* \rangle| \leq 0.1 \text{ and } i > m^*, \\[4pt]
	\frac{\xi'_i}{|\langle b_i, y^*\rangle|}, & \text{otherwise},
\end{cases}
\]
where $\xi'_i$ are uniformly distributed in $[0, 1]$, and $b = y^* + A x^*$. Only the first case of $\alpha_i$ differs from the LASSO generator, to which the construction reduces for $\lambda_2 = 0$. By construction, $x^*$ is the solution of \eqref{p:elanet}, and
\[
F_\star = \frac{1}{2}\|y^*\|^2 + \frac{\lambda_2}{2}\|x^*\|^2 + \|x^*\|_1.
\]

\section{Details of the MNIST experiment}\label{sec:appmlp}

The radius $r=0.5$ was fixed before the experiments, at about six times the amplitude $1/\sqrt{128}\approx0.088$ of PyTorch's default initialization of the output layer (every entry of $W_\ell$ and $b_\ell$ is drawn uniformly from $[-1/\sqrt{m_\ell},1/\sqrt{m_\ell}]$, with fan-in $m_1=784$ and $m_2=128$, so that $w_0$ is strictly feasible); it is not tuned and is the same for every method, candidate and seed. The constraint becomes active only on the output layer: no coordinate touches the boundary during the first few hundred iterations, after which the fraction of entries of $W_2$ on the boundary grows to $59\%$ and $68\%$ for BGD with the adaptive and the halving radii, against $31\%$ and $45\%$ for PGD with constant stepsize and with line search, at $B=5000$, whereas at most $0.1\%$ of the entries of $W_1$, which holds $98.6\%$ of the coordinates, are active (Figure~\ref{fig:mnist-stepsize}, right). All runs use single precision on a CPU, PyTorch $2.1.2$.

\subsection{Cost model}
\label{app:mlp-cost}

We count the cost as the number of forward plus backward passes. Since the gradient with respect to the input is not needed, the two passes have almost the same number of floating-point operations for this architecture (ratio $1.01$; the measured time ratio is about $1.3$). The objective value used by the adaptive and the halving rules comes from the forward pass needed for the gradient, and, for the LASSO, from the residual $Ax_k-b$ needed for $\nabla f(x_k)$; hence BGD and PGD with constant stepsize both cost $2$ units per iteration, and the scalar search of BGD only adds proximal evaluations, which are reported separately. PGD with line search costs $j+1$ units at an iteration testing $j$ trial points, the forward pass of the accepted trial being reused.

\subsection{Algorithmic details}
\label{app:mlp-algorithms}

\paragraph{Line search.}
With $\gamma_u=\gamma_d=2$ and $L\leftarrow L_0$, the trial
$T_L=\Proj(w_k-\nabla f(w_k)/L)$ is accepted if
\begin{equation}
	f(T_L)\;\le\;f(w_k)+\langle\nabla f(w_k),T_L-w_k\rangle
	+\tfrac{L}{2}\|T_L-w_k\|^{2}+\mathrm{tol}_k ,
	\label{eq:ls-accept}
\end{equation}
and otherwise $L\leftarrow\gamma_uL$; upon acceptance $L_k:=L$ and
$L\leftarrow\max(L_0,L_k/\gamma_d)$. The safeguards $L\le10^{12}$ and at most
$100$ trials per iteration were never triggered. The tolerance
$\mathrm{tol}_k=10^{-6}\max(1,|f(w_k)|)$, where $10^{-6}\approx16u$ and $u=2^{-24}$ is the unit roundoff, absorbs the cancellation in the difference of the two sides of \eqref{eq:ls-accept} in single precision, which would otherwise cause spurious rejections; it is inactive in practice, since the median slack in \eqref{eq:ls-accept} over the accepted steps of the selected run is $-1.5\cdot10^{-5}$, an order of magnitude below $-\mathrm{tol}_k$.
The line search averages $1.8$ trial points per iteration on the tuning seed, and $1.95$ in median over the evaluation seeds.

\paragraph{BGD subproblem.}
Here $g$ is the indicator of $[-r,r]^{d}$, so the subproblem is the
minimization of a linear function over the intersection of a ball and a box;
it is solved, up to a relative tolerance of $10^{-2}$ on the radius, by
the procedure of Appendix~\ref{sec:sub-box}. Every step obtained in Case~2 of Appendix~\ref{sec:sub-general}, which is every step of our runs, is thus a projected gradient step
$w_{k+1}=\Proj\bigl(w_k-\gamma_k\nabla f(w_k)\bigr)$ whose stepsize $\gamma_k$
is determined implicitly by $t_k$ and equals $t_k/\|\nabla f(w_k)\|$ whenever no
coordinate is clipped along the step; Figure~\ref{fig:mnist-stepsize} plots this
quantity alongside the constant $\gamma$ of PGD and the $1/L_k$ of the line
search. For the selected $t_0$ the search of Appendix~\ref{sec:subproblem} used
$10.4$ evaluations of $\varphi$ on average over the evaluation seeds and at most
$12$, the theoretical bound \eqref{eq:sub-N} for these brackets; over all ten candidates and all evaluation seeds
the cap $T_m$ was never reached, every computed radius was within $1\%$ of the prescribed one, and
no corner step --- the case in which the whole box minimizer falls inside
the ball --- occurred in any run: the distance from $w_k$ to
that corner stays close to $150$ throughout, well above the largest radius in
the grid, $t_0=100$.

\paragraph{Radius schedule.}
The radius is given by \eqref{eq:adaptive} with $T=10$ and $l_t=10^{-6}$, the loss being compared with its value at the previous checkpoint, so that Theorem~\ref{th:non-convex-withoutPL} applies with $a_k=\nicefrac{t_k}{\gamma}$. The safeguard fires rarely for the selected $t_0=11.11$: once on the tuning seed and on $5$ of the $10$ evaluation seeds, and never on the other $5$. The coefficient $t_0 2^{-h_k}\geq 5.56$ therefore stays far above the lower bound $l_t$, and the radius essentially follows the prescribed decay $t_0 s_k$. The halving rule \eqref{eq:halving} is used with $T=10$ and the same subproblem solver. The training loss selects $(t_0, c) = (1, 3^{-5})$, where $3^{-5}$ is the smallest value of the grid, whereas the validation loss and the validation error select $(10, 3^{-1})$ and $(1, 3^{-4})$; for the selected parameters, the radius is halved $7$ times on the tuning seed and on each evaluation seed, the search used $10.4$ evaluations of $\varphi$ on average and at most $12$, the cap $T_m$ was never reached over the whole grid, and no corner step occurred.

\subsection{Additional results}
\label{app:mlp-figures}

\begin{figure}[t!]
	\centering
	\includegraphics[width=0.4\linewidth]{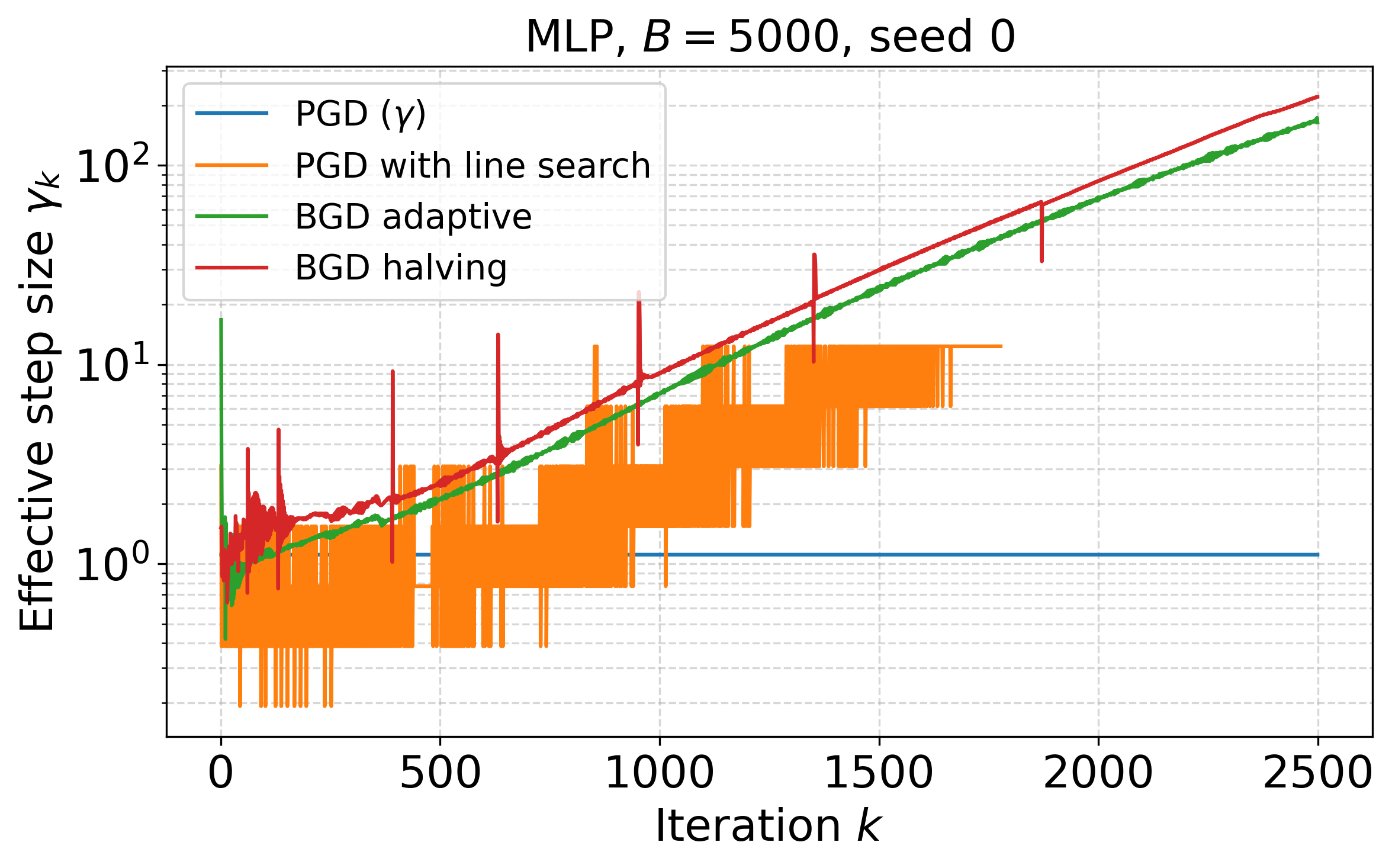}
	\hfill
	\includegraphics[width=0.56\linewidth]{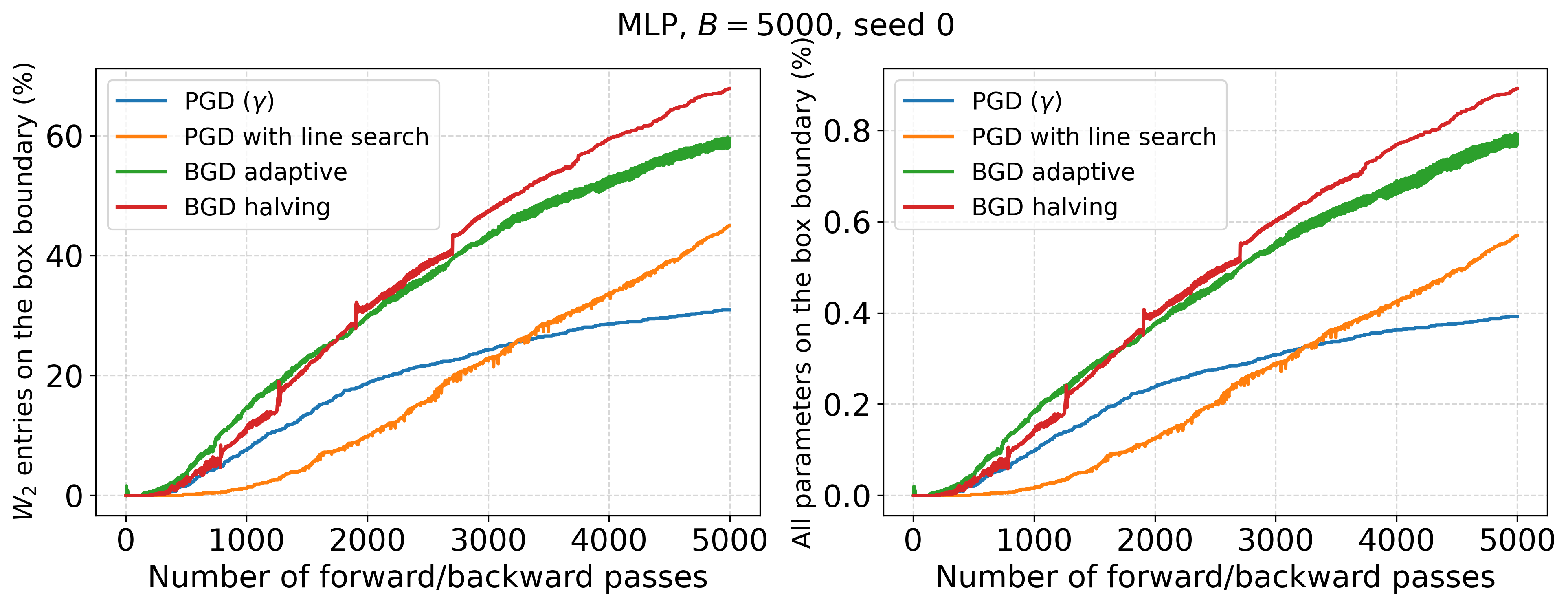}
	\caption{Tuning seed, parameters selected by the training loss, $B=5000$, for PGD with constant stepsize and with line search and for BGD with the adaptive and the halving radii. Left: effective stepsize $\gamma_k$ vs.\ iteration. Right: fraction of parameters on the boundary of the box, for the entries of $W_2$ and for all parameters.}
	\label{fig:mnist-stepsize}
\end{figure}

Figure~\ref{fig:mnist-stepsize} shows, for the tuning seed and the parameters selected by the training loss, the effective stepsize, defined as the scalar $\gamma_k$ for which $w_{k+1}=\Proj(w_k-\gamma_k\nabla f(w_k))$, together with the fraction of parameters on the boundary of the box. The scales are very different: PGD with line search stays within $[0.19,12.3]$ around the constant stepsize $\gamma=1.111$, whereas the effective stepsize of BGD with the adaptive radii, after a transient over the first few iterations, grows steadily from about $1$ to $170$ over the $2500$ iterations, while the radius decreases to $t_k=0.047$. With the halving radii, the effective stepsize grows similarly, from about $1.5$ to $220$, while the radius is halved $7$ times, from $t_0=1$ to $2^{-7}\approx0.0078$. This growth is how a method that takes a prescribed step length adapts automatically once the landscape flattens, and it is the mechanism behind the gap in Figure~\ref{fig:mnist}. The isolated spikes in the figure are caused by the halvings: before a halving, the iterates oscillate across a valley and the objective stops decreasing; the first step with the halved radius lands close to the bottom of the valley, where the gradient norm is several times smaller, so that the effective stepsize of the next step jumps up, and it returns to its trend once the gradient norm has settled.

\section{Proofs for the master theorems}

\subsection{Proof of Theorem \ref{th:another-ad}}

Fix \( k \geq 0 \), and for the purposes of this proof, let us use the shorter notation \( \mathcal{B}_k := \mathcal{B}(x_k, t_k) \). If $\|x_k-y\| \leq t_k$, then from the definition of $x_{k+1}$, $m_k(x_{k+1}) \leq m_k(y)$, which contradicts the assumption that $m_k(x_{k+1}) > m_k(y)$. Thus {\em (i)} holds. Let us now establish claims {\em (ii)}, {\em (iii)}, {\em (iv)}, and {\em (v)}.

\textbf{Step 1: KKT conditions for a convex optimization problem.}
The function $m_k + \iota_{\mathcal{B}_k}$ is proper, since
$x_k \in \dom g \cap \mathcal{B}_k \neq \emptyset$, and closed, being the sum of
a continuous linear term, of the closed function $g$, and of the indicator of
the closed set $\mathcal{B}_k$. Its domain is contained in the bounded set
$\mathcal{B}_k$, so all its sublevel sets are compact and $m_k(\cdot)$ attains
its minimum over $\mathcal{B}_k$; note that $\dom g \cap \mathcal{B}_k$ itself
need not be compact, as $\dom g$ need not be closed, and that convexity is not
needed for this argument. Let $x_{k+1}$ be a minimizer, not necessarily unique.
By Fermat's rule, $x_{k+1}$ minimizes $m_k(\cdot)$ over $\mathcal{B}_k$ if and
only if $0 \in \partial \left( m_k + \iota_{\mathcal{B}_k} \right)(x_{k+1})$.
Since $x_k \in \dom g \cap {\rm int}\, \mathcal{B}_k$, the Moreau--Rockafellar
sum rule applies \citep[Theorem~23.8]{rockafellar1997convex} and yields the
equivalent first-order optimality condition
\[
0 \in \nabla f(x_k) + \partial g(x_{k+1}) + \mathcal{N}_{\mathcal{B}_k}(x_{k+1}),
\]
where $\mathcal{N}_\mathcal{C}(z)$ denotes the normal cone of a set
$\mathcal{C}$ at a point $z$. The normal cone of a ball is explicit:
$\mathcal{N}_{\mathcal{B}_k}(z) = \{0\}$ if $\|z - x_k\| < t_k$, and
$\mathcal{N}_{\mathcal{B}_k}(z) = \{ \lambda (z - x_k) : \lambda \geq 0 \}$ if
$\|z - x_k\| = t_k$. Hence $x_{k+1}$ is characterized by the KKT conditions:
there exists $\lambda_k \geq 0$ such that
\begin{equation}\label{eq:kkt}
	\nabla f(x_k) + \partial g(x_{k+1}) + \lambda_k(x_{k+1} - x_k) = 0,
	\quad
	\|x_{k+1} - x_k\| \leq t_k,
	\quad
	\lambda_k \left( \|x_{k+1} - x_k\| - t_k \right) = 0,
\end{equation}
where $\partial g(x_{k+1})$ denotes a subgradient of $g$ at $x_{k+1}$, and the
last relation is the complementary slackness condition.

\textbf{Step 2: Deductions from the KKT conditions.} Taking the inner product of the stationarity condition in (\ref{eq:kkt}) with \( x_{k+1} - y \) gives the identity

\begin{equation}\label{eq:dekkt}
	\langle \nabla f(x_k), x_{k+1} - y \rangle + \langle \partial g(x_{k+1}), x_{k+1} - y \rangle + \lambda_k \langle x_{k+1} - x_k, x_{k+1} - y \rangle = 0. 
\end{equation}

Since $m_k(\cdot)$ is convex, we have 
$$
m_k(y) \geq m_k(x_{k+1}) + \langle \nabla f(x_k) + \partial g(x_{k+1}), y - x_{k+1} \rangle. 
$$
Then from the assumption that $m_k(x_{k+1}) > m_k(y)$, we get 
$$
\langle \nabla f(x_k) + \partial g(x_{k+1}), x_{k+1} -y  \rangle \geq m_k(x_{k+1}) - m_k(y) >0. 
$$ 
Combining the above inequality and (\ref{eq:dekkt}), we arrive at $\lambda_k \langle x_k - x_{k+1}, x_{k+1} - y \rangle > 0$. Then we have 

\begin{equation}\label{eq:dis-inner}
	\lambda_k > 0, \quad 	\langle x_k - x_{k+1}, x_{k+1} - y \rangle > 0. 
\end{equation}

\textbf{Step 3: Distance inequality.} Since $\lambda_k >0$, the complementary slackness condition in \eqref{eq:kkt} shows that $x_{k+1}$ lies on the boundary of the ball $\mathcal{B}_k$, i.e., 
\begin{equation}\label{eq:boundarytk}
	\|x_{k+1} - x_k\| = t_k.
\end{equation}

So, identity {\em (ii)} holds. Moreover, from (\ref{eq:dis-inner}), we can obtain 

\begin{eqnarray*}
	\|x_{k+1} - y\|^2 
	&=& \|x_k - y - (x_k - x_{k+1})\|^2 \\
	&= &\|x_k - y\|^2 - 2\langle x_k - x_{k+1}, x_{k+1} - y \rangle - \|x_k - x_{k+1}\|^2 \\
	&\overset{(\ref{eq:dis-inner}) + (\ref{eq:boundarytk})}{<}& \|x_k - y\|^2 - t_k^2.
\end{eqnarray*}

So, inequality {\em (iii)} holds. Next we prove claim {\em (iv)} by contradiction. If $x_{k+1}$ is not the unique solution of $\min_{z \in \mathcal{B}(x_k, t_k)} m_k(z)$, then there exists ${\tilde x}_{k+1} \neq x_{k+1}$ such that ${\tilde x}_{k+1} \in \mathcal{B}(x_k, t_k)$ and $m_k({\tilde x}_{k+1}) = m_k(x_{k+1})$. Define $x(\lambda) \eqdef \lambda x_{k+1} + (1-\lambda) {\tilde x}_{k+1}$ for any $0\leq \lambda \leq 1$. Then from the convexity of $m_k(\cdot)$, $x(\lambda)$ is also a solution of $\min_{z \in \mathcal{B}(x_k, t_k)} m_k(z)$. Specifically, $m_k(x(\lambda)) = m_k(x_{k+1}) > m_k(y)$. Then from {\em (ii)}, we have $\|x(\lambda) - x_k\| = t_k$. However, by {\em (ii)}, also $\|x_{k+1} - x_k\| = \|{\tilde x}_{k+1} - x_k\| = t_k$, and from the strong convexity of $\|\cdot\|^2$, for $0<\lambda <1$, 
\begin{eqnarray*}
	\|x(\lambda) - x_k\|^2 &=& \lambda \|x_{k+1} - x_k\|^2 + (1-\lambda) \|{\tilde x}_{k+1} - x_k\|^2 - \lambda(1-\lambda) \|x_{k+1} - {\tilde x}_{k+1}\|^2 \\
	&=& t_k^2 - \lambda(1-\lambda) \|x_{k+1} - {\tilde x}_{k+1}\|^2 < t_k^2, 
\end{eqnarray*}
which is a contradiction. Hence, the hypothesis is not true, and {\em (iv)} holds. 

It remains to prove claim {\em (v)}. From (\ref{eq:kkt}) and (\ref{eq:dis-inner}), we know $x_{k+1}$ satisfies the first-order optimality condition of the following strongly-convex optimization problem 
$$
\min_{z \in \R^d}  g(z) +  \langle \nabla f(x_k), z \rangle + \frac{\lambda_k}{2}\|z-x_k\|^2, 
$$
which is exactly the definition of $ \operatorname{Prox}_{\gamma_k g} \left(  x_k - \gamma_k \nabla f(x_k)  \right)$ with $\gamma_k = \frac{1}{\lambda_k}>0$. Furthermore, from (\ref{eq:kkt}) and (\ref{eq:dis-inner}), we can obtain $\|\nabla f(x_k) + \partial g(x_{k+1})\|>0$ and $\lambda_k = \frac{\|\nabla f(x_k) + \partial g(x_{k+1})\|}{t_k}$.

\subsection{Proof of Theorem \ref{th:Tk}}

Let $\mathcal{C} \eqdef \{ z \,|\, m_k(z) \leq m_k(y) \}$. Since $m_k$ is proper, closed and convex and $y \in \dom g$, the set $\mathcal{C} \subseteq \dom g$ is nonempty, closed and convex. Since $m_k$ attains its minimum over $\mathcal{B}(x_k, t)$ at $x_{k+1}(t)$, we have $t \in \cT_k(y)$, i.e., $m_k(x_{k+1}(t)) > m_k(y)$, if and only if $m_k(z) > m_k(y)$ for all $z \in \mathcal{B}(x_k, t)$, i.e., if and only if $\mathcal{B}(x_k, t) \cap \mathcal{C} = \emptyset$.

{\em (i)} If $m_k(x_k) \leq m_k(y)$, then $x_k \in \mathcal{C}$, so every ball around $x_k$ meets $\mathcal{C}$ and $\cT_k(y) = \emptyset$. If $m_k(x_k) > m_k(y)$, then $x_k \notin \mathcal{C}$, and since $\mathcal{C}$ is closed, $b \eqdef \operatorname{dist}(x_k, \mathcal{C}) > 0$. A ball $\mathcal{B}(x_k, t)$ misses $\mathcal{C}$ if and only if $t < b$, hence $\cT_k(y) = (0, b) \neq \emptyset$. Both claims of (i) follow.

{\em (ii)} By the proof of (i), $\cT_k(y) = (0,b)$, and $b \leq \|x_k - y\|$ since $y \in \mathcal{C}$. Let $p \eqdef \Pi_{\mathcal{C}}(x_k)$, so that $\|p - x_k\| = b$ and $m_k(p) \leq m_k(y)$. In fact $m_k(p) = m_k(y)$: otherwise, for $z_\lambda \eqdef p + \lambda(x_k - p)$ with $\lambda \in (0,1)$, convexity gives $m_k(z_\lambda) \leq (1-\lambda)m_k(p) + \lambda m_k(x_k) \leq m_k(y)$ for $\lambda$ small enough, so that $z_\lambda \in \mathcal{C}$ with $\|z_\lambda - x_k\| = (1-\lambda) b < b$, a contradiction. Consequently, every minimizer $z$ of $m_k$ over $\mathcal{B}(x_k, b)$ satisfies $m_k(z) \leq m_k(p) = m_k(y)$, i.e., $z \in \mathcal{C}$ with $\|z - x_k\| \leq b$, so $z$ is a point of $\mathcal{C}$ nearest to $x_k$; by the uniqueness of the projection, $z = p$. Hence $x_{k+1}(b) = p$ is the unique minimizer, $\|x_{k+1}(b) - x_k\| = b$, and $m_k(x_{k+1}(b)) = m_k(y)$. Since $y \in \mathcal{C}$, the projection inequality $\langle x_k - p, y - p \rangle \leq 0$ gives
$$
\|x_k - y\|^2 = \|x_k - p\|^2 + \|p - y\|^2 + 2\langle x_k - p, p - y\rangle \geq b^2 + \|x_{k+1}(b) - y\|^2.
$$
Since $\cT_k(y)$ depends on $y$ only through $m_k(y)$ and $m_k(x_{k+1}(b)) = m_k(y)$, we have $\cT_k(x_{k+1}(b)) = \cT_k(y) = (0,b)$. Finally, for $0 < t < b$, Theorem~\ref{th:another-ad}~{\em (iii)} with reference point $x_{k+1}(b) \in \dom g$ gives $\|x_{k+1}(t) - x_{k+1}(b)\|^2 < \|x_k - x_{k+1}(b)\|^2 - t^2 = b^2 - t^2$.

{\em (iii)} Since \( x_\star \) minimizes \(F\), the first-order optimality condition for (\ref{p:p}) gives $\nabla f(x_\star) + \partial g(x_\star) = 0$, where $\partial g(x_\star)$ denotes a subgradient of $g$ at $x_\star$. Then, for $t \in \cT_k^\star$, we have 
$$
g(x_{k+1}(t)) \geq g(x_\star) + \langle \partial g(x_\star), x_{k+1}(t) - x_\star \rangle = g(x_\star) - \langle \nabla f(x_\star), x_{k+1}(t) - x_\star \rangle, 
$$
which indicates that
\begin{eqnarray*}
	m_k(x_{k+1}(t)) - m_k(x_\star) &=& g(x_{k+1}(t)) - g(x_\star) + \langle \nabla f(x_k), x_{k+1}(t) - x_\star \rangle \\ 
	&\geq&  \langle \nabla f(x_k) - \nabla f(x_\star), x_{k+1}(t) - x_\star \rangle \\ 
	&>&0, 
\end{eqnarray*}
where the last inequality comes from the definition of $\cT_k^\star$. Therefore, $t \in \cT_k$ and $\cT_k^\star \subseteq \cT_k$.

\section{Proofs for Section \ref{sec:BGrecover}}

\begin{lemma}\label{lm:prox}
	If $y = \operatorname{Prox}_{\gamma g} \left(  x_k - \gamma \nabla f(x_k)  \right)$ for some $\gamma>0$, then for any $z \in \R^d$, we have 
	$$
	m_k(z) - m_k(y) \geq \frac{1}{\gamma} \langle y-x_k, y-z \rangle. 
	$$
\end{lemma}

\begin{proof}

	First, since $g$ is proper, closed and convex, the operator
	$\operatorname{Prox}_{\gamma g}(\cdot)$ is well defined, and by its definition
	$y = \operatorname{Prox}_{\gamma g}\bigl( x_k - \gamma \nabla f(x_k) \bigr)$
	is the unique solution of the strongly convex optimization problem
	$$
	\min_{z \in \R^d} \left\{ g(z) + \langle \nabla f(x_k), z - x_k \rangle
	+ \frac{1}{2\gamma} \|z-x_k\|^2 \right\}. 
	$$
	By Fermat's rule, $y$ is
	characterized by $0 \in \partial \left( \cdot \right)(y)$, where the
	subdifferential is that of the objective above; since its last two terms are
	finite and differentiable on all of $\R^d$, the subdifferential splits without
	any qualification condition, and we obtain
	$$
	\partial g(y) + \nabla f(x_k) + \frac{1}{\gamma} (y-x_k) = 0,
	$$
	where $\partial g(y)$ is a subgradient of $g$ at $y$. Then, from the convexity
	of $g$, for any $z \in \R^d$,
	\begin{eqnarray*}
		g(z) &\geq& g(y) + \langle \partial g(y), z-y \rangle \\
		&=& g(y) - \langle \nabla f(x_k), z-y \rangle
		- \frac{1}{\gamma} \langle y-x_k, z-y \rangle,
	\end{eqnarray*}
	which implies that
	$$
	m_k(z) - m_k(y) \geq \frac{1}{\gamma} \langle y-x_k, y-z \rangle .
	$$
	
\end{proof}

We write $y(\gamma) \eqdef \operatorname{Prox}_{\gamma g} \left( x_k - \gamma \nabla f(x_k) \right)$ for $\gamma > 0$, so that $y_k = y(\gamma)$ and $y_k^{(\alpha)} = y(\alpha\gamma)$.

\begin{lemma}\label{lm:prox-mono}
	For $0 < \gamma' \leq \gamma''$,
	$$
	\frac{\gamma'}{\gamma''}\, \|y(\gamma'') - x_k\| \leq \|y(\gamma') - x_k\| \leq \|y(\gamma'') - x_k\|,
	$$
	and the first inequality is an equality if and only if $y(\gamma') = x_k + \frac{\gamma'}{\gamma''}\, (y(\gamma'') - x_k)$.
\end{lemma}

\begin{proof}
	Let $u \eqdef y(\gamma'') - x_k$, $v \eqdef y(\gamma') - x_k$ and $\alpha \eqdef \nicefrac{\gamma'}{\gamma''} \in (0, 1]$. Lemma \ref{lm:prox} with $(y, z) = (y(\gamma''), y(\gamma'))$ and with $(y, z) = (y(\gamma'), y(\gamma''))$ gives
	\begin{align*}
	\gamma'' \bigl( m_k(y(\gamma')) - m_k(y(\gamma'')) \bigr) + \langle u, v - u \rangle &\geq 0, \\
	\gamma' \bigl( m_k(y(\gamma'')) - m_k(y(\gamma')) \bigr) + \langle v, u - v \rangle &\geq 0.
	\end{align*}
	Multiplying the first inequality by $\alpha$ and adding the second, the model values cancel and we obtain $(1+\alpha) \langle u, v \rangle \geq \alpha\|u\|^2 + \|v\|^2$. By the Cauchy--Schwarz inequality, $\|v\|^2 - (1+\alpha)\|u\|\,\|v\| + \alpha\|u\|^2 \leq 0$, i.e., $(\|v\| - \alpha\|u\|)(\|v\| - \|u\|) \leq 0$, which gives $\alpha\|u\| \leq \|v\| \leq \|u\|$. If $\|v\| = \alpha\|u\|$, the Cauchy--Schwarz inequality is an equality with $\langle u, v\rangle \geq 0$, so $v$ is a nonnegative multiple of $u$, hence $v = \alpha u$; the converse is clear.
\end{proof}

\subsection{Proof of Theorem \ref{th:recoverPGD}}

{\em (i)} Since $y_k \neq x_k$ is the unique minimizer of $m_k(z) + \frac{1}{2\gamma}\|z - x_k\|^2$, we have $m_k(y_k) < m_k(x_k)$. Let $\mathcal{C} \eqdef \{ z \,|\, m_k(z) \leq m_k(y_k) \}$. For $z \in \mathcal{C}$, Lemma \ref{lm:prox} gives $\langle x_k - y_k, z - y_k \rangle \leq \gamma \bigl( m_k(z) - m_k(y_k) \bigr) \leq 0$, hence
$$
\|z - x_k\|^2 = \|z - y_k\|^2 + \|y_k - x_k\|^2 + 2\langle z - y_k, y_k - x_k \rangle \geq \|y_k - x_k\|^2 = b_k^2.
$$
Thus $\mathcal{B}(x_k, t) \cap \mathcal{C} = \emptyset$ for $t < b_k$, i.e., $t \in \cT_k(y_k)$ (see the proof of Theorem \ref{th:Tk}), whereas $b_k \notin \cT_k(y_k)$ since $y_k \in \mathcal{C} \cap \mathcal{B}(x_k, b_k)$; hence $\cT_k(y_k) = (0, b_k)$. By Theorem \ref{th:Tk} {\em (ii)}, $x_{k+1}(b_k)$ is unique and $\|x_{k+1}(b_k) - y_k\|^2 \leq \|x_k - y_k\|^2 - b_k^2 = 0$, so $x_{k+1}(b_k) = y_k$.

{\em (ii)} Lemma \ref{lm:prox-mono} with $\gamma' = \alpha\gamma$ and $\gamma'' = \gamma$ gives $t^{(\alpha)} = \|y_k^{(\alpha)} - x_k\| \geq \alpha\|y_k - x_k\| > 0$, with equality if and only if $y_k^{(\alpha)} = (1-\alpha)x_k + \alpha y_k$. Since $y_k^{(\alpha)} \neq x_k$, part {\em (i)} applied to $y_k^{(\alpha)} = y(\alpha\gamma)$ gives $x_{k+1}(t^{(\alpha)}) = y_k^{(\alpha)}$.

{\em (iii)} By {\em (i)}, $t \in \cT_k(y_k)$, so Theorem \ref{th:another-ad} {\em (v)} gives $x_{k+1}(t) = y(\gamma_t)$ for some $\gamma_t > 0$, with $\|y(\gamma_t) - x_k\| = t = \alpha\|y_k - x_k\|$. If $\gamma_t \geq \gamma$, the second inequality of Lemma \ref{lm:prox-mono} gives $t \geq \|y_k - x_k\| > t$, a contradiction; hence $\gamma_t < \gamma$, and the first inequality of Lemma \ref{lm:prox-mono} with $\gamma' = \gamma_t$ and $\gamma'' = \gamma$ gives $t \geq \frac{\gamma_t}{\gamma}\|y_k - x_k\|$, i.e., $\gamma_t \leq \alpha\gamma$.

\section{Proofs for the convex case}

The following elementary lemma captures the counting argument shared by Theorems \ref{th:convex-f-prac}, \ref{th:convex-G-prac} and \ref{th:non-convex-withoutPL}.

\begin{lemma}\label{lm:dichotomy}
	Let $\{e_k\}_{k\geq0} \subseteq [0, +\infty)$, let $\{\delta_k\}_{k\geq0}$ be a positive sequence, and let $\cK \subseteq \mathbb{N}$.
	{\em (i)} If $\sum_{k=0}^{+\infty} \delta_k = +\infty$ and $e_{k+1} < e_k - \delta_k$ for every $k \notin \cK$, then $\cK$ is infinite.
	{\em (ii)} If $e_{k+1} \leq e_k - \delta_k$ for every $k \geq 0$, then $\frac{1}{K}\sum_{k=0}^{K-1} \delta_k \leq \frac{e_0}{K}$ for every $K \geq 1$.
\end{lemma}

\begin{proof}
	{\em (i)} If $\cK$ were finite, then for $k_1 > \max\cK$ (or $k_1 = 0$ if $\cK = \emptyset$), summing the inequality over $k = k_1, \dots, K$ would give $\sum_{k=k_1}^K \delta_k < e_{k_1} - e_{K+1} \leq e_{k_1}$ for every $K \geq k_1$, contradicting $\sum_k \delta_k = +\infty$. {\em (ii)} Summing over $k = 0, \dots, K-1$ gives $\sum_{k=0}^{K-1}\delta_k \leq e_0 - e_K \leq e_0$.
\end{proof}

\subsection{Proof of Theorem \ref{th:convex-f}}

Assume $F(x_k) > F(x_\star)$. From the convexity of $f$, 
\begin{equation}\label{eq:fconvex-inq}
	f(x_\star) \geq f(x_k) + \langle \nabla f(x_k), x_\star - x_k \rangle. 
\end{equation}
Then we have 
$$
m_k(x_k) - m_k(x_\star) = g(x_k) - g(x_\star) + \langle \nabla f(x_k), x_k - x_\star \rangle \geq F(x_k) - F(x_\star) >0. 
$$ 
From Theorem \ref{th:Tk}, we know $\cT_k = (0, b_k)$ with $0<b_k \leq \|x_k - x_\star\|$, and 
\begin{equation}\label{eq:xk+1-bk2}
	\|x_{k+1} - x_\star\|^2 \leq \|x_k - x_\star\|^2 - b_k^2. 
\end{equation}
Since $f$ is $L$-smooth, we have 
\begin{eqnarray*}
	f(x_{k+1}) + g(x_{k+1}) &\leq& g(x_{k+1}) +  f(x_k) + \langle \nabla f(x_k), x_{k+1} - x_k \rangle + \frac{L}{2} \|x_{k+1} - x_k\|^2 \\ 
	&=& g(x_\star) +  f(x_k) + \langle \nabla f(x_k), x_\star - x_k \rangle + \frac{Lb_k^2}{2} \\ 
	&\overset{(\ref{eq:fconvex-inq})}{\leq}&  f(x_\star) + g(x_\star) +  \frac{Lb_k^2}{2}, 
\end{eqnarray*}
where we use $m_k(x_{k+1}) = m_k(x_\star)$ and $\|x_{k+1} - x_k\|=b_k$ by Theorem \ref{th:Tk} {\em (ii)} in the first equality. Then we get $b_k \geq \frac{\sqrt{2(F(x_{k+1}) - F(x_\star))}}{\sqrt{L}}$. 
Combining the above inequality and (\ref{eq:xk+1-bk2}) yields that 
\begin{equation}\label{eq:xk+1-deltafk+1}
	\|x_{k+1} - x_\star\|^2 \leq \|x_k - x_\star\|^2 - \frac{2(F(x_{k+1}) - F(x_\star))}{L}. 
\end{equation}

Now assume $F(x_k) = F(x_\star)$. Then $t_k=0$, and $x_{k+1} = x_k$ and $F(x_{k+1}) = F(x_\star)$. Thus inequality (\ref{eq:xk+1-deltafk+1}) also holds. 

Summing the inequality (\ref{eq:xk+1-deltafk+1}) for $k= 0, \dots, K-1$, we get 
$$
\|x_K - x_\star\|^2 \leq \|x_0 - x_\star\|^2 - \frac{2}{L} \sum_{k=1}^K (F(x_k) - F(x_\star)). 
$$

Since the left-hand side is nonnegative, this implies that 
$$
\frac{1}{K} \sum_{k=1}^K (F(x_k) - F(x_\star)) \leq \frac{L\|x_0 - x_\star\|^2}{2K}. 
$$
By Jensen’s inequality, we have 
$$
F({\hat x}_K) - F(x_\star) \leq \frac{1}{K} \sum_{k=1}^K (F(x_k) - F(x_\star)) \leq \frac{L\|x_0 - x_\star\|^2}{2K}. 
$$

\subsection{Proof of Theorem \ref{th:convex-f-prac}}

{\em (i)} For any $k\geq 1$ and $k \notin \cK$, we prove that $t_{k-1} \in \cT_{k-1}$ by contradiction. If $t_{k-1} \notin \cT_{k-1}$, then $m_{k-1}(x_k) \leq m_{k-1}(x_\star)$, and by the $L$-smoothness of $f$, we can obtain 
\begin{eqnarray*}
	f(x_k) + g(x_k)&\leq& g(x_k) + f(x_{k-1}) + \langle \nabla f(x_{k-1}), x_k - x_{k-1} \rangle + \frac{L}{2} \|x_k - x_{k-1}\|^2 \\ 
	&\leq& g(x_\star) +  f(x_{k-1}) + \langle \nabla f(x_{k-1}), x_\star - x_{k-1} \rangle + \frac{Lt_{k-1}^2}{2} \\ 
	&\leq&  f(x_\star) + g(x_\star) + a_k^2, 
\end{eqnarray*}
where the last inequality comes from the convexity of $f$. However, $k \notin \cK$ implies that $F(x_k) - F(x_\star) > a_k^2$, which is a contradiction. Hence the hypothesis is not true, and $t_{k-1} \in \cT_{k-1}$ for all $k\geq 1$ and $k \notin \cK$. Then by Theorem \ref{th:another-ad} {\em (iii)}, we have 
$$
\|x_k - x_\star\|^2 < \|x_{k-1} - x_\star\|^2 - t_{k-1}^2 = \|x_{k-1} - x_\star\|^2 -  \frac{2}{L} a_k^2, 
$$
for all $k\geq 1$ and $k \notin \cK$. 

{\em (ii)} By {\em (i)}, Lemma \ref{lm:dichotomy}(i) applies with $e_k = \|x_k - x_\star\|^2$, $\delta_k = \frac{2}{L}a_{k+1}^2$ and the set $\{ k \,|\, k+1 \in \cK \}$, which is therefore infinite, and so is $\cK$.

\subsection{Proof of Corollary \ref{co:convex-f-prac}}

If $\min_{0 \leq k \leq K} F(x_k) - F(x_\star) > a^2$ for some $K \geq \frac{L\|x_0-x_\star\|^2}{2a^2}$, then $k \notin \cK$ for all $1 \leq k \leq K$, and summing the inequality of Theorem \ref{th:convex-f-prac} {\em (i)} over these $k$ gives $0\leq \|x_K - x_\star\|^2 < \|x_0 - x_\star\|^2 - \frac{2Ka^2}{L}$, i.e., $K < \frac{L\|x_0-x_\star\|^2}{2a^2}$, a contradiction. The second statement is the case $a = \sqrt{\epsilon}$, for which $t_k = \sqrt{\nicefrac{2}{L}}\, a_{k+1} = \sqrt{\nicefrac{2\epsilon}{L}}$.

\subsection{Proof of Theorem \ref{th:convex-halving}}\label{sec:proof-halving}

The analysis of the halving rule \eqref{eq:halving} rests on the following two lemmas, which are also used in the strongly convex and nonconvex cases. For $L' \geq L$, we call a step from $x_k$ with radius $t$ \emph{short} (with respect to $L'$) if $t \leq \|y'_k - x_k\|$, where $y'_k \eqdef \operatorname{Prox}_{g/L'}(x_k - \nabla f(x_k)/L')$, and \emph{long} otherwise; note that $\|y'_k - x_k\| = \|G_{\nicefrac{1}{L'}}(x_k)\|/L'$.

\begin{lemma}\label{lm:short-long}
	Let Assumption \ref{as:basic} hold, and assume that $f$ is $L$-smooth on an open convex set containing $\dom g$. Let $L' \geq L$, $x_k \in \dom g$, $t > 0$, and let $x_{k+1}$ be obtained from $x_k$ by one step of Broximal GD with radius $t$.
	\begin{itemize}
		\item[(i)] If the step is short, then $F(x_{k+1}) \leq F(x_k) - \frac{L'}{2}t^2$.
		\item[(ii)] If the step is long, then $F(x_{k+1}) \leq F(x_k) + \frac{L'}{2}t^2$; if moreover $f$ is convex on the above set, then for every $z \in \dom g$,
		$$
		F(x_{k+1}) \leq F(z) + L' t\, \|z - x_k\| + \frac{L'}{2}t^2,
		$$
		and if $f$ is $\mu$-strongly convex on the above set, with $\mu > 0$, then $F(x_{k+1}) - F_\star < \frac{L'}{2}\left( 1 + \frac{L'}{\mu} \right) t^2$.
	\end{itemize}
\end{lemma}

\begin{proof}
	Since $f$ is also $L'$-smooth, claim (i) and the first part of (ii) are Lemma \ref{lm:non-convex-step} {\em (ii)}--{\em (iv)} with $\gamma = \nicefrac{1}{L'}$, for which $r_k = \|y'_k - x_k\|$. For the second part of (ii), since $\|y'_k - x_k\| < t$, the point $y'_k$ belongs to $\mathcal{B}(x_k, t)$, so that $m_k(x_{k+1}) \leq m_k(y'_k)$. Together with the smoothness of $f$ and $\|x_{k+1} - x_k\| \leq t$, this gives
	$$
	F(x_{k+1}) \leq F(x_k) + m_k(x_{k+1}) - m_k(x_k) + \frac{L'}{2}t^2 \leq F(x_k) + m_k(y'_k) - m_k(x_k) + \frac{L'}{2}t^2.
	$$
	Since $y'_k$ minimizes the $L'$-strongly convex function $z \mapsto m_k(z) + \frac{L'}{2}\|z - x_k\|^2$, for every $z \in \dom g$,
	\begin{align*}
	m_k(y'_k) &\leq m_k(z) + \frac{L'}{2}\left( \|z - x_k\|^2 - \|z - y'_k\|^2 - \|y'_k - x_k\|^2 \right) \\
	&\leq m_k(z) + L' \|z - x_k\|\,\|y'_k - x_k\| \leq m_k(z) + L' t\,\|z - x_k\|,
	\end{align*}
	where the second inequality uses $\|z - x_k\|^2 - \|z - y'_k\|^2 = 2\langle z - x_k, y'_k - x_k\rangle - \|y'_k - x_k\|^2$. Finally, by the convexity of $f$, $F(x_k) + m_k(z) - m_k(x_k) = f(x_k) + \langle \nabla f(x_k), z - x_k\rangle + g(z) \leq F(z)$. If $f$ is $\mu$-strongly convex, the last inequality improves to $F(x_k) + m_k(z) - m_k(x_k) \leq F(z) - \frac{\mu}{2}\|z - x_k\|^2$; taking $z = x_\star$ and writing $\tilde d \eqdef \|y'_k - x_k\| < t$ and $R \eqdef \|x_\star - x_k\|$, we obtain from the second inequality above
	$$
	F(x_{k+1}) - F_\star \leq L' \tilde d R - \frac{\mu}{2}R^2 + \frac{L'}{2}t^2 \leq \frac{L'^2 \tilde d^2}{2\mu} + \frac{L'}{2}t^2 < \frac{L'}{2}\left( 1 + \frac{L'}{\mu} \right) t^2,
	$$
	where the second inequality maximizes over $R \geq 0$.
\end{proof}

\begin{lemma}\label{lm:window}
	Let Assumption \ref{as:basic} hold, and assume that $f$ is $L$-smooth on an open convex set containing $\dom g$. Let $c > 0$, $L_c \eqdef \max\{L, 2c\}$, and let $x_{a+1}, \dots, x_{a+T}$ be obtained from $x_a \in \dom g$ by $T$ steps of Broximal GD with the same radius $\tau > 0$.
	\begin{itemize}
		\item[(i)] For every $a \leq l \leq a+T$, $F(x_l) \leq F(x_a) + \frac{l-a}{2}L_c\tau^2$.
		\item[(ii)] If $F(x_a) - F(x_{a+T}) < cT\tau^2$, then there is $j \in \{a, \dots, a+T-1\}$ such that the step from $x_j$ is long with respect to $L_c$, i.e., $\|G_{\nicefrac{1}{L_c}}(x_j)\| < L_c\tau$. Let $j$ be the largest such index. If moreover $f$ is convex on the above set, then for every $x_\star \in \cX_\star$,
		$$
		F(x_{a+T}) - F_\star \leq L_c \tau\, \|x_j - x_\star\| + \frac{L_c}{2}\tau^2,
		$$
		and if $f$ is $\mu$-strongly convex on the above set, with $\mu > 0$, then $F(x_{a+T}) - F_\star < \frac{L_c}{2}\left( 1 + \frac{L_c}{\mu} \right)\tau^2$.
	\end{itemize}
\end{lemma}

\begin{proof}
	{\em (i)} By Lemma \ref{lm:short-long} with $L' = L_c$, each step increases $F$ by at most $\frac{L_c}{2}\tau^2$. {\em (ii)} If all the $T$ steps were short, Lemma \ref{lm:short-long} {\em (i)} would give $F(x_a) - F(x_{a+T}) \geq T\frac{L_c}{2}\tau^2 \geq cT\tau^2$; hence some step is long. With $j$ the last such index, the steps from $x_{j+1}, \dots, x_{a+T-1}$ are short and decrease $F$ by Lemma \ref{lm:short-long} {\em (i)}, so that $F(x_{a+T}) \leq F(x_{j+1})$, and Lemma \ref{lm:short-long} {\em (ii)} with $L' = L_c$ (and $z = x_\star$ in the convex case) bounds $F(x_{j+1}) - F_\star$.
\end{proof}

We now prove Theorem \ref{th:convex-halving}. We call the indices $iT$, $i \geq 0$, at which the test in \eqref{eq:halving} compares objective values, \emph{checkpoints}. The $i$-th window consists of the steps of the iterations $iT, \dots, (i+1)T-1$, which lead from $x_{iT}$ to $x_{(i+1)T}$; we say that it \emph{starts} at the checkpoint $iT$ and \emph{ends} at the checkpoint $(i+1)T$, and we call it \emph{halving} if $\tau_{i+1} = \nicefrac{\tau_i}{2}$ and \emph{surviving} otherwise, and we write $\Phi_k \eqdef F(x_k) - F_\star$ and $\rho_j \eqdef t_0 2^{-j}$, $j \geq 0$, for the successive values of the radius (to be distinguished from $t_k$, the radius used at iteration $k$).\footnote{The symbol $\rho_j$ should not be confused with the parameter $\rho$ of the test problems in Section~\ref{sec:exp}.} A surviving window satisfies $\Phi_{(i+1)T} \leq \Phi_{iT} - cT\tau_i^2$. Since $\Phi_k \geq 0$, each radius $\rho_j$ is used in finitely many windows, the last of which is the only halving one; let $k_j$ be the checkpoint at which this window ends, and $k_{-1} \eqdef 0$.

\emph{Distance to a solution.} Let $x_\star$ be the projection of $x_0$ onto the nonempty closed convex set $\cX_\star$, and $R_k \eqdef \|x_k - x_\star\|$. If $t_k \in \cT_k(x_\star)$, then $R_{k+1} < R_k$ by Theorem \ref{th:another-ad} {\em (iii)}. If $t_k \notin \cT_k(x_\star)$, then $R_{k+1} \leq R_k + t_k$, and $m_k(x_{k+1}) \leq m_k(x_\star)$, so that the $L$-smoothness and the convexity of $f$ give, as in the proof of Theorem \ref{th:convex-f-prac},
$$
F(x_{k+1}) \leq F(x_k) + m_k(x_{k+1}) - m_k(x_k) + \frac{L}{2}t_k^2 \leq F(x_k) + m_k(x_\star) - m_k(x_k) + \frac{L}{2}t_k^2 \leq F_\star + \frac{L}{2}t_k^2.
$$
Suppose that $t_k \notin \cT_k(x_\star)$ for some iteration $k$ with radius $t_k = \rho_j$, and let $k'$ be the first such iteration. Let $p$ be the checkpoint at which the window containing the step of iteration $k'$ ends, i.e., $p = (i+1)T$ if $iT \leq k' < (i+1)T$. Since each of the at most $T-1$ steps between $k'+1$ and $p$ increases $F$ by at most $\frac{L_c}{2}\rho_j^2$ (Lemma \ref{lm:short-long} with $L' = L_c$), we have $\Phi_p \leq \frac{L}{2}\rho_j^2 + \frac{T-1}{2}L_c\rho_j^2 \leq \frac{T}{2}L_c\rho_j^2$. Since each surviving window decreases the gap by at least $cT\rho_j^2$, at most $\nicefrac{L_c}{2c}$ surviving windows with radius $\rho_j$ follow $p$, and then the halving window. Hence at most
$$
M \eqdef T + \frac{L_c}{2c}T + T = T\left( 2 + \frac{L_c}{2c} \right)
$$
iterations with radius $\rho_j$ take place from $k'$ on, so that $R_k$ increases by at most $M\rho_j$ while the radius is $\rho_j$. Summing over $j$, we obtain
\begin{equation}\label{eq:halving-Rbar}
	\|x_k - x_\star\| \leq \|x_0 - x_\star\| + M\sum_{j \geq 0}\rho_j = \operatorname{dist}(x_0, \cX_\star) + 2Mt_0 = \bar R, \qquad k \geq 0.
\end{equation}

\emph{Gap after a stage.} Let
\begin{equation}\label{eq:halving-Ej}
	E_j \eqdef L_c\rho_j\bar R + \frac{L_c}{2}\rho_j^2, \qquad j \geq 0.
\end{equation}
By Lemma \ref{lm:window} {\em (ii)} applied to the halving window with radius $\rho_j$ and by \eqref{eq:halving-Rbar}, $\Phi_{k_j} \leq E_j$. Moreover, $\Phi_k \leq E_j$ for every $k \geq k_j$. Indeed, the steps from $x_k$ with $k \geq k_j$ use radii $\rho_i$ with $i > j$; a short step decreases $F$ by Lemma \ref{lm:short-long} {\em (i)}, and a long step from $x_k$ gives $\Phi_{k+1} \leq L_c\rho_i\|x_k - x_\star\| + \frac{L_c}{2}\rho_i^2 \leq E_i \leq E_j$ by Lemma \ref{lm:short-long} {\em (ii)} and \eqref{eq:halving-Rbar}.

\emph{Length of the stages.} The radius $t_0$ is used from the checkpoint $0$ on (included), and the radius $\rho_j$, $j \geq 1$, from the checkpoint $k_{j-1}$ on (included), at which the gap is at most $E_{j-1}$. Since each surviving window decreases the gap by at least $cT\rho_j^2$, their number is at most $\nicefrac{\Phi_0}{cTt_0^2}$ for $j = 0$ and $\nicefrac{E_{j-1}}{cT\rho_j^2}$ for $j \geq 1$. Adding the halving window and using $\rho_{j-1} = 2\rho_j$, we obtain $k_0 \leq N_0$ and
$$
k_j - k_{j-1} \leq T + \frac{E_{j-1}}{c\,\rho_j^2} = T + \frac{2L_c\bar R}{c\,\rho_j} + \frac{2L_c}{c} = N_j, \qquad j \geq 1,
$$
which is {\em (i)}.

\emph{Claim (ii).} If $K \geq \sum_{j=0}^{J}N_j$, then $K \geq k_J$, and therefore $\Phi_K \leq E_J = L_c\rho_J\left( \bar R + \frac{1}{2}\rho_J \right)$.

\emph{Claim (iii).} Let $J = J_\epsilon$, which is nonnegative since $\epsilon \leq L_c\bar Rt_0$. Then $\rho_J \leq \nicefrac{8\epsilon}{9L_c\bar R}$. Since $\bar R \geq T(4 + \nicefrac{L_c}{c})t_0 \geq 4t_0 \geq 4\rho_J$, we have $\bar R + \frac{1}{2}\rho_J \leq \frac{9}{8}\bar R$, and {\em (ii)} gives $F(x_K) - F_\star \leq \frac{9}{8}L_c\rho_J\bar R \leq \epsilon$ for $K \geq \sum_{j=0}^{J}N_j$. Moreover, $2^J < \nicefrac{9L_c\bar Rt_0}{4\epsilon}$, so that
$$
\sum_{j=1}^{J}N_j = J\left(T + \frac{2L_c}{c}\right) + \frac{L_c\bar R}{c\,t_0}\left( 2^{J+2} - 4 \right) < J\left(T + \frac{2L_c}{c}\right) + \frac{9L_c^2\bar R^2}{c\,\epsilon}.
$$

\subsection{Proof of Theorem \ref{th:convex-G-prac}}

Write $v_k \eqdef \nabla f(x_k) + \partial g(x_k)$ and $\Delta_k \eqdef F(x_k) - F(x_\star) \geq 0$. We first show that if $\Delta_k > 0$, then
\begin{equation}\label{eq:G-admissible}
	(0, \nicefrac{\Delta_k}{\|v_k\|}) \subseteq \cT_k .
\end{equation}
Indeed, $\Delta_k > 0$ implies $x_k \notin \cX_\star$, hence $\|v_k\| > 0$. By the convexity of $f$, $\langle \nabla f(x_k), x_k - x_\star \rangle \geq f(x_k) - f(x_\star)$, so that
\begin{equation}\label{eq:mkxk-mkxstar}
	m_k(x_k) - m_k(x_\star) = g(x_k) - g(x_\star) + \langle \nabla f(x_k), x_k - x_\star \rangle \geq \Delta_k ,
\end{equation}
and for $0 < t < \nicefrac{\Delta_k}{\|v_k\|}$, using the convexity of $g$, the Cauchy--Schwarz inequality and $\|x_{k+1}(t) - x_k\| \leq t$,
\begin{eqnarray*}
	m_k(x_{k+1}(t)) - m_k(x_\star) &=& g(x_{k+1}(t)) - g(x_k) + \langle \nabla f(x_k), x_{k+1}(t) - x_k\rangle \\
	&& +\, m_k(x_k) - m_k(x_\star) \\
	&\overset{(\ref{eq:mkxk-mkxstar})}{\geq}& \langle v_k, x_{k+1}(t) - x_k \rangle + \Delta_k \ \geq\ \Delta_k - t\|v_k\| \ >\ 0,
\end{eqnarray*}
i.e., $t \in \cT_k$.

{\em (i)} For $k \notin \cK$, we have $\Delta_k > a_k > 0$ and $t_k = \nicefrac{a_k}{\|v_k\|} < \nicefrac{\Delta_k}{\|v_k\|}$, so $t_k \in \cT_k$ by \eqref{eq:G-admissible}. Then Theorem \ref{th:another-ad} {\em (iii)} and $\|v_k\| \leq G$ give
$$
\|x_{k+1} - x_\star\|^2 < \|x_k - x_\star\|^2 - t_k^2 \leq \|x_k - x_\star\|^2 - \frac{a_k^2}{G^2}.
$$

{\em (ii)} By {\em (i)}, Lemma \ref{lm:dichotomy}(i) applies with $e_k = \|x_k - x_\star\|^2$ and $\delta_k = \nicefrac{a_k^2}{G^2}$, so $\cK$ is infinite.

{\em (iii)} If $\Delta_k = 0$ or $\|v_k\| = 0$, then $t_k = 0$, $x_{k+1} = x_k$, and $\Delta_k = 0$ in both cases (in the second by the convexity of $F$), so $\|x_{k+1} - x_\star\|^2 \leq \|x_k - x_\star\|^2 - \nicefrac{\Delta_k^2}{G^2}$ holds trivially. Otherwise $t_k = \nicefrac{\Delta_k}{\|v_k\|} > 0$, and by \eqref{eq:G-admissible} either $t_k \in \cT_k$ or $t_k = \sup \cT_k$, so Theorem \ref{th:another-ad} {\em (iii)} or Theorem \ref{th:Tk} {\em (ii)} gives
$$
\|x_{k+1} - x_\star\|^2 \leq \|x_k - x_\star\|^2 - t_k^2 \leq \|x_k - x_\star\|^2 - \frac{\Delta_k^2}{G^2}.
$$
Lemma \ref{lm:dichotomy}(ii) with $e_k = \|x_k - x_\star\|^2$ and $\delta_k = \nicefrac{\Delta_k^2}{G^2}$ yields the claim.

\subsection{Proof of Corollary \ref{co:convex-G-prac}}

If $\min_{0 \leq k \leq K-1} F(x_k) - F(x_\star) > a$ for some $K \geq \frac{G^2\|x_0-x_\star\|^2}{a^2}$, then $k \notin \cK$ for all $0 \leq k \leq K-1$, and summing the inequality of Theorem \ref{th:convex-G-prac} {\em (i)} over these $k$ gives $0\leq \|x_K - x_\star\|^2 < \|x_0 - x_\star\|^2 - \frac{Ka^2}{G^2}$, i.e., $K < \frac{G^2\|x_0-x_\star\|^2}{a^2}$, a contradiction. The second statement is the case $a = \epsilon$.

\section{Proofs for the strongly convex case}

\subsection{Proof of Theorem \ref{th:strongly-convex-theory}}

First, assume $\nabla f(x_k) \neq \nabla f(x_\star)$. Then by Lemma~C.2 of \citet{richtarik2026local}, strict convexity ensures that $x_k \neq x_\star$. Thus $t_k>0$. For any $0<t<t_k$, we have 
\begin{eqnarray*}
	&&\langle \nabla f(x_k) - \nabla f(x_\star), x_{k+1}(t) - x_\star \rangle \\
	&=& \langle \nabla f(x_k) - \nabla f(x_\star), x_k - x_\star \rangle  + \langle \nabla f(x_k) - \nabla f(x_\star), x_{k+1}(t) - x_k \rangle \\ 
	&\geq& \langle \nabla f(x_k) - \nabla f(x_\star), x_k - x_\star \rangle  - t \|\nabla f(x_k) - \nabla f(x_\star)\| \\ 
	&\geq& \theta \|\nabla f(x_k) - \nabla f(x_\star)\| \|x_k - x_\star\| -  t \|\nabla f(x_k) - \nabla f(x_\star)\| \\ 
	&>& \|\nabla f(x_k) - \nabla f(x_\star)\| (\theta \|x_k-x_\star\| - t_k)\\
	&=& 0, 
\end{eqnarray*}
where we use the Cauchy–Schwarz inequality and $\|x_{k+1}(t)-x_k\| \leq t$ in the first inequality, use Lemma C.3 of \citet{richtarik2026local} in the second inequality, and use $t<t_k$ in the last inequality. Therefore, $t \in \cT_k^\star$ for any $0<t<t_k$. From $\cT_k^\star \subseteq \cT_k$ in Theorem \ref{th:Tk} {\em (iii)}, we know either $t_k \in \cT_k$ or $t_k = \sup \cT_k$. Then, from Theorem \ref{th:another-ad} {\em (iii)} and Theorem \ref{th:Tk} {\em (ii)}, we have 
$$
\|x_{k+1} - x_\star\|^2 \leq \|x_k - x_\star\|^2 - t_k^2 = (1 - \theta^2)\|x_k - x_\star\|^2 = \frac{(L-\mu)^2}{(L+\mu)^2}\|x_k - x_\star\|^2
$$
establishing the claimed contraction. 

If $\nabla f(x_k) = \nabla f(x_\star)$, then by Lemma~C.2 of \citet{richtarik2026local}, strict convexity ensures that $x_k = x_\star$. Hence, $x_{k+1} = x_k=x_\star$, and the above inequality also holds. Unrolling the above recursion, we can get the bound for $\|x_K-x_\star\|^2$.

Theorems~\ref{th:strongly-convex-prac} and \ref{th:non-convex-PL} share the following scalar argument. 

\begin{lemma}\label{lm:envelope}
Let $\{e_k\}_{k\geq0}\subseteq[0,+\infty)$, ${\bar C}>0$, $r\in(0,1)$, $\kappa>0$ and $\omega\geq1$ satisfy $\frac{\kappa r}{1-r}>\omega$, and let $\cK\eqdef\{k \,|\, e_k\leq {\bar C}r^k\}$. Assume that {\em (a)} $e_{k+1}<e_k-\kappa{\bar C}r^k$ for every $k\notin\cK$, and {\em (b)} $e_k\leq\omega{\bar C}r^{k-1}$ for every $k\geq1$ with $k-1\in\cK$. Then
{\em (i)} if $\frac{\kappa{\bar C}}{1-r}>e_0$, the set $\cK$ is infinite;
{\em (ii)} if $\frac{\kappa{\bar C}}{1-r}\geq e_0$, then $e_K\leq\frac{\kappa{\bar C}r^K}{1-r}$ for every $K\geq0$.
\end{lemma}

\begin{proof}
Summing {\em (a)} over $k=i,\dots,K-1$ shows that, whenever $k\notin\cK$ for $i\leq k\leq K-1$,
\begin{equation}\label{eq:env-sum}
e_K\leq e_i-\frac{\kappa{\bar C}(r^i-r^K)}{1-r},
\end{equation}
with strict inequality if $i<K$.

{\em (i)} Suppose that $\cK$ is finite. If $\cK=\emptyset$, then \eqref{eq:env-sum} with $i=0$ and $e_K\geq0$ give $\frac{\kappa{\bar C}(1-r^K)}{1-r}<e_0$ for all $K\geq1$, and letting $K\to+\infty$ contradicts $\frac{\kappa{\bar C}}{1-r}>e_0$. Otherwise, let $m=1+\max\cK$; then \eqref{eq:env-sum} with $i=m$ and {\em (b)} give $\frac{\kappa{\bar C}(r^m-r^K)}{1-r}\leq e_m\leq\omega{\bar C}r^{m-1}$ for all $K>m$, and letting $K\to+\infty$ yields $\frac{\kappa r}{1-r}\leq\omega$, a contradiction.

{\em (ii)} Fix $K\geq0$. If $K\in\cK$, then $e_K\leq{\bar C}r^K\leq\frac{\kappa{\bar C}r^K}{1-r}$, since $\frac{\kappa}{1-r}>\frac{\omega}{r}\geq1$. Otherwise, if no $k<K$ belongs to $\cK$, then \eqref{eq:env-sum} with $i=0$ gives $e_K\leq e_0-\frac{\kappa{\bar C}}{1-r}+\frac{\kappa{\bar C}r^K}{1-r}\leq\frac{\kappa{\bar C}r^K}{1-r}$. Finally, if $j\eqdef\max\{k<K \,|\, k\in\cK\}$ exists, then $k\notin\cK$ for $j<k\leq K$, and \eqref{eq:env-sum} with $i=j+1$ together with {\em (b)} gives
$$
e_K\leq\omega{\bar C}r^j-\frac{\kappa{\bar C}(r^{j+1}-r^K)}{1-r}={\bar C}r^j\left(\omega-\frac{\kappa r}{1-r}\right)+\frac{\kappa{\bar C}r^K}{1-r}<\frac{\kappa{\bar C}r^K}{1-r}.
$$
\end{proof}

\subsection{Proof of Theorem \ref{th:strongly-convex-prac}}

{\em (i)} For any $k \notin \cK$, $\|x_k - x_\star\| > Cq^k$. Writing $d_k \eqdef \nabla f(x_k) - \nabla f(x_\star)$, we have
\begin{eqnarray*}
	\langle d_k, x_{k+1} - x_\star \rangle &=& \langle d_k, x_k - x_\star \rangle  + \langle d_k, x_{k+1} - x_k \rangle 
	\ \geq\ \langle d_k, x_k - x_\star \rangle  - t_k \|d_k\| \\ 
	&\geq& \theta \|d_k\| \|x_k - x_\star\| -  t_k \|d_k\| 
	\ >\ \|d_k\| (\theta Cq^k - t_k) \ =\ 0, 
\end{eqnarray*}
where we use the Cauchy–Schwarz inequality and $\|x_{k+1}-x_k\| \leq t_k$ in the first inequality, and use Lemma C.3 of \citet{richtarik2026local} in the second inequality. From Theorem \ref{th:Tk} {\em (iii)}, we know $t_k \in \cT_k$. Then by Theorem \ref{th:another-ad} {\em (iii)}, we have 
$$
\|x_{k+1} - x_\star\|^2 < \|x_k - x_\star\|^2 - \theta^2 C^2 q^{2k}. 
$$

{\em (ii)} For $k\geq 1$ such that $k-1 \in \cK$, we have 
\begin{eqnarray*}
	\|x_k - x_\star\|^2 &=& \|x_k - x_{k-1} + x_{k-1} - x_\star\|^2\\ 
	&\leq& (1+a)\|x_k - x_{k-1}\|^2 + (1+\frac{1}{a})\|x_{k-1} - x_\star\|^2 \\ 
	&\leq& (1+a) t_{k-1}^2 + (1+\frac{1}{a})C^2q^{2k-2} \\ 
	&=& (1+a) \theta^2C^2 q^{2k-2} + (1+\frac{1}{a})C^2q^{2k-2} \\ 
	&=& (1+2\theta+\theta^2) C^2 q^{2k-2}, 
\end{eqnarray*}
where we use Young's inequality in the first inequality, use $\|x_k-x_{k-1}\|\leq t_{k-1}$ and $k-1 \in \cK$ in the second inequality, and use $a = \frac{1}{\theta}$ in the last equality.

{\em (iii)--(iv)} By (i) and (ii), Lemma~\ref{lm:envelope} applies with $e_k = \|x_k - x_\star\|^2$, ${\bar C} = C^2$, $r = q^2$, $\kappa = \theta^2$ and $\omega = (1+\theta)^2$. Indeed, $\frac{\kappa{\bar C}}{1-r} = \frac{C^2}{1-\alpha}$, and $\frac{\kappa r}{1-r} = \frac{q^2}{1-\alpha} = \frac{1}{1-\alpha} - \theta^2 > (1+\theta)^2$ by the assumption $1-\alpha < \frac{1}{1+2\theta+2\theta^2}$. Hence (iii) follows from Lemma~\ref{lm:envelope}(i), and Lemma~\ref{lm:envelope}(ii) gives $\|x_K - x_\star\|^2 \leq \frac{C^2 q^{2K}}{1-\alpha} \leq \epsilon^2$ whenever $K \geq \frac{1}{1-q^2} \ln \frac{\theta^2 C^2}{\epsilon^2(1-q^2)}$, since $\frac{\theta^2}{1-q^2} = \frac{1}{1-\alpha}$ and $\ln \frac{1}{q^2} \geq 1-q^2$.

For $\alpha = \frac{2\theta + 3\theta^2}{1+2\theta+3\theta^2}$ and $C = \sqrt{1-\alpha} \|x_0-x_\star\|$, we have 
$$
\frac{1}{1-q^2} = \frac{1}{(1-\alpha)\theta^2} = \frac{1+2\theta+3\theta^2}{\theta^2} = \frac{1}{\theta^2} + \frac{2}{\theta} + 3, 
$$
and 
$$
\frac{\theta^2 C^2}{\epsilon^2(1-q^2)} = \frac{C^2}{\epsilon^2 (1-\alpha)} =  \frac{(1-\alpha)\|x_0-x_\star\|^2}{\epsilon^2 (1-\alpha)} = \frac{\|x_0-x_\star\|^2}{\epsilon^2},  
$$
which completes the proof.

\subsection{Proof of Theorem \ref{th:sc-halving}}

We use the notation of Appendix~\ref{sec:proof-halving}, and let $E_j \eqdef \frac{L_c}{2}\left( 1 + \frac{L_c}{\mu} \right)\rho_j^2$, $j \geq 0$. By Lemma \ref{lm:window} {\em (ii)} applied to the halving window with radius $\rho_j$, $\Phi_{k_j} \leq E_j$. Moreover, $\Phi_k \leq E_j$ for every $k \geq k_j$: the steps from $x_k$ with $k \geq k_j$ use radii $\rho_i$ with $i > j$, a short step decreases $F$ by Lemma \ref{lm:short-long} {\em (i)}, and a long step from $x_k$ gives $\Phi_{k+1} \leq E_i \leq E_j$ by Lemma \ref{lm:short-long} {\em (ii)}. As in the proof of Theorem \ref{th:convex-halving}, the number of surviving windows with radius $\rho_j$ is at most $\nicefrac{\Phi_0}{cTt_0^2}$ for $j = 0$ and $\nicefrac{E_{j-1}}{cT\rho_j^2} = \nicefrac{2L_c(1 + L_c/\mu)}{cT}$ for $j \geq 1$, so that $k_0 \leq N_0$ and $k_j - k_{j-1} \leq N$ for $j \geq 1$. Hence $K \geq N_0 + JN$ implies $K \geq k_J$ and
$$
\Phi_K \leq E_J = \frac{L_c t_0^2}{2}\left( 1 + \frac{L_c}{\mu} \right)4^{-J}.
$$
Finally, $F$ is $\mu$-strongly convex and $0 \in \partial F(x_\star)$, so that $\frac{\mu}{2}\|x_K - x_\star\|^2 \leq \Phi_K$.

\section{Proofs for the non-convex case}\label{sec:proof-nonconvex}

We first establish a technical lemma, which describes a single step of BGD with an arbitrary radius. Parts (iii) and (iv) are similar to Theorem D.16 of \citet{richtarik2026local}. 

\begin{lemma}\label{lm:non-convex-step}
	Let Assumption \ref{as:basic} hold. Let $f$ be $L$-smooth on an open set containing $\dom  g$. Fix some $\gamma \in (0, \nicefrac{1}{L}]$ and let \( x_{k+1} \) be obtained from \( x_k \) after one step of Broximal GD with radius \( t_k>0 \), i.e.,
	\[
	x_{k+1} \in \argmin_{z \in \mathcal{B}(x_k, t_k)} m_k(z) \eqdef g(z) + \langle \nabla f(x_k), z \rangle.
	\]
	Denote $D_k \eqdef \cD_g(x_k, \gamma)$ and $r_k \eqdef \gamma \|G_\gamma(x_k)\|$. Then the following results hold. 
	\begin{itemize}
		\item[(i)] $D_k \geq \|G_\gamma(x_k)\|^2$. 
		\item[(ii)] $F(x_{k+1}) \leq F(x_k) - \min\left\{ t_k \sqrt{D_k}, \frac{\gamma}{2} D_k \right\} + \frac{1}{2\gamma} t_k^2$; in particular, $F(x_{k+1}) \leq F(x_k) + \frac{1}{2\gamma} t_k^2$. 
		\item[(iii)] If $t_k < r_k$, then $F(x_{k+1}) < F(x_k) - \frac{1}{2\gamma} t_k^2$. 
		\item[(iv)] If $t_k = r_k$, then $F(x_{k+1}) \leq F(x_k) - \frac{\gamma}{2} D_k \leq F(x_k) - \frac{\gamma}{2}\|G_\gamma(x_k)\|^2$. 
	\end{itemize}
\end{lemma}

\begin{proof}
	Denote $y_k \eqdef \operatorname{Prox}_{\gamma g} \left( x_k - \gamma \nabla f(x_k) \right)$ and $\hat m_k(z) \eqdef m_k(z) - m_k(x_k) = \langle \nabla f(x_k), z - x_k \rangle + g(z) - g(x_k)$. By the definition of $G_\gamma$, $\|y_k - x_k\| = \gamma\|G_\gamma(x_k)\| = r_k$. Since $y_k$ is the minimizer in \eqref{eq:Dg} with $x = x_k$, we have 
	\begin{equation}\label{eq:hk-yk}
		\hat m_k(y_k) + \frac{1}{2\gamma} r_k^2 = -\frac{\gamma}{2} D_k. 
	\end{equation}
	
	{\em (i)} Choosing $z = x_k$ and $y = y_k$ in Lemma \ref{lm:prox}, we obtain $\hat m_k(y_k) = m_k(y_k) - m_k(x_k) \leq -\frac{1}{\gamma} r_k^2$. Combining this with \eqref{eq:hk-yk} yields $-\frac{\gamma}{2} D_k \leq -\frac{1}{2\gamma} r_k^2$, i.e., 
	\begin{equation}\label{eq:Dg-G}
		D_k \geq \frac{r_k^2}{\gamma^2} = \|G_\gamma(x_k)\|^2. 
	\end{equation}
	
	Next, since $\|x_{k+1} - x_k\| \leq t_k$, from the $L$-smoothness of $f$ and $L \leq \frac{1}{\gamma}$, we have 
	\begin{eqnarray}
		F(x_{k+1}) &\leq& g(x_{k+1}) + f(x_k) + \langle \nabla f(x_k), x_{k+1}-x_k \rangle + \frac{L}{2} \|x_{k+1} - x_k\|^2 \notag \\ 
		&\leq& F(x_k) + \hat m_k(x_{k+1}) + \frac{1}{2\gamma} t_k^2. \label{eq:Fxk+1-hk}
	\end{eqnarray}
	We bound $\hat m_k(x_{k+1})$ in two cases. If $t_k \geq r_k$, then $y_k \in \cB(x_k, t_k)$, and from the definition of $x_{k+1}$ and \eqref{eq:hk-yk},
	\begin{equation}\label{eq:case-large}
		\hat m_k(x_{k+1}) \leq \hat m_k(y_k) = -\frac{\gamma}{2} D_k - \frac{1}{2\gamma} r_k^2. 
	\end{equation}
	If $t_k < r_k$, then ${\tilde y}_k \eqdef x_k + \frac{t_k}{r_k}(y_k - x_k) \in \cB(x_k, t_k)$, and from the convexity of $\hat m_k(\cdot)$, $\hat m_k(x_k) = 0$ and \eqref{eq:hk-yk}, 
	\begin{equation}\label{eq:case-small}
		\hat m_k(x_{k+1}) \leq \hat m_k({\tilde y}_k) \leq \frac{t_k}{r_k} \hat m_k(y_k) = -\frac{t_k \gamma D_k}{2 r_k} - \frac{t_k r_k}{2\gamma}. 
	\end{equation}
	
	{\em (ii)} In the case $t_k \geq r_k$, \eqref{eq:case-large} gives $\hat m_k(x_{k+1}) \leq -\frac{\gamma}{2} D_k$. In the case $t_k < r_k$, the inequality $a + b \geq 2\sqrt{ab}$ applied to \eqref{eq:case-small} gives $\hat m_k(x_{k+1}) \leq -t_k\sqrt{D_k}$. In both cases $\hat m_k(x_{k+1}) \leq - \min\{ t_k \sqrt{D_k}, \frac{\gamma}{2} D_k \}$, and the claim follows from \eqref{eq:Fxk+1-hk}. 
	
	{\em (iii)} From \eqref{eq:Dg-G}, $\gamma^2 D_k \geq r_k^2$, thus $\frac{t_k \gamma D_k}{2 r_k} \geq \frac{t_k r_k}{2\gamma} > \frac{t_k^2}{2\gamma}$, where we use $t_k < r_k$. Then \eqref{eq:case-small} gives $\hat m_k(x_{k+1}) < -\frac{1}{\gamma} t_k^2$, and the claim follows from \eqref{eq:Fxk+1-hk}. 
	
	{\em (iv)} Plugging \eqref{eq:case-large} with $t_k = r_k$ into \eqref{eq:Fxk+1-hk} gives $F(x_{k+1}) \leq F(x_k) - \frac{\gamma}{2} D_k$, and the second inequality follows from (i). 
\end{proof}

\subsection{Proof of Theorem \ref{th:non-convex-PL-theory}}

For $t_k = \gamma  \|G_{\gamma} (x_k)\|$, Lemma \ref{lm:non-convex-step}(iv) and Assumption \ref{as:proximalPL} (recall $x_k \in \dom g$) give 
$$
F(x_{k+1}) - F_\star \leq F(x_k) - F_\star - \frac{\gamma}{2}  \cD_g(x_k, \gamma) \leq (1-\gamma \mu) (F(x_k) - F_\star) 
$$
for all $k \geq 0$ with $G_\gamma(x_k) \neq 0$. If $G_\gamma(x_k) = 0$, then $t_k = 0$, $x_{k+1} = x_k$, and the minimizer in \eqref{eq:Dg} is $y = x_k$, so that $\cD_g(x_k, \gamma) = 0$ and Assumption \ref{as:proximalPL} gives $F(x_k) = F_\star$; the inequality then holds trivially. Iterating this inequality yields the result.

\subsection{Proof of Theorem \ref{th:non-convex-PL}}

We exclude the trivial case $F \equiv F_\star$ on $\dom g$. Then $\gamma\mu \leq 1$: taking any $x_0 \in \dom g$ with $F(x_0) > F_\star$, Theorem~\ref{th:non-convex-PL-theory} gives $0 \leq F(x_1) - F_\star \leq (1-\gamma\mu)(F(x_0) - F_\star)$. Hence $q \in (0,1)$.

{\em (i)} For $k\notin \cK$, Assumption \ref{as:proximalPL} gives $D_k \eqdef \cD_g(x_k, \gamma) \geq 2\mu (F(x_k) - F_\star) > 2\mu Cq^k$. Hence 
$$
t_k \sqrt{D_k} > \gamma \sqrt{\frac{\mu Cq^k}{2}} \cdot \sqrt{2\mu Cq^k} = \gamma\mu Cq^k \qquad \text{and} \qquad \frac{\gamma}{2} D_k > \gamma\mu Cq^k. 
$$
Then Lemma \ref{lm:non-convex-step}(ii) gives 
$$
F(x_{k+1}) < F(x_k) - \gamma\mu Cq^k + \frac{1}{2\gamma} t_k^2 = F(x_k) - \gamma\mu Cq^k + \frac{1}{4}\gamma\mu Cq^k = F(x_k) - \frac{3}{4}\gamma\mu Cq^k. 
$$

{\em (ii)} By Lemma \ref{lm:non-convex-step}(ii), $F(x_k) \leq F(x_{k-1}) + \frac{1}{2\gamma}t_{k-1}^2 = F(x_{k-1}) + \frac{\gamma\mu}{4}Cq^{k-1}$. Since $k-1 \in \cK$, we obtain $F(x_k) - F_\star \leq \left(1+\frac{\gamma\mu}{4}\right) Cq^{k-1}$. 

{\em (iii)--(iv)} By (i) and (ii), Lemma~\ref{lm:envelope} applies with $e_k = F(x_k) - F_\star$, ${\bar C} = C$, $r = q$, $\kappa = \frac{3\gamma\mu}{4}$ and $\omega = 1+\frac{\gamma\mu}{4}$. Indeed, $\frac{\kappa{\bar C}}{1-r} = \frac{C}{1-\alpha}$, and $\frac{\kappa r}{1-r} = \frac{q}{1-\alpha} = \frac{1}{1-\alpha} - \frac{3\gamma\mu}{4} > 1+\frac{\gamma\mu}{4}$ by the assumption $1-\alpha < \frac{1}{1+\gamma\mu}$. Hence (iii) follows from Lemma~\ref{lm:envelope}(i), and Lemma~\ref{lm:envelope}(ii) gives $F(x_K) - F_\star \leq \frac{Cq^K}{1-\alpha} \leq \epsilon$ whenever $K \geq \frac{1}{1-q} \ln \frac{3\gamma \mu C}{4\epsilon(1-q)}$, since $\frac{3\gamma\mu}{4(1-q)} = \frac{1}{1-\alpha}$ and $\ln\frac{1}{q} \geq 1-q$.

Finally, for $\alpha = \frac{3\gamma\mu}{2+3\gamma \mu}$, we have $1-\alpha = \frac{2}{2+3\gamma\mu} < \frac{1}{1+\gamma\mu}$ and $\frac{1}{1-q} = \frac{4}{3(1-\alpha)\gamma\mu} = \frac{2(2+3\gamma\mu)}{3\gamma\mu} = \frac{4}{3\gamma\mu} + 2$, and for $C = (1-\alpha) (F(x_0) - F_\star)$, we have $\frac{3\gamma\mu C}{4\epsilon(1-q)} = \frac{C}{\epsilon (1-\alpha)} = \frac{F(x_0) - F_\star}{\epsilon}$, which completes the proof.

\subsection{Proof of Theorem \ref{th:pl-halving}}

We use the notation of Appendix~\ref{sec:proof-halving}, with $F_\star = \inf F$, and write $D_k \eqdef \cD_g(x_k, \gamma)$. By Lemma \ref{lm:non-convex-step} {\em (ii)}, a step from $x_k$ with radius $\tau$ satisfies
\begin{equation}\label{eq:pl-halving-step}
	F(x_{k+1}) \leq F(x_k) - \min\left\{ \tau\sqrt{D_k}, \frac{\gamma}{2}D_k \right\} + \frac{\tau^2}{2\gamma}.
\end{equation}
We call the step \emph{small} if $\min\{ \tau\sqrt{D_k}, \frac{\gamma}{2}D_k \} < s\tau^2$. A step that is not small decreases $F$ by at least $s\tau^2 - \frac{\tau^2}{2\gamma} = c\tau^2$, by \eqref{eq:pl-halving-step}. A small step satisfies $\sqrt{D_k} < s\tau$ or $D_k < \frac{2s}{\gamma}\tau^2$, i.e., $D_k < \kappa\tau^2$, so that Assumption \ref{as:proximalPL} and \eqref{eq:pl-halving-step} give
\begin{equation}\label{eq:pl-halving-small}
	F(x_{k+1}) - F_\star \leq F(x_k) - F_\star + \frac{\tau^2}{2\gamma} \leq \frac{D_k}{2\mu} + \frac{\tau^2}{2\gamma} < \left( \frac{\kappa}{2\mu} + \frac{1}{2\gamma} \right)\tau^2.
\end{equation}

Let $E_j \eqdef \left( \frac{\kappa}{2\mu} + \frac{1}{2\gamma} \right)\rho_j^2$, $j \geq 0$. If a window with radius $\rho_j$ contained no small step, it would decrease $F$ by at least $cT\rho_j^2$ and would not be halving. Hence the halving window with radius $\rho_j$ contains a small step; the steps after the last one are not small and decrease $F$, so that $\Phi_{k_j} \leq E_j$ by \eqref{eq:pl-halving-small}. In the same way, $\Phi_k \leq E_j$ for every $k \geq k_j$: the steps from $x_k$ with $k \geq k_j$ use radii $\rho_i$ with $i > j$, a step that is not small decreases $F$, and a small step from $x_k$ gives $\Phi_{k+1} \leq E_i \leq E_j$ by \eqref{eq:pl-halving-small}.

As in the proof of Theorem \ref{th:convex-halving}, the number of surviving windows with radius $\rho_j$ is at most $\nicefrac{\Phi_0}{cTt_0^2}$ for $j = 0$ and $\nicefrac{E_{j-1}}{cT\rho_j^2} = \nicefrac{4E_{j-1}}{cT\rho_{j-1}^2}$ for $j \geq 1$; adding the halving window, we obtain $k_0 \leq N_0$ and $k_j - k_{j-1} \leq T + \frac{4}{c}\left( \frac{\kappa}{2\mu} + \frac{1}{2\gamma} \right) = N$ for $j \geq 1$. Hence $K \geq N_0 + JN$ implies $K \geq k_J$ and $F(x_K) - F_\star \leq E_J = \left( \frac{\kappa}{2\mu} + \frac{1}{2\gamma} \right)t_0^2 4^{-J}$.

\subsection{Proof of Theorem \ref{th:non-convex-withoutPL}}

{\em (i)} If $k\notin \cK$, then $\|G_{\gamma} (x_k)\| >a_k$, which implies that $t_k = \gamma a_k < \gamma \|G_{\gamma} (x_k)\|$. Since $\cT_k(y_k) = (0, \gamma \|G_{\gamma} (x_k)\| )$ from Theorem \ref{th:recoverPGD}, we know $t_k \in \cT_k(y_k)$. Furthermore, from Lemma \ref{lm:non-convex-step}(iii), we have 
$$
F(x_{k+1}) < F(x_k) - \frac{1}{2\gamma} t_k^2 = F(x_k) - \frac{\gamma}{2} a_k^2. 
$$

{\em (ii)} By {\em (i)}, Lemma \ref{lm:dichotomy}(i) applies with $e_k = F(x_k) - F_\star$ and $\delta_k = \frac{\gamma}{2}a_k^2$, so $\cK$ is infinite.

{\em (iii)} For $a_k =  \|G_{\gamma} (x_k)\|$, we have $t_k = \gamma  \|G_{\gamma} (x_k)\|$. If $t_k > 0$, Lemma \ref{lm:non-convex-step}(iv) applies; if $t_k = 0$, then $x_{k+1} = x_k$ and $G_\gamma(x_k) = 0$. In both cases, we can obtain 
$$
F(x_{k+1}) - F_\star \leq F(x_k) - F_\star - \frac{\gamma}{2}  \|G_{\gamma} (x_k)\|^2, 
$$
for all $k\geq 0$, and Lemma \ref{lm:dichotomy}(ii) with $e_k = F(x_k) - F_\star$ and $\delta_k = \frac{\gamma}{2}\|G_\gamma(x_k)\|^2$ yields the claim.

\subsection{Proof of Corollary \ref{co:non-convex-withoutPL}}

If $\min_{0 \leq k \leq K-1} \|G_{\gamma} (x_k)\| > a$ for some $K \geq \frac{2(F(x_0) - F_\star)}{\gamma a^2}$, then $k \notin \cK$ for all $0 \leq k \leq K-1$, and summing the inequality of Theorem \ref{th:non-convex-withoutPL} {\em (i)} over these $k$ gives $0\leq F(x_K) - F_\star < F(x_0) - F_\star - \frac{K\gamma a^2}{2}$, i.e., $K < \frac{2(F(x_0) - F_\star)}{\gamma a^2}$, a contradiction. The second statement is the case $\gamma = \nicefrac{1}{L}$ and $a = \epsilon$, for which $t_k = \gamma a = \nicefrac{\epsilon}{L}$.

\subsection{Proof of Theorem \ref{th:nc-halving}}

We use the notation of Appendix~\ref{sec:proof-halving}, with $F_\star = \inf F$. Lemmas \ref{lm:short-long} and \ref{lm:window} do not require the convexity of $f$, except for their last statements, which are not used here. Surviving windows decrease $F$, and by Lemma \ref{lm:window} {\em (i)} the halving window with radius $\rho_j$, which is unique for each $j$, increases $F$ by at most $\frac{T}{2}L_c\rho_j^2$; hence $F(x_{iT}) \leq F(x_0) + \frac{T}{2}L_ct_0^2\sum_{j \geq 0}4^{-j} = F(x_0) + \frac{2}{3}TL_ct_0^2$ at every checkpoint, and the same telescoping shows that the total decrease over all surviving windows is at most $\bar\Phi$. Let $n_j$ be the number of surviving windows with radius $\rho_j$. Since each of them decreases $F$ by at least $cT\rho_j^2 \geq cT\rho_J^2$ for $j \leq J$, we have $\sum_{j=0}^{J}n_j \leq \nicefrac{\bar\Phi}{cT\rho_J^2}$, and therefore
$$
k_J = T\sum_{j=0}^{J}(n_j + 1) \leq (J+1)T + \frac{\bar\Phi}{c\,\rho_J^2} = (J+1)T + \frac{4^{J}\,\bar\Phi}{c\,t_0^2}.
$$
By Lemma \ref{lm:window} {\em (ii)}, the halving window with radius $\rho_J$ contains an iterate $x_l$, $l < k_J$, with $\|G_{\nicefrac{1}{L_c}}(x_l)\| < L_c\rho_J$, and $\|G_{\nicefrac{1}{L}}(x_l)\| \leq \|G_{\nicefrac{1}{L_c}}(x_l)\|$ by Lemma \ref{lm:prox-mono}, since $\nicefrac{1}{L_c} \leq \nicefrac{1}{L}$. This proves the first claim. For the second one, take $J = J_\epsilon \geq 0$; then $L_c\rho_J \leq \epsilon$ and $4^{J} < \nicefrac{4L_c^2t_0^2}{\epsilon^2}$.

\end{document}